\documentclass[Afour,sageh,times]{sagej}

\usepackage{lineno}
\usepackage{amsmath}
\usepackage{amssymb}
\usepackage{xcolor}
\usepackage{caption}
\usepackage{subcaption}
\usepackage{rotating}
\usepackage{listings}
\usepackage{algorithm}
\usepackage[noend]{algpseudocode}

\graphicspath{{pics/}}

\newcommand\trd[1]{\textcolor{red}{#1}} 

\newcommand{\divr}{\nabla\cdot}
\newcommand{\qb}{\mathbf q}
\newcommand{\Lm}{\mathcal L}
\newcommand{\Fm}{\mathcal F}
\newcommand{\Sm}{\mathcal S}
\newcommand{\Rm}{\mathcal R}
\newcommand{\Nm}{\mathcal N}

\newcommand{\dm}{tt} 

\newcommand{\qtt}{\mathbf q^{\dm}}
\newcommand{\utt}{\mathbf u^{\dm}}
\newcommand{\ptt}{P^{\dm}}
\newcommand{\rtt}{\rho^{\dm}}
\newcommand{\ttt}{\theta^{\dm}}

\newcommand{\qtte}{\mathbf q^{\dm}_{e}}
\newcommand{\utte}{\mathbf u^{\dm}_{e}}
\newcommand{\ptte}{P^{\dm}_{e}}
\newcommand{\rtte}{\rho^{\dm}_{e}}
\newcommand{\ttte}{\theta^{\dm}_{e}}

\newcommand{\Utt}{\mathbf U^{\dm}}
\newcommand{\Ttt}{\Theta^{\dm}}

\newcommand{\Utte}{\mathbf U^{\dm}_{e}}
\newcommand{\Ttte}{\Theta^{\dm}_{e}}

\newcommand{\gb}{\mathbf g}

\newcommand{\nsr}{\emph{standard} }
\newcommand{\sr}{\emph{schur} }
\newcommand{\comment}[1] {}

\def\volumeyear{2017}
\begin{document}

\runninghead{Abdi, Giraldo, Constantinescu, Carr III, Wilcox and Warburton}

\title{Acceleration of the Implicit-Explicit Non-hydrostatic Unified Model of the Atmosphere (NUMA) on Manycore Processors}

\author{Daniel S. Abdi\affilnum{1} and Francis X. Giraldo\affilnum{1} and Emil M. Constantinescu\affilnum{3} and Lester E. Carr III\affilnum{1}and Lucas C. Wilcox\affilnum{1} and Timothy C. Warburton \affilnum{2} }

\affiliation{\\
\affilnum{1}Department of Applied Mathematics, Naval Postgraduate School, USA\\
\affilnum{2}Department of  Mathematics, Virginia Tech University, USA\\
\affilnum{3}Argonne National Laboratory, USA}

\corrauth{Daniel S. Abdi, Naval Postgraduate School
Monterey, CA 93943, USA.}

\email{dsabdi@nps.edu}

\begin{abstract}
We present the acceleration of an IMplicit-EXplicit (IMEX) non-hydrostatic atmospheric model on manycore processors such as GPUs and Intel's MIC architecture. IMEX time integration methods sidestep the constraint imposed by the Courant-Friedrichs-Lewy condition on explicit methods through corrective implicit solves within each time step. In this work, we implement and evaluate the performance of IMEX on manycore processors relative to explicit methods. Using 3D-IMEX at Courant number C=15 , we obtained a speedup of about 4X relative to an explicit time stepping method run with the maximum allowable C=1. Moreover, the unconditional stability of IMEX with respect to the fast waves means the speedup can increase significantly with the Courant number as long as the accuracy of the resulting solution is acceptable. We show a speedup of 100X at C=150 using 1D-IMEX to demonstrate this point. Several improvements on the IMEX procedure were necessary in order to outperform our results with explicit methods: a) reducing the number of degrees of freedom of the IMEX formulation by forming the Schur complement; b) formulating a horizontally-explicit vertically-implicit (HEVI) 1D-IMEX scheme that has a lower workload and potentially better scalability than 3D-IMEX; c) using high-order polynomial preconditioners to reduce the condition number of the resulting system; d) using a direct solver for the 1D-IMEX method by performing and storing LU factorizations once to obtain a constant cost for any Courant number. Without \emph{all} of these improvements, explicit time integration methods turned out to be difficult to beat. We discuss in detail the IMEX infrastructure required for formulating and implementing efficient methods on manycore processors. Several parametric studies are conducted to demonstrate the gain from each of the above mentioned improvements. Finally, we validate our results with standard benchmark problems in numerical weather prediction and evaluate the performance and scalability of the IMEX method using up to 4192 GPUs and 16 Knights Landing processors.
\end{abstract}


\maketitle

\section{Introduction}


The multiscale dynamics in the atmosphere supports different types of wave motion. The slow processes are best captured with explicit time stepping methods; however, the fast waves present in the atmosphere, such as acoustic and gravity waves, place severe limitations on the maximum allowable time step that can be taken by explicit methods. The Courant number is given by
\begin{equation}
C=\frac{c_{max} \Delta t}{\Delta x}
\end{equation}
where $c_{max}$ is the maximum wave speed in the system, $\Delta x$ is the smallest grid spacing and $\Delta t$ is the  time step. Most explicit time stepping methods require $C \le 1$.  To deal with this restriction,
various approaches have been developed such as a) reducing $c_{max}$ by filtering fast-moving waves. For example, hydrostatic models eliminate all vertical motion, however, the validity of these models reach their limit at about 10 km resolution. In this work we use the Non-hydrostatic Unified Model of the Atmosphere (NUMA) that has been designed from the start to work with resolutions in the non-hydrostatic regime. Other approaches for filtering the acoustic waves from the equations include the, e.g., Boussinesq approximations,  anelastic and pseudo-incompressible systems (see e.g. \citep{durran1989}).  Other approaches for increasing the maximum allowable time-step include b) using larger grid spacing $\Delta x$ for the fast-moving waves e.g. the \citet{turkel1979} method,  and c) treating the fast waves implicitly thereby making the method unconditionally stable with regard to the fast-moving waves. In this paper, we consider this last approach via IMEX time-integrators.

The dynamics of the atmosphere supports waves of different temporal scales such as acoustic waves, gravity waves and Rossby waves. Acoustic waves are the fastest waves in the atmosphere but they have little to no effect on the large-scale processes in the linear regime. As a result, it is inefficient to  use explicit methods because one would be forced to use small time steps due to a physical phenomenon that is essentially inconsequential. Gravity waves are the next fastest moving waves in the atmosphere, however, the energy carried by these waves is a very small percentage of the total energy. This is an important fact that determines the accuracy of methods for the solution of `stiff' equations that have fast decaying modes which will inevitably contain little energy \cite{fulton2004}. The Rossby waves (contained in the advection terms) are the next fastest waves; these waves are often treated using explicit time-stepping methods. 

In this work, we use semi-implicit methods to ameliorate the stringent time-step restriction imposed by the CFL condition, however, here we shall refer to them as IMEX methods. The fast acoustic and gravity waves are treated implicitly while the Rossby waves are treated explicitly; thus, the maximum stable time step is governed by the relatively slow Rossby waves.  Semi-implicit (SI) schemes have been widely used in numerical weather predication (NWP) starting with \cite{kwizak1971,robert1972}. Their use have been credited with a six-fold increase in computational efficiency of modern weather prediction systems \cite{restelli:thesis}. In \cite{giraldo:2006b}, a hybrid Eulerian-Lagrangian semi-implicit method is proposed to absorb the advective terms in the Lagrangian (material) derivative, and thus ensure a virtual disappearance of the CFL condition. Other operator-splitting methods include the split-explicit method \cite{gadd1978,klemp2007} used in NWP, and projection and fractional-step methods commonly used for solving the incompressible Navier-Stokes equations. The split-explicit methods \comment{also known as horizontally-explicit vertically-implicit (HEVI), }do a directional splitting in which both fast and slow waves in the horizontal direction are advanced explicitly using different time steps with possible sub-cycling for the fast waves. The method has a potential for good performance on massively parallel systems of manycore and multicore architectures \cite{norman2013,michalakesManish2008}.

One could also treat all waves implicitly, however, fully implicit methods require the solution of non-linear systems of equations which is computationally expensive, difficult to implement and often poorly scalable. In \cite{chao2014}, a scalable fully implicit non-hydrostatic atmospheric model is discussed. The equations are solved using a Jacobian-free Newton-Krylov (JFNK) algorithm. We should note that in the IMEX approach used in this work, the implicitly discretized terms are linear thereby resulting in a constant Jacobian matrix; this avoids the need for a non-linear solver. In \citep{Yang2016}, an ultra-scalable (upto 10.5 million cores) fully-implicit non-hydrostatic atmospheric model is presented to solve a global NWP problem with a record  horizontal resolution of 488 m. The solver uses a novel hybrid domain-decomposed multi-grid preconditioner to accelerate convergence of the iterative solution, and a parallel ILU preconditioner suited for many-core processors.  In \cite{archibald2015}, acceleration of fully implicit time stepping methods for solving the shallow water equations on the GPU is described. The community atmospheric model (CAM), which is a spectral element (SE) code like our model, is used for that study. The implicit solver in CAM-SE that makes calls to the Trilinos library was replaced with GPU kernels to accelerate the residual calculation in the JFNK solver by $\leq$ 3X.

Implementation of IMEX methods on manycore processors is far more challenging than that of explicit methods, hence, there is not a lot of literature on the subject especially regarding its use in numerical weather prediction. IMEX methods require substantial infrastructure that is difficult to implement on manycore processors  with the possibility that it might not perform better than explicit methods. There is also the issue of additional accuracy obtained from using small time steps that further motivates the use of explicit methods. We have already discussed the implementation of a scalable explicit solver in our previous work \cite{abdi2016}, which we will use as a base-line for comparison with our IMEX solvers. Other implementations on GPUs that need to be mentioned include the work in \cite{yang2013} whereby a highly scalable hybrid CPU-GPU algorithm for solving the shallow water equations using explicit time stepping is discussed.  A partitioning scheme  and communication strategy for optimal load balancing of work between the CPU and GPU is proposed for a simulation on a cubed sphere grid. With regard to scalable IMEX time stepping, the work of \citet{Muller2015} on a GPU implementation of an anisotropic elliptic PDE solver for atmospheric modeling is a significant step. They use two scalable solvers, a matrix-free conjugate gradient method and a geometric multi-grid method, to solve the resulting linear system of equations. They showed good scalability using up to 16384 GPUs of the Titan supercomputer; a peak bandwidth usage of 40\% and peak floating point operations rate of 3\% is reported that suggests the memory-bound nature of their code. In \cite{shi2012}, a GPU implementation of matrix-free iterative solvers for the solution of a 3D Helmholtz problem in NWP is discussed along with the infrastructure required for preconditioning.  In \cite{you2015}, a GPU accelerated semi-implicit alternating direction implicit (ADI) method is discussed for solving the incompressible and compressible Navier-Stokes equations. They mention that inversion of multiple tridiagonal matrices on the GPU is the major challenge to getting good performance. In \cite{ong2015}, parallel semi-implicit time integrators are implemented on an NVIDIA GPU using the cuBLAS library. They observed that a fourth order solution using an IMEX time integrator on 4 GPUs  took the same amount of time as a first order solution on a single GPU. Clearly, it is non-trivial to develop efficient IMEX time-integrators for manycore processors.  Let us now describe our IMEX time-integration approach which has yielded a substantial improvement in speedup over our explicit time-integrators.

\section{IMEX time integration}
The time-integration statement can be written in compact  form
\begin{equation}
\label{timeint}
\frac{\partial \qb}{\partial t} = \Rm(\qb)
\end{equation}
where $\Rm(\qb)$ represents all the terms aside from the temporal derivative, $\qb$ contains the unknown values of the state variables at the current time step. Using the method of lines, we discretize the right hand side operator in space and convert Eq.\ \eqref{timeint} to a system of ordinary differential equations, for which several time stepping algorithms are available.

To derive the IMEX method, we seek a linear operator $\Lm(\qb)$ that contains the problematic fast waves and formulate the splitting for the right-hand side as
\begin{equation}
\label{sieqn}
\begin{split}
\Rm(\qb) = \underbrace{ \Rm(\qb^*) - \Lm(\qb^*) }_\text{explicit} + \underbrace{\Lm(\qb) }_\text{implicit} = \underbrace{\Nm(\qb^*)}_\text{nonlinear}+ \underbrace{\Lm(\qb)}_\text{linear} \\
= \underbrace{\Rm(\qb^*)}_\text{predictor}+ \underbrace{\big[ \Lm(\qb) - \Lm(\qb^*) \big]}_\text{corrector}
 \end{split}
\end{equation}
where $\qb^*$  represents values at previous time steps or stages. We can view the operator splitting as an explicit/implicit, nonlinear/linear or a predictor/corrector approach.  Terms computed from $\qb^*$ are explicit while those computed from $\qb$ are implicit.

In operational NWP, the aspect ratio between the horizontal and vertical grids is often very high, typically in the order of 1000. This has led to methods which exploit this property by using different time stepping methods  for the horizontal and vertical components, known as horizontally-explicit vertically-implicit (HEVI) \cite{ikawa1988,lock2014,weller2013,bao2015,giraldo2013} schemes, henceforth referred as 1D-IMEX methods. In these methods, only the vertically propagating fast waves are modeled implicitly; the horizontally propagating fast waves, hopefully, do not raise the horizontal Courant number above 1. In the idealistic case of an aspect ratio of 1, 1D-IMEX  should not bring any benefit because one would not be able to increase the time step without violating the horizontal CFL condition. The 1D-IMEX method is formulated with the following splitting of the right hand side
\begin{equation}
\label{sieqn1d}
\Rm(\qb) = \underbrace{\Rm(\qb^*)}_\text{predictor}+ \underbrace{\big[ \Lm_V(\qb) - \Lm_V(\qb^*) \big]}_\text{corrector}
\end{equation}
where $\Lm_V(\qb) = \Lm(\qb)\cdot \mathbf{\hat{r}}$ is the vertical/radial component of the linear operator $\Lm(\qb)$, and $\mathbf{\hat{r}}$ is the unit vector in the radial direction. In parallel 1D-IMEX, the domain decomposition is done only in the horizontal direction, hence, implicit treatment in the vertical direction does not require additional communication between processors. This allows 1D-IMEX to scale as well as explicit methods  on massively parallel supercomputers (see, e.g., \citep{andreas2016}).

 In the following, we consider two classes of IMEX time stepping methods: linear multi-step and multi-stage methods.

\subsection{Linear multi-step methods}
A generic implicit linear K-step method is given as follows
\begin{equation}
\label{tieqn}
\qb^{n+1} =\sum_{k=0}^{K-1} \alpha_k \qb^{n-k} + \chi \Delta t \sum_{k=-1}^{K-1} \gamma_k  \Rm(\qb^{n-k}) 
\end{equation}
where $\qb^{n}$ denotes the solution at time level $n\Delta t$, for n = 0,1,..., and $\alpha_k,\gamma_k,\chi$ are constants defining the method. 
The fact that the summation of the $\Rm(\qb)$ terms starts from $k=-1$ implies the method is implicit. 
\comment{
For the fully implicit methods $\gamma_k = 0$ for all $k\ge0$ and $\gamma_{-1}=1$, thus
\[
\qb^{n+1} = \sum_{k=0}^{K-1} \alpha_k \qb^{n-k} + \chi \Delta t \big[\Rm(\qb^{n+1})\big].
\]
}
To derive an IMEX BDF method, we split $\Rm(\qb^{n+1})$ into explicit and implicit components using Eq.\ \eqref{sieqn}
\begin{equation}
\label{sibdf}
\qb^{n+1} = \sum_{k=0}^{K-1} \alpha_k \qb^{n-k} + \chi \Delta t \big[\Nm(\qb^{n+1})_{e} + \Lm(\qb^{n+1})\big].
\end{equation}
We need to find approximations for the explicit component without degrading the order of the combined IMEX method. We achieve this using Lagrange interpolation of the appropriate order
\[
\Nm(\qb^{n+1})_{e} = \sum_{k=0}^{K-1} \beta_k \Nm(\qb^{n-k}),
\]
where the coefficients are obtained from
\[
\beta_k =  \prod_{\substack{j=0}}^{K-1}  \frac{t_{n+1}- t_{n-j}}{t_{n-k} - t_{n-j}}.
\]

Substituting the approximation for the explicit term in Eq.\ \eqref{sibdf} yields
\[
\begin{aligned}
\begin{split}
\qb^{n+1} &= \sum_{k=0}^{K-1} \alpha_k \qb^{n-k} +  \chi \Delta t \bigg[\Lm(\qb^{n+1}) + \sum_{k=0}^{K-1} \beta_k \Nm(\qb^{n-k})\bigg] \\
                &= \sum_{k=0}^{K-1} \alpha_k \qb^{n-k} + \chi \Delta t \bigg[ \Lm(\qb^{n+1})\\
                &+ \sum_{k=0}^{K-1} \beta_k \big[ \Rm(\qb^{n-k}) - \Lm(\qb^{n-k}) \big]  \bigg] .
\end{split}
\end{aligned}
\]
Defining the explicit part that estimates  $\qb^{n+1}$ as
\[
\qb^{n+1}_{e} = \sum_{k=0}^{K-1} \alpha_k \qb^{n-k} + \chi \Delta t  \sum_{k=0}^{K-1} \beta_k  \Rm(\qb^{n-k}) 
\]
we get the predictor-corrector equation
\[
\begin{aligned}
\qb^{n+1} &=\qb^{n+1}_{e} + \chi \Delta t \bigg[ \Lm(\qb^{n+1}) - \sum_{k=0}^{K-1} \beta_k \Lm(\qb^{n-k}) \bigg] \\
                &=\underbrace{\qb^{n+1}_{e}}_\text{predictor} + \underbrace{\chi \Delta t \Lm\bigg[ \qb^{n+1} - \sum_{k=0}^{K-1} \beta_k \qb^{n-k} \bigg]}_\text{corrector}\\
\end{aligned}
\]
where the last step is possible due to the linear property of the operator $\Lm$. This exposes the IMEX method as a predictor-corrector approach where we first obtain an estimate using
a standard explicit time marching algorithm, and then correct the result with terms built from the implicit linear operator $\Lm$. Introducing the following variables
\[
\begin{aligned}
& \qtt = \qb^{n+1} - \sum_{k=0}^{K-1} \beta_k \qb^{n-k} \\
& \qtte = \qb^{n+1}_{e} - \sum_{k=0}^{K-1} \beta_k \qb^{n-k}\\
\end{aligned}
\]
and $\lambda = \chi \Delta t $, we simplify the formulation  to
\begin{equation}
\label{imex}
\qtt = \qtte + \lambda \Lm(\qtt).
\end{equation}
For example, the coefficients for the two-step BDF2 method are $\alpha_0 = 4/3, \alpha_1 = -1/3, \chi=2/3, \beta_0 =  2$ and $\beta_1 = -1$. This IMEX BDF method
was first proposed by \citet{karniadakis1991} for the incompressible Navier-Stokes equations, and later used by \citet{shenwang1999} for solving the primitive equations of the atmosphere.

\subsection{Linear multi-stage methods}
We consider IMEX Runge-Kutta methods -- often denoted as Additive Runge-Kutta (ARK) methods when used in the IMEX context. These methods
use two different time integrators for the stiff and non-stiff parts. ARK methods can be represented compactly with the following double Butcher tableaux \citep{butcher2003}:
\begin{equation}
\label{ark}
\centering
\begin{aligned}
\bigg(
\begin{tabular}{c|c}
$c_i$ & $a_{ij}$\\
\hline
&$b_j$
\end{tabular},
\begin{tabular}{c|c}
$\tilde{c}_i$ & $\tilde{a}_{ij}$\\
\hline
&$\tilde b_j$
\end{tabular}
; i,j=1,\cdots,s \bigg),\\
c_i =\sum_{j=1}^{s} a_{ij} , \tilde c_i =\sum_{j=1}^{s} \tilde a_{ij}
\end{aligned}
\end{equation}

where $s$ is the number of stages, $a_{ij},b_j,c_i$ define the explicit integrator for the non-stiff terms, and $\tilde a_{ij},\tilde b_j,\tilde c_i$ define the implicit integrator for the stiff terms.
For our work, we consider singly diagonally implicit s-stage ARK methods (SDIRK) represented with stage computations $i=1,\cdots,s$ for the partitioned system in Eq.\ \eqref{sieqn}
\[
\begin{aligned}
\qb^{(i)} &= \qb^n + \Delta t \sum_{j=1}^{i-1} a_{ij} \Nm(\qb^{(j)})\\
	     &\,\,\qquad + \Delta t \sum_{j=1}^{i} \tilde a_{ij} \Lm(\qb^{(j)}) \\
             &= \qb^n + \Delta t \sum_{j=1}^{i-1} a_{ij} \Nm(\qb^{(j)})\\
             &\,\,\qquad + \Delta t \sum_{j=1}^{i-1} \tilde a_{ij} \Lm(\qb^{(j)}) + \tilde a_{ii} \Delta t \Lm(\qb^{(i)}) \\
             &= \qb^n + \Delta t \sum_{j=1}^{i-1} a_{ij} \Rm(\qb^{(j)})\\
             &\,\,\qquad + \Delta t \sum_{j=1}^{i-1} (\tilde a_{ij} - a_{ij}) \Lm(\qb^{(j)}) + \tilde a_{ii} \Delta t \Lm(\qb^{(i)})\\
             &= \qb^n + \Delta t \sum_{j=1}^{i-1} a_{ij} \Rm(\qb^{(j)})\\
             &\,\,\qquad + \tilde a_{ii} \Delta t \Lm \bigg( \qb^{(i)} + \sum_{j=1}^{i-1} \frac{\tilde a_{ij} - a_{ij}}{\tilde a_{ii}} \qb^{(j)} \bigg)
\end{aligned}
\]
where we used the linearity of $\Lm(\qb)$ in the last step. We also note that in SDIRK methods $\tilde a_{ii} = \tilde a_{jj}, i,j\le s$. 

By defining the explicit estimates  $\qb^{(i)}$ at the current stage as
\[
\qb^{(i)}_e = \qb^{n} + \Delta t \sum_{j=1}^{i-1} a_{ij} \Rm(\qb^{(j)}),
\]
we get the predictor-corrector equation
\[
\qb^{(i)} = \underbrace{\qb^{(i)}_e}_\text{predictor} + \underbrace{\tilde a_{ii} \Delta t \Lm \bigg( \qb^{(i)}  +  \sum_{j=1}^{i-1} \frac{\tilde a_{ij} - a_{ij}}{\tilde a_{ii}} \qb^{(j)} \bigg)}_\text{corrector}.
\]
We introduce the following variables
\[
\begin{aligned}
& \qtt = \qb^{(i)} + \sum_{j=1}^{i-1} \frac{\tilde{a}_{ij} - a_{ij}}{\tilde a_{ii}} \qb^{(j)} \\
& \qtte = \qb^{(i)}_{e} + \sum_{j=1}^{i-1} \frac{\tilde{a}_{ij} - a_{ij}}{\tilde a_{ii}} \qb^{(j)},
\end{aligned}
\]
and $\lambda = \tilde a_{ii} \Delta t$, which results in the simplified equation for the stages
\[
\qtt = \qtte + \lambda \Lm(\qtt).
\]
The final step is to complete the time marching by combining all the stage derivatives as
\[
\qb^{n+1} = \qb^{n} + \Delta t \sum_{i=1}^{s} b_{i} \Nm(\qb^{(i)})+ \Delta t \sum_{i=1}^{s} \tilde b_{i} \Lm(\qb^{(i)}).
\]
To ensure conservation of linear invariants, we assume $b=\tilde b$ following \cite{giraldo2013}. 

\section{IMEX formulation of the governing equations}
\label{governing}
The dynamics of non-hydrostatic atmospheric processes are governed by the compressible Euler equations. The equation sets can be written in various conservative and non-conservative forms (see e.g. \citet{giraldo:2008a}). Among those, we consider one conservative set (Set2C) and another more efficient but non-conservative set (Set2NC). First, we present the equation sets and  then extract the linear operator $\Lm$ for the IMEX methods.  The derivation of the linear operator is rather straightforward once we formulate the equation sets using perturbed states from a hydrostatically balanced reference state as $\qb(\mathbf x, t) = \qb_0(\mathbf x) + \qb'(\mathbf x, t) $. Further details on both the IMEX time-integrators used and in the implicit formulation of both cG and dG methods can be found in \cite{giraldo2013} and \cite{restelli:2009}.

\subsection{Equation Set2NC}
The five prognostic variables that comprise $\qb$ are ${(\rho, \mathbf u ^\top, \theta)}^\top$, where $\rho$ is density, $\theta$ is potential temperature and $\mathbf u=(u,v,w)^\top$ are the velocity components , and the superscript $\top$ denotes the transpose operator. We write the governing equations in the following way

\begin{equation}
\begin{aligned}
& \frac{\partial\rho}{\partial t} + \divr(\rho\mathbf u) = 0\\
& \frac{\partial\mathbf u}{\partial t} + \mathbf u \cdot \nabla \mathbf u + \frac{1}{\rho} \nabla P + \gb = 0\\
& \frac{\partial\theta}{\partial t} + \mathbf u \cdot \nabla\theta  = 0
\end{aligned}
\end{equation}
where the pressure in the momentum equation is obtained from the equation of state 
\begin{equation}
P=P_0\left(\frac{\rho R \theta}{P_0}\right)^{\gamma}
\end{equation}
and $R = c_p - c_v$ and $\gamma = \frac{c_p}{c_v}$ for given specific heat of pressure and volume of $c_p$ and  $c_v$,  respectively. 

For better numerical stability, the density, pressure and potential temperature variables are split into  background and perturbation components. The time-invariant background components are often obtained by assuming hydrostatic equilibrium and a neutral atmosphere. Let us define the decomposition of the variables into background and perturbation components as follows

\[
\begin{aligned}
\rho(\mathbf x, t)=\rho_0(\mathbf x) + \rho'(\mathbf x, t)\\
\theta(\mathbf x, t)=\theta_0(\mathbf x) + \theta'(\mathbf x, t)\\
P(\mathbf x, t)=P_0(\mathbf x) + P'(\mathbf x, t) .
\end{aligned}
\]
In this section we follow the derivation of the IMEX method for Equation Set2NC from \cite{giraldo:2010b,giraldo2013}.
We re-write the equation set in terms of the perturbation components as follows

\begin{equation}
\begin{aligned}
& \frac{\partial\rho'}{\partial t} + \mathbf u \cdot \nabla \rho' + \mathbf u \cdot \nabla \rho_0 + (\rho' + \rho_0)\nabla \cdot \mathbf u= 0\\
& \frac{\partial\mathbf u}{\partial t} + \mathbf u \cdot \nabla \mathbf u + \frac{1}{\rho_0+\rho'} \nabla P'  + \frac{\rho'}{\rho_0+\rho'}\gb = 0\\
& \frac{\partial\theta'}{\partial t} + \mathbf u \cdot \nabla\theta' + \mathbf u \cdot \nabla \theta_0 = 0 .
\end{aligned}
\end{equation}
Then, we construct the linear operator $\Lm$ for the IMEX procedure such that it includes the fastest waves in the system, namely the acoustic and gravity waves, as follows

\begin{equation}
\Lm(\qb)=-\begin{pmatrix} 
\mathbf u \cdot \nabla \rho_0 +  \rho_0 \nabla \cdot \mathbf u \\\\
\frac{1}{\rho_0} \nabla P'  + \frac{\rho'}{\rho_0}\gb \\\\
\mathbf u \cdot \nabla \theta_0
\end{pmatrix}.
\end{equation}
The perturbation pressure $P'$ is obtained from a linearization of the equation of state in the following way

\begin{equation}
P'=\frac{\gamma P_0}{\rho_0} \rho' + \frac{\gamma P_0}{\theta_0} \theta'.
\end{equation}
Substituting the linear operator in the IMEX formulation of Eq.\ \eqref{imex}, we obtain

\begin{equation}
\label{imex2nc}
\begin{pmatrix} \rtt \\\\ \utt\\\\ \ttt \end{pmatrix} =
\begin{pmatrix} \rtte \\\\ \utte \\\\ \ttte \end{pmatrix}
-\lambda\begin{pmatrix} 
\utt\cdot \nabla \rho_0 +  \rho_0 \nabla \cdot \utt\\\\
\frac{1}{\rho_0} \nabla \ptt + \frac{\rtt}{\rho_0}\gb \\\\
\utt\cdot \nabla \theta_0
\end{pmatrix}
\end{equation}
where
\begin{equation}
\label{2ncp}
\ptt=\frac{\gamma P_0}{\rho_0} \rtt + \frac{\gamma P_0}{\theta_0} \ttt = G_0 \rtt + H_0 \ttt.
\end{equation}
Similarly we define 
\[
\ptte=G_0 \rtte + H_0 \ttte.
\]

The system represented by Eqs.\ \eqref{imex2nc} -\ \eqref{2ncp} is the \nsr IMEX form. This form contains a total of $5N_p$ degrees of freedoms to solve the 3D Euler equations.
The number of degrees of freedom can be reduced to just $N_p$ by solving for pressure via the \sr complement.
\comment{reformulating the equations using a technique known in the literature by many names, including block LU decomposition, solving an equivalent
pseudo-Helmholtz problem, and solving the Schur complement of the system.
From here on we refer to this form as the \sr  form. The \sr form solves for a single pressure-like variable from an equivalent system in the form of an advection-reaction equation. }

To get the \sr form, first we make successive substitutions starting with $\ttt$ into the pressure equation
\[
\rtt = \frac{1}{G_0} \{ \ptt- H_0 \big[ \ttte - \lambda (\utt\cdot \nabla \theta_0) \big] \}.
\]
Then, we substitute this equation into the momentum equation eliminating $\rtt$ and expressing the momentum in terms of only pressure $\ptt$
\begin{equation}
\label{utteqn}
\begin{aligned}
\utt&= \utte - \lambda \{\frac{1}{\rho_0} \nabla \ptt + \frac{\gb}{\rho_0} \frac{1}{G_0} \{ \ptt -\\
     &\qquad\qquad H_0 \big[ \ttte - \lambda (\utt\cdot \nabla \theta_0) \big] \}\} \\
\utt&= \mathcal{A}^{-1} \big[ \utte - \lambda \{\frac{1}{\rho_0} \nabla \ptt + \frac{\gb}{G_0\rho_0}  ( \ptt- H_0 \ttte )\} \big] \\
\utt&= \underbrace{\mathcal{A}^{-1} \big[ \utte + \lambda \{ \frac{\gb}{G_0\rho_0}  (H_0 \ttte )\} \big]}_{\utt_{a}} \\
     &- \underbrace{\mathcal{A}^{-1} \big[  \lambda \{\frac{1}{\rho_0} \nabla \ptt + \frac{\gb}{G_0\rho_0}  ( \ptt)\} \big]}_{\utt_{p}} \\
\end{aligned}
\end{equation}
where
\[
\mathcal A = \mathbf I + \lambda^2 \frac{\gb}{\theta_0} {(\nabla \theta_0)}^\top.
\]
Substituting $\rtt$ and $\ttt$ into the pressure equation yields the \sr complement for pressure
\[
\ptt= \ptte - \lambda \mathbf{F_0} \cdot \utt- \lambda \rho_0 G_0 \nabla \cdot \utt\\
\]
where $\mathbf F_0 = G_0 \nabla \rho_0 + H_0 \nabla \theta_0$. 
Then, substituting $\utt= \utt_{a} - \utt_{p}$ yields
\begin{equation}
\label{helm2nc}
\begin{aligned}
\ptt - \lambda \mathbf{F_0} \cdot \utt_{p} &- \lambda \rho_0 G_0 \nabla \cdot \utt_{p} \\
&= \\
\ptte - \lambda \mathbf{F_0} \cdot \utt_{a} &- \lambda \rho_0 G_0 \nabla \cdot \utt_{a}
\end{aligned}
\end{equation}
where $\utt_{p}(\ptt,\nabla \ptt)$ is given in Eq.\ \eqref{utteqn}.

\subsection{Equation Set2C}
The five prognostic variables for this equation set are ${(\rho, \mathbf U ^\top, \Theta)}^\top$, where $\rho$ is density, $\mathbf{U}=\rho \mathbf u$,  $\Theta=\rho \theta$, $\theta$ is potential temperature and $\mathbf u=(u,v,w)^\top$ are the velocity components. We write the governing equations in the following way

\begin{equation}
\begin{aligned}
& \frac{\partial\rho}{\partial t} + \nabla \cdot \mathbf{U}  = 0\\
& \frac{\partial \mathbf U}{\partial t} + \nabla \cdot \left( \frac{ \mathbf U \otimes \mathbf U}{\rho} + P\mathbf I_3 \right) + \rho \gb = 0\\
& \frac{\partial \Theta}{\partial t} + \nabla \cdot \left( \frac{\Theta \mathbf U}{\rho} \right) = 0
\end{aligned}
\end{equation}
where $\mathbf I_3$ is the rank-3 identity matrix. Below, we follow the derivation of the IMEX formulation (\sr and \nsr forms) for Set2C presented in \cite{giraldo:2010b}.

Splitting into background and perturbation components 
\[
\begin{aligned}
\rho(\mathbf x, t)=\rho_0(\mathbf x) + \rho'(\mathbf x, t)\\
\Theta(\mathbf x, t)=\Theta_0(\mathbf x) + \Theta'(\mathbf x, t)\\
P(\mathbf x, t)=P_0(\mathbf x) + P'(\mathbf x, t)
\end{aligned}
\]
and re-writing the equation set using perturbation components gives
\begin{equation}
\label{eulereq}
\begin{aligned}
& \frac{\partial\rho'}{\partial t} + \nabla \cdot \mathbf U = 0\\
& \frac{\partial \mathbf U}{\partial t} + \nabla \cdot \left( \frac{\mathbf U \otimes \mathbf U}{\rho} + P'\mathbf I_3 \right) + \rho' \gb = 0\\
& \frac{\partial \Theta'}{\partial t} + \nabla \cdot \left( \frac{\Theta \mathbf U}{\rho} \right) = 0.
\end{aligned}
\end{equation}
where, for convenience, we have assumed that the background state is in hydrostatic balance and have eliminated those terms - although in principle this need not be done for more complicated background states.
We construct the linear operator $\Lm$ such that it includes the acoustic and gravity waves as follows

\begin{equation}
\Lm(\qb)=-\begin{pmatrix} 
\divr \mathbf U  \\\\
\nabla P'  + \rho' \gb  \\\\
\divr(\frac{\Theta_0}{\rho_0}\mathbf U)
\end{pmatrix}
\end{equation}
with the pressure linearized as follows
\[
P'=\frac{\gamma P_0}{\Theta_0} \Theta'.
\]

Substituting the linear operator in the IMEX formulation of Eq.\ \eqref{imex}, we obtain
\begin{equation}
\label{imex2c}
\begin{pmatrix} \rtt \\\\ \Utt \\\\ \Ttt \end{pmatrix} =
\begin{pmatrix} \rtte \\\\ \Utte \\\\ \Ttte \end{pmatrix}
-\lambda\begin{pmatrix} 
\divr \Utt  \\\\
\nabla \ptt + \rtt \gb  \\\\
\divr(\frac{\Theta_0}{\rho_0} \Utt)
\end{pmatrix}
\end{equation}
and
\[
\ptt=\frac{\gamma P_0}{\Theta_0} \Ttt = F_0 \Ttt.
\]
Similarly we define 
\[
\ptte=F_0 \Ttte.
\]
To get the \sr form, first we substitute $\Ttt$ into the pressure equation
\[
\ptt=F_0 \big[ \Ttte - \lambda \nabla \cdot (G_0 \Utt) \big] = \ptte - F_0 \lambda \nabla \cdot (G_0 \Utt)
\]
where $G_0=\Theta_0/\rho_0$. Subtracting $ G_0 \rtt$ from $\Ttt$, we get
\[
\Ttt - G_0 \rtt = \Ttte - G_0 \rtte - \lambda \Utt \cdot \nabla G_0.
\]
Substituting the pressure equation into the above equation
\[
\rtt = \frac{1}{F_0 G_0} \ptt+ \frac{\lambda}{G_0} \Utt \cdot \nabla G_0 - \frac{1}{G_0} \Ttte + \rtte.
\]
Substituting $\rtt$ in the momentum equation
\begin{equation}
\label{utteqn2}
\begin{aligned}
\Utt &= \mathcal A^{-1} \big[ \Utte - \lambda \nabla \ptt- \\
&\lambda \frac{\gb}{F_0 G_0} \ptt  - \lambda \gb \big(  \rtte - \frac{1}{G_0} \Ttte \big)   \big] \\
\Utt &= \underbrace{\mathcal A^{-1} \big[ \Utte   - \lambda \gb \big(  \rtte - \frac{1}{G_0} \Ttte \big)   \big]}_{\Utt_{\tt a}} \\
      &- \underbrace{\mathcal A^{-1} \big[ \lambda \big(\nabla \ptt+ \frac{\gb}{F_0 G_0} \ptt\big)  \big]}_{\Utt_{p}}
\end{aligned}
\end{equation}
where
\[
\mathcal A = \mathbf I + \lambda^2 \frac{\gb}{G_0} {(\nabla G_0)}^\top.
\]
Substituting $\Utt = \Utt_{a} - \Utt_{p}$ into the pressure equation yields the \sr  complement for pressure
\begin{equation}
\label{helm2c}
\begin{aligned}
\ptt- &F_0\lambda \nabla \cdot (G_0 \Utt_{p}) \\
&=\\
\ptte - &F_0\lambda \nabla \cdot (G_0 \Utt_{a})
\end{aligned}
\end{equation}
where $\Utt_{p}(\ptt,\nabla \ptt)$ is given in Eq.\ \eqref{utteqn2}.

\section{Spatial discretization}
We use  Set2NC primarily with the continuous Galerkin (cG) spatial discretization and the conservative Set2C with the discontinuous Galerkin (dG) discretization.
In this section, we shall discuss  the discretization of only the linear operators $\Lm(\qb)$ in Eqs.\ \eqref{imex2nc} and \eqref{imex2c} and refer the reader to our previous work \cite{abdiGiraldo2015,abdi2016} for  the discussion on discretization of the rest of the terms.

We begin by separating the linear operator $\Lm(\qb)$  into flux and source terms as follows
\[
\Lm(\qb) = -\divr \Fm(\qb) + \Sm(\qb).
\]

Assuming a domain $\Omega \in \mathbb{R}^3$ with boundary $\Gamma$, the discretization of the linear operator using the Galerkin procedure is given as follows
\begin{equation}
\label{cg}
\int_{\Omega} \psi \Lm(\qb) d\Omega= -\int_{\Omega} \psi \nabla\cdot \Fm(\qb) d\Omega + \int_{\Omega} \psi \Sm(\qb)  d\Omega
\end{equation}
where $\psi$ is the test function.

The value of separating $\Lm(\qb)$ into flux and source terms is in the discontinuous Galerkin (dG) method where face fluxes are required for imposing weak coupling between elements. We put the gravity terms in the source term $\Sm(\qb)$ of the linearized operator  and then define the  Rusanov flux  as follows
\[
\Fm(\qb)^* = \{\Fm(\qb)\} - \mathbf{\hat{n}}\frac{|\widehat{c}|}{2} [\![\qb ]\!]
\]
where $|\widehat{c}|$ is the speed of sound, $\{\}$ represents an average and $[\![ ]\!]$ represents an outward jump across a face of an element, and $\mathbf{\hat{n}}$ is the outward pointing unit normal. 
We apply the dG method to the \nsr IMEX form in a straightforward manner as
\[
\begin{split}
\int_{\Omega} \psi \Lm(\qb) d\Omega = \int_{\Omega} \nabla \psi \cdot \Fm(\qb) d\Omega - \int_{\Gamma} \psi \Fm(\qb)^* \cdot \hat{\mathbf n}  d\Gamma \\
+ \int_{\Omega} \psi \Sm(\qb)  d\Omega.
\end{split}
\]
The strong form dG, which has the same volume integral terms as Eq.\ \eqref{cg}, is obtained by taking a second integration by parts as follows
\begin{equation}
\begin{split}
\int_{\Omega} \psi \Lm(\qb) d\Omega = -\int_{\Omega} \psi \nabla\cdot \Fm(\qb) d\Omega + \int_{\Omega} \psi \Sm(\qb)  d\Omega \\
- \int_{\Gamma} \psi \hat{\mathbf n} \cdot  (\Fm(\qb)^* - \Fm(\qb) )  d\Gamma.
\end{split}
\end{equation}

\section{Implementation on manycore processors}
In this section, we describe the implementation of the infrastructure required for conducting IMEX time integration on manycore processors. In our previous work \cite{abdi2016}, we presented a GPU acceleration of NUMA and its scalability on tens of thousands of GPUs using explicit time integration; here, we shall only discuss the new infrastructure required for enabling IMEX time stepping and refer the reader to our previous work for a complete view of the IMEX time-integration approach (i.e., the explicit part). To summarize the new additions, we need: a) kernels for evaluating left-hand side and right-hand side operators for the implicit terms described in Sec.\ \ref{governing}. This will be done for both 1D- and 3D- IMEX time integration, for both \sr and \nsr forms, for both Set2NC and Set2C, and for both cG and dG. We will also need b) different types of iterative solvers and preconditioners for solving the linearly discretized  system. The system resulting from the 1D-IMEX discretization can also be conveniently solved using direct methods because the Jacobian matrix for a `column' is small. In fact, this is the preferred choice when running on the CPU; however, direct solvers may not be clear winners for solving small matrices on the GPU. In both cases, the addition of implicit solves at every time step/stage negatively impacts performance on one node and also scalability on multi-node clusters.

The implementation of IMEX on manycore processors is broken down into different kernels, of which the major ones are: a) volume kernels for evaluating the left- and right- hand sides of both IMEX forms;  b) surface kernels for computing flux integrals at the trace of elements;  and c) kernels for extracting results of all prognostic variables  from the sole pressure variable used in the \sr form. We also need kernels for solving the resulting system of equations iteratively, or directly in the case of 1D-IMEX.

\subsection{Volume kernels}
The volume kernels for the \nsr forms involve first order terms, while those for the \sr forms involve second-order terms. Because we are primarily using iterative methods for solving the implicit problem, the second order terms are evaluated by first computing $\nabla P$ and then $\nabla\cdot\mathbf{f}(P,\nabla P)$ in the spirit of  the local discontinuous Galerkin method. Rewriting the general IMEX problem of Eq.\ \eqref{imex} as 
\[
(\mathbf{I} - \chi\Delta t \Lm)\qtt = \qtte,
\]
we separate the left- and right- hand side terms. The volume kernels for evaluating the right-hand-side terms of the \nsr forms, shown in Eqs.\ \eqref{imex2nc} and \eqref{imex2c}, are rather straightforward to implement because they are simply estimates $\qtte$ obtained by using explicit time stepping in the predictor-corrector approach. Therefore, these kernels basically set the right-hand side terms to $(\rtte, \utte, \ttte)^\top$ and $(\rtte, \Utte, \Ttte)^\top$ for Set2NC and Set2C, respectively. On the other hand, the left-hand side term $(\mathbf{I} - \chi\Delta t \Lm)\qtt$ requires either constructing the Jacobian matrix  in the case of the direct solution of the 1D-IMEX problem, or evaluating the gradient of pressure, $\nabla \ptt$, divergence of velocity, $\nabla\cdot\Utt$, and other terms included in the linear operator in the case of an iterative solution.

\begin{algorithm*}
\caption{Algorithm for computing gradient or divergence in 1D-IMEX.}
\label{gradientOutline1d}
\begin{algorithmic}
\\
\Procedure{GradDiv}{$\qtt$, Jr, GD} 
    \State \textbf{Shared} sD[$N_q$][$N_q$]   \Comment{sD are $\nabla\psi$ at LGL nodes pre-loaded in shared mem.}
    \State \textbf{Shared} sq[$N_q$][$N_q$][$N_q$] \Comment{Shared space for collaborative computation}
    \State Local memory fence
    \For{k,j,i $\in \{0 \dots N_q$\}} \Comment{Load field variables into shared memory}
        \State sq[k][j][i] = $\qtt_{ijk}$
    \EndFor
    \State Local memory fence
    \For{k,j,i $\in \{0 \dots N_q$\}}
        \State qr=0 \Comment{Compute local gradients}
        \For{n $\in \{0 \dots N_q$\}}
             \State qr += sD[k][n]$\times$sq[n][j][i]
        \EndFor
        \State GD$_{ijk}$ = qr $\times$ Jr$_{ijk}$\Comment{Jr are coefficients of the 1D-Jacobian matrix}
    \EndFor
\EndProcedure

\end{algorithmic}
\end{algorithm*}

For the 1D-IMEX method, gradient and divergence in the radial directions are one and the same. The algorithm for this purpose is shown in Alg.\ \ref{gradientOutline1d}. On the other hand, 3D-IMEX requires computing gradient in all three Cartesian directions of a tensor-product element as discussed in detail in \citet{abdi2016}, and the algorithm repeated here in Alg.\ \ref{gradientOutline3d}. We should note that we evaluate nodal contributions to the left- and right- hand sides element by element for both 1D- and 3D- IMEX; alternatively, we could have used a column by column approach in the case of 1D-IMEX where each column extends through all the vertical layers. In fact, to use the direct solver in 1D-IMEX in which we build matrices for each column, the grid should be generated as a set of columns. \comment{Note that the definition of a `column' is only applicable for a vertically conforming grid.}
\comment{that we make sure of using the p6est adaptive mesh refinement (AMR) library; p6est  basically generates a 2D grid whose nodes are extruded to form columns.}

\begin{algorithm*}
\caption{GPU algorithms for computing gradient or divergence in 3D-IMEX.}
\label{gradientOutline3d}
\begin{algorithmic}
\\
\Procedure{GradDiv}{$\qtt$, GD, J, sD, compute}
    \State \textbf{Shared} sD[$N_q$][$N_q$]  \Comment{sD are $\nabla\psi$ at LGL nodes pre-loaded in shared mem.}
    \State \textbf{Shared} sq[$N_q$][$N_q$][$N_q$] \Comment{Shared space for collaborative computation}
    \State Local memory fence
    \For{k,j,i $\in \{0 \dots N_q$\}} \Comment{Load field variables into shared memory}
        \State sq[k][j][i] = $\qtt_{ijk}$
    \EndFor
    \State Local memory fence
    \For{k,j,i $\in \{0 \dots N_q$\}}
        \State qx=0; qy=0; qz=0; \Comment{Compute local gradients}
        \For{n $\in \{0 \dots N_q$\}}
        	     \State qx += sD[i][n]$\times$sq[k][j][n]
             \State qy += sD[j][n]$\times$sq[k][n][i]
             \State qz += sD[k][n]$\times$sq[n][j][i]
        \EndFor
	\If{compute = GRAD}   
	     \State  GD$_{ijk\_x}$ \text{ } = (qx $\times$ Jrx$_{ijk}$ + qy $\times$ Jsx$_{ijk}$ + qz $\times$ Jtx$_{ijk}$) \Comment{Js are the 9 coefficients of the jacobian matrix J}
        	     \State  GD$_{ijk\_y}$ \text{ } = (qx $\times$ Jry$_{ijk}$ + qy $\times$ Jsy$_{ijk}$ + qz $\times$ Jty$_{ijk}$)
        	     \State	GD$_{ijk\_z}$ \text{ } = (qx $\times$ Jrz$_{ijk}$ + qy $\times$ Jsz$_{ijk}$ + qz $\times$ Jtz$_{ijk}$)
	\ElsIf{compute = DIVX}
	      \State  GD$_{ijk\_w}$ \text{ } = (qx $\times$ Jrx$_{ijk}$ + qy $\times$ Jsx$_{ijk}$ + qz $\times$ Jtx$_{ijk}$)
	\ElsIf{compute = DIVY}
	      \State  GD$_{ijk\_w}$ += (qx $\times$ Jry$_{ijk}$ + qy $\times$ Jsy$_{ijk}$ + qz $\times$ Jty$_{ijk}$)
	\ElsIf{compute = DIVZ}
	      \State  GD$_{ijk\_w}$ += (qx $\times$ Jrz$_{ijk}$ + qy $\times$ Jsz$_{ijk}$ + qz $\times$ Jtz$_{ijk}$)
	\EndIf
    \EndFor
\EndProcedure
\\
\Procedure{Grad}{q, GD, J, sD}  \Comment{Compute gradient of a scalar field}
     \State \textbf{call} \textproc{GradDiv}(q, GD, J, sD, GRAD)
\EndProcedure
\\
\Procedure{Div}{$\qb$, GD, sD, J} \Comment{Compute divergence of a vector field}
     \State \textbf{call} \textproc{GradDiv}($\qb\cdot$x, GD, J, sD, DIVX)
     \State \textbf{call} \textproc{GradDiv}($\qb\cdot$y, GD, J, sD, DIVY)
     \State \textbf{call} \textproc{GradDiv}($\qb\cdot$z, GD, J, sD, DIVZ)
\EndProcedure
\\
\end{algorithmic}
\end{algorithm*}

The volume kernel for the left-hand side of the \nsr form is shown in Alg.\ \ref{lhsStandard}. First, we load the geometric factors, field variables  and reference values to thread private memory. The major computation in this kernel is computing the gradient of the pressure and divergence of the velocity field. Gradient and divergence of reference fields are instead computed once at startup and passed to the kernel.
\begin{algorithm*}
\caption{Volume kernel for computing left-hand side of the 3D-IMEX \nsr form.}
\label{lhsStandard}
\begin{algorithmic}
\\
\Procedure{LHS-Standard-Set2C}{$\qtt$,J,D,$G_0$,$\nabla G_0$} 
    \State \textbf{Shared} sD[$N_q$][$N_q$]  \Comment{sD are $\nabla\psi$ at LGL nodes.}
    \State \textbf{Shared} sq[$N_q$][$N_q$][$N_q$] \Comment{Shared space for collaborative computation}
    \State \textbf{Private} $Q^{tt}$ \Comment{Private variable for current values of prognostic fields; $Q^{tt} = \{\Utt,\rtt,\Ttt\}$}
    \For{j,i $\in \{0 \dots N_q$\}}
        \State sD[j][i] = D[i][j]
    \EndFor
    \For{k,j,i $\in \{0 \dots N_q$\}}  \Comment{ Read $\qtt_{ijk}$  to thread private memory}
        \State $Q^{tt}$ = $\qtt_{ijk}$
    \EndFor
    \State \textbf{call} \textproc{Div}($\Utt$, $\divr\Utt$, J, sD)
    \State \textbf{call} \textproc{Grad}($\ptt$, $\nabla\ptt$, J, sD)
    \For{k,j,i $\in \{0 \dots N_q$\}}
    	\State LHS$\cdot\rtt$ = $\rtt + \lambda\nabla\cdot\Utt$
	\State LHS$\cdot\Utt$ = $\Utt + \lambda(\nabla\ptt + \rtt \gb)$
	\State LHS$\cdot\Ttt$ = $\Ttt + \lambda(G_0\divr\Utt + \nabla G_0 \cdot \Utt)$
    \EndFor
\EndProcedure
\end{algorithmic}
\end{algorithm*}

Both the right- and left- hand side evaluations of the \sr form involve non-trivial operations as shown in Algs.\ \ref{rhsSchur} and \ref{lhsSchur}. In the case of a direct solution procedure, the \sr form has only pressure as the solution variable reducing the size of the Jacobian matrix significantly. Its benefit with regard to iterative methods is that the condition number of the \sr form is much better than the \nsr form (and the eigenvalues are all real), thereby, leading to fewer iterations for convergence (see \cite{giraldo:2010b,giraldo2013}). An additional kernel for `extracting' the five prognostic variables from pressure is required in the case of the \sr form, which  adds some cost to the \sr form.
\begin{algorithm*}
\caption{Volume kernel for computing right-hand side of the 3D-IMEX \sr form.}
\label{rhsSchur}
\begin{algorithmic}
\\
\Procedure{RHS-Schur-Set2C}{$\qtte$,$\mathbf{p}_e^{tt}$,J,D,$G_0$,$\nabla G_0$} 
    \State \textbf{Shared} sD[$N_q$][$N_q$]  \Comment{sD are $\nabla\psi$ at LGL nodes.}
    \State \textbf{Shared} sq[$N_q$][$N_q$][$N_q$] \Comment{Shared space for collaborative computation}
    \State \textbf{Private} $Q_e^{tt}$,$P_e^{tt}$ \Comment{Private variable for explicit estimate of prognostic fields; $Q_e^{tt} = \{\Utte,\rtte,\Ttte\}$}
    \For{j,i $\in \{0 \dots N_q$\}}
        \State sD[j][i] = D[i][j]
    \EndFor
    \For{k,j,i $\in \{0 \dots N_q$\}}  \Comment{ Read $\ptt_{e\_ijk}$  to thread private memory}
        \State $Q_e^{tt}$ = $\qtt_{e\_ijk}$
        \State $P_e^{tt}$ = $\mathbf{p}_{e\_ijk}^{tt}$
    \EndFor
    \For{k,j,i $\in \{0 \dots N_q$\}}
    	\State Compute $\mathcal A = \mathbf I + \lambda^2 \frac{\gb}{G_0} {(\nabla G_0)}^\top$
	\State Compute $\Utt_{\tt a} = \mathcal A^{-1} \big[ \Utte   - \lambda \gb \big(  \rtte - \frac{1}{G_0} \Ttte \big)   \big]$
    \EndFor
    \State \textbf{call} \textproc{Div}($G_0\Utt_{\tt a}$, $\nabla \cdot$ ($G_0\Utt_{\tt a})$, J, sD)
    \For{k,j,i $\in \{0 \dots N_q$\}}
    	\State RHS$\cdot\ptt$ = $\ptte- F_0\lambda \nabla \cdot (G_0 \Utt_{\tt a}) $
    \EndFor
\EndProcedure

\end{algorithmic}
\end{algorithm*}

\begin{algorithm*}
\caption{Volume kernel for computing left-hand side of the 3D-IMEX \sr form.}
\label{lhsSchur}
\begin{algorithmic}
\\
\Procedure{LHS-Schur-Set2C}{$\mathbf{p}^{tt}$,J,D,$G_0$,$\nabla G_0$} 
    \State \textbf{Shared} sD[$N_q$][$N_q$]  \Comment{sD are $\nabla\psi$ at LGL nodes.}
    \State \textbf{Shared} sq[$N_q$][$N_q$][$N_q$] \Comment{Shared space for collaborative computation}
     \State \textbf{Private} $P^{tt}$ \Comment{Private variable for pressure}
    \For{j,i $\in \{0 \dots N_q$\}}
        \State sD[j][i] = D[i][j]
    \EndFor
    \For{k,j,i $\in \{0 \dots N_q$\}}  \Comment{ Read $\ptt_{ijk}$ to thread private memory}
        \State $P^{tt}$ = $\mathbf{p}^{tt}_{ijk}$
    \EndFor
    \State \textbf{call} \textproc{Grad}($\ptt$, $\nabla \ptt$, J, sD)
    \For{k,j,i $\in \{0 \dots N_q$\}}
    	\State Compute $\mathcal A = \mathbf I + \lambda^2 \frac{\gb}{G_0} {(\nabla G_0)}^\top$
	\State Compute $\Utt_p = \mathcal A^{-1} \big[ \lambda \big(\nabla \ptt+ \frac{\gb}{F_0 G_0} \ptt\big)  \big]$
    \EndFor
    \State \textbf{call} \textproc{Div}($G_0\Utt_p$, $\nabla\cdot (G_0\Utt_p)$, J, sD)
    \For{k,j,i $\in \{0 \dots N_q$\}}
    	\State LHS$\cdot\ptt$ = $\ptt- F_0\lambda \nabla \cdot (G_0 \Utt_{p}) $
    \EndFor
\EndProcedure

\end{algorithmic}
\end{algorithm*}

\subsection{Surface kernel}
In the case of dG-IMEX, we need to compute flux integrals of the implicit problem on the trace of the elements. Similar to the standard surface kernel used for the right-hand side (see \citet{abdi2016} for details), we use the Rusanov flux. The only difference is that here we use the reference fields, from which the Jacobian of the left-hand side is built, to compute the maximum wave speed $c$. The maximum wave speed used for computing the right-hand side  fluxes has an additional term $|\mathbf {\hat{n}} \cdot \mathbf u| $ besides the speed of sound. The procedure for computing the flux integrals of the \nsr form dG-IMEX is shown in Alg.\ \ref{fluxStandard}.

Computing the flux term for the \sr form dG-IMEX is complicated because only the pressure variable is used during solution. This poses a problem of computing fluxes for pressure and intermediate velocity-like variables that will, hopefully, recover the same solution as that of the \nsr form. Fluxes for the right- and left- hand side terms of the \sr form come from application of the divergence operator in $\nabla \cdot (G_0 \Utt_{p})$ and $\nabla \cdot (G_0 \Utt_{a})$, respectively. This can be viewed as a Local Discontinuous Galerkin (LDG) approach where $\Utt$ is the auxiliary variable. Therefore, when we compute $\nabla\ptt$, that is used in constructing $\Utt_p(\ptt,\nabla\ptt)$, we have to apply fluxes. Unfortunately, using centered fluxes of \citet{bassi1997}  does not converge for the IMEX problem. We should note here that we can use the same kernel for computing the \sr form dG-IMEX fluxes for both the right- and left- hand side, except for the fact that we do not need to compute $\nabla\ptt$ for the right-hand side term.

\begin{algorithm*}
\caption{Surface kernel for the dG-IMEX \nsr form}
\label{fluxStandard}
\begin{algorithmic}
\Procedure{surfaceKernel}{$\qtt,S$}
    \For{id=0 \textbf{to} 6} \Comment{Six faces of the hexahedron}
        \For{j,i $\in \{0 \dots N_q$\}}
           \State Load face normal $\mathbf {\hat{n}}$ and lift coefficient $L$ \Comment{$L = \frac{w_{ij} J_{ij}}{w_{ijk} J_{ijk}}$}
           \State Load $\qtt_+$ and $\qtt_-$ for current node and adjoining node in the other element 
           \State Compute sound speed $|c| = \sqrt{\gamma \ptt_0 / \rtt_0}$
           \State Compute Rusanov flux $ \Fm(\qb)^* = \{\Fm(\qb)\} - \mathbf {\hat{n}}\frac{|c|}{2} [\![\qb ]\!]$
           \State S += $\lambda L \mathbf{\hat{n}}\cdot (\Fm(\qb)^* - \Fm(\qb))$
           \State Global memory fence
        \EndFor
     \EndFor
\EndProcedure
\end{algorithmic}
\end{algorithm*}

\subsection{Direct stiffness summation}
For the continuous Galerkin method, we need to apply the direct stiffness summation operator in places where fluxes are computed in the dG counterpart (see \cite{abdiGiraldo2015} for details). This ensures a strong coupling of the elements in the cG method. The DSS operator needs to be applied not only for the right-hand side term but also the left-hand side term to ensure convergence of the iterative solution. In the case of a direct solution to the 1D-IMEX problem, we need to modify the Jacobian matrix as described in \cite{abdiGiraldo2015}.

\subsection{1D-IMEX kernels}
For conducting global simulations on the sphere using 1D-IMEX, we need some additional kernels. First of all, vector fields in the Cartesian coordinate system need to be rotated to and from the spherical coordinate system before evaluating the right-hand side and after completing the iterative solution, respectively. Hence, we build and store a small rotation matrix $A \in \mathbb R^{3\times3}$ and its inverse at each node. The second set of kernels specific to the 1D-IMEX are required for the direct solution of the IMEX problem. Direct solution necessitates use of the `column by column' approach in which we build and store the Jacobian matrix for each column. The nodes in a column are coupled only with those nodes in the same column, therefore the global Jacobian matrix is block-diagonal if we ensure consecutive node numbering for nodes in the same column. The coupling of a node with horizontally-adjacent nodes is handled in the explicit update stage of, e.g., the ARK method. Building the Jacobian matrix of a column requires repeated evaluations of the left-hand term as shown in Alg.\ \ref{Jacobian}. Since the number of nodes in one column could be too small for a whole GPU device, it is important to be able to process multiple columns while building the Jacobian as well as later during the direct solution stage. To evaluate the Jacobian matrix coefficients (`influence' coefficients) of a degree of freedom (e.g. $\rtt$) at a node $i$, first we construct a vector $\qtt_i$ with $\rtt_i=1$  at that node and all other degrees of freedom set to 0. We can evaluate influence coefficients of  degrees of freedom at the same level in all columns simultaneously, because this influence does not diffuse to adjacent columns. This virtual independence of columns during the implicit solve stage results in a block-diagonal global Jacobian matrix which can be solved very easily. The Jacobian matrix needs to be computed only once at startup because of our choice of using a constant (in time) background state that results in the linear property of the operator $\Lm(\qb)$, hence the cost of building the Jacobian matrix is not critical to performance. After that we compute the LU decomposition (see Alg.\ \ref{ludcmp}) of the block-diagonal matrix and store it to accelerate the direct solution of the system for different right-hand sides. Parallelizing the LU decomposition and the ensuing forward-backward substitution is difficult on manycore processors due to the a) sequential nature of the algorithm and b) small size of the matrices. We relieve the second problem by solving multiple columns simultaneously. We should mention that using optimized libraries for solving block-diagonal matrices on the GPU is probably the best approach. Direct solution of  the 3D-IMEX problem on the GPU may also be feasible because even though the matrix is not block-diagonal it is tightly banded.

\begin{algorithm*}
\caption{Algorithm for building the Jacobian of 1D-IMEX}
\label{Jacobian}
\begin{algorithmic}
\Procedure{BuildJacobian}{$\qtt,S$}
     \For{nd $ \in \{0 \dots N_r$\}}         \Comment{$N_r$ is the number of nodes in a column}
        \For{dof $ \in \{0 \dots N_{dof}$\}}  \Comment{$N_{dof}$=5 for the \nsr form and $N_{dof}$=1 for the \sr form.}
           \\
           \State Construct $\qtt$ with the value set to 1 at $(dof,nd)$ and 0 every where else in the columns
           \State Compute the left-hand side LHS($\qtt$,$\dots$) to get `influence coefficients' of the specific $dof$
           \State Copy the `influence coefficients' to the appropriate locations in the Jacobian matrix
        \EndFor
     \EndFor
     \\
     \State Compute and store the LU decomposition of the block-diagonal matrix
\EndProcedure
\end{algorithmic}
\end{algorithm*}

\subsection{Iterative solvers}
The  system of equations resulting from the IMEX discretization are not necessarily symmetric; this excludes several solvers and preconditioners designed for symmetric-positive definite (SPD)
matrices. NUMA uses Krylov subspace methods such as the Generalized Minimal Residual (GMRES), stabilized bi-conjugate Gradient (BiCGstab), and simple iterative schemes such as the Richardson iteration. 
With the use of the proper preconditioner, the Richardson method could be competitive with the Krylov methods. The set of preconditioners used in NUMA are specifically designed to exploit the characteristics
of the matrices resulting from both the \nsr and \sr forms which have different eigenspectrum (see \cite{Carr2016,Carr2012} for details). Our requirement on the preconditioners are  that they be: a) amenable to vectorization (e.g. GPUs) and scalable to massively parallel systems (supercomputers); b) should be easy to apply the action of the preconditioner; and that c) the reduction in the number of iterations outweigh the cost of constructing and applying the preconditioner. For instance, requirements (a) and (b) forced the first author to use diagonal preconditioners, such as Jacobi and diagonal-ILU, in the parallel version of an incompressible Navier-Stokes solver in \cite{abdi9}. Element-based preconditioning of finite-elements (e.g see \cite{Augarde2006}) is an attractive approach that satisfies all three requirements. In \cite{Carr2012}, element-based spectrally optimized (EBSO)  sparse approximate inverse (SAI) preconditioners are investigated for use in the massively-parallel atmospheric model NUMA. The spectrally optimized preconditioner outperformed other equally parallel SAI precondioners such as the low-order Chebyshev generalized least-squares polynomials and an element-based variant of the Frobenius norm optimization procedure. In \cite{Carr2016}, polynomial based non-linear least squares optimized (PBNO) preconditioners are investigated for use in NUMA and shown to outperform generalized linear least squares (GLS) polynomial preconditioners. It is shown that a high order PBNO preconditioning of Richardson iteration makes the method competitive with the Krylov methods when run in serial mode. This is rather good news for the Richardson iteration because its dot-product free nature suggests that even better performance maybe obtained when run in parallel mode.

The sparse approximate inverse (SAI) for the PBNO preconditioner is given by
\[
s(A) = \sum_{i=0}^k c_i A^i
\]
where $A$ is the Jacobian matrix with the eigenspectrum in $[\lambda_{min}, \lambda_{max}]$. Because the Jacobian matrix is constant for the linear IMEX problem, the coefficients of the preconditioner, $c_i$ and $\overline \lambda$, are computed once at start up to the desired accuracy. The procedure for applying the PBNO preconditioner on a stabilized bi-conjugate gradient solver is shown in Alg.\ \ref{bicgstab}. We note that both the preconditioner and iterative solver make repeated calls to the left-hand side evaluator -- which is equivalent to a matrix-vector product. Most of the operations in the solver and the preconditioner are easy to implement on many-core and multi-core processors; optimized BLAS libraries (e.g. cuBLAS) can be used to implement the algorithm efficiently. Besides the left-hand side evaluations, the most time consuming operations are the dot-product operations, $\textproc{Dot}$ and $\textproc{TwoNorm}$, which are often bottlenecks of performance on both many- and multi-core processors. Often times the reduction operations are implemented efficiently using binary-tree algorithms that are $\mathcal{O}(log_2 N)$. PBNO preconditioned Richardson iteration is dot-product free which makes  it an attractive alternative to Krylov solvers due to its simplicity and efficiency of implementation.

\begin{algorithm*}
\caption{PBNO preconditioned stabilized bi-conjugate gradient}
\label{bicgstab}
\begin{algorithmic}
\Procedure{BiCGstabPBNO}{$q,b,\lambda$}
\State $q = 0$
\State  $r_{n}$ = \textbf{call} \textproc{PBNO}($b,\lambda$) \Comment{Apply PBNO preconditioner}
\State ro = \textbf{call} \textproc{TwoNorm}($r_{n}$) \Comment{Compute the $L_2$ norm of $r_n$.}
\State p = $r_{n}$
\State $r_o$ = $r_{n}$
\While {not converged}
\State t =\textbf{call} \textproc{LHS}($p,\lambda$)
\State  $Ap$ = \textbf{call} \textproc{PBNO}($t,\lambda$) 
\State $a_r = \frac{\textproc{Dot}(r_{n},r_o)}{\textproc{Dot}(Ap,r_o)}$ \Comment{Dot products involve expensive reduction operation.}
\State s = $r_n - a_r  Ap$
\State t =\textbf{call} \textproc{LHS}($s,\lambda$)
\State  $As$ = \textbf{call} \textproc{PBNO}($t,\lambda$) 
\State $w_r = \frac{\textproc{Dot}(As,s)}{\textproc{Dot}(As,As)}$
\State $r_e = a_r  p + w_r  s$
\State $q = q + r_e$
\State $r_e = s - w_r  As$
\State rn = \textbf{call} \textproc{TwoNorm}($r_e$)
\If{$\frac{rn}{r0} \le eps$} \Comment{eps is specified tolerance for terminating iterations}
\State \textbf{exit}
\EndIf
\State $b_r = \frac{\textproc{Dot}(r_e,r_o)}{\textproc{Dot}(r_n,r_o)} \times \frac{a_r}{w_r}$
\State $p = r_e + b_r (p - w_r  Ap )$
\State $r_n = r_e$
\EndWhile
\EndProcedure
\\
\Procedure{PBNO}{$r,\lambda$}
\State $r =  \overline{\lambda} r $ \Comment{The scaling factor is $\overline{\lambda} = \frac{2}{\lambda_{min} + \lambda_{max}} $}
\State $r^* = c_0 r$
\For{i $ \in \{1 \dots P$\}} \Comment{P is the order of the preconditioner}
\State m =\textbf{call} \textproc{LHS}($r^*,\lambda$)
\State $r^* = \overline{\lambda} m + c_i r$
\EndFor
\EndProcedure
\end{algorithmic}
\end{algorithm*}

\begin{algorithm*}
\caption{LU decomposition of small matrices on GPU. Multiple columns are processed simultaneously}
\label{ludcmp}
\begin{algorithmic}
\Procedure{LUdecompose}{$A$}
\State $n_{id}= e_{id} * N_q + ne_{id}$ \Comment{$e_{id}$=element id, $ne_{id}$ = node id in the element.}
\For{k $ \in \{0 \dots M-1$\}} \Comment{The matrices are MxM}
\State $E =min(k+N_b, M)$ \Comment{$N_b$ is the bandwidth of the matrix}
\For{$i = k+1+n_{id}; \text{   } i < E; \text{   } i = i + N_t$}  \Comment{\textbf{ParFor} using $N_t$ threads.}
\State $A_{ip} \mathrel{/=} A_{pp}$
\EndFor
\State Local memory fence
\For{$i = k+1+n_{id}; \text{   } i < E; \text{   } i = i + N_t$}  \Comment{\textbf{ParFor} using $N_t$ threads.}
\For{j $ \in \{k+1 \dots E$\}}
\State $A_{ij} \mathrel{-=} A_{pj} * A_{ik}$
\EndFor
\EndFor
\State Local memory fence
\EndFor
\EndProcedure
\end{algorithmic}
\end{algorithm*}
            
 \begin{table*}
\centering
\caption{Double Butcher tableau for the second-order ARK scheme of \citet{giraldo2013}.}
\begin{tabular}{c|ccc}
$0$ & $0$ &  & \\
$2-\sqrt{2}$ & $2-\sqrt{2}$ & 0&  \\
$1$ & $1-\frac{1}{6}(3+2\sqrt{2})$ & $\frac{1}{6}(3+2\sqrt{2})$ & 0\\
\hline 
&$-\frac{1}{2\sqrt{2}}$ & $-\frac{1}{2\sqrt{2}}$ & $1-\frac{1}{2\sqrt{2}}$
\end{tabular},
\begin{tabular}{c|ccc}
$0$ & $0$ &  & \\
$2-\sqrt{2}$ & $1-\frac{1}{\sqrt{2}}$ & $1-\frac{1}{\sqrt{2}}$&  \\
$1$ & $-\frac{1}{2\sqrt{2}}$ & $-\frac{1}{2\sqrt{2}}$ & $1-\frac{1}{\sqrt{2}}$\\
\hline 
&$-\frac{1}{2\sqrt{2}}$ & $-\frac{1}{2\sqrt{2}}$ & $1-\frac{1}{\sqrt{2}}$
\end{tabular}
\label{ark2table}
\end{table*}

%
%

\section{Performance of IMEX on manycore processors}
In this section, we evaluate the single node performance of the IMEX scheme on two different manycore architectures, namely, the Graphic Processing Unit (GPU) and Intel's Knights Landing (KNL). In order to use only one kernel language to take advantage of different types of processor architectures, we use OCCA\footnote{OCCA stands for Open Concurrent Compute Abstraction; further information on OCCA can be found at \cite{medina2014}.}. The OCCA kernels translate to CUDA/OpenCL code running in one-node-per-thread mode for the former, while they translate to OpenMP code with the one-element-per-node approach for the latter. We should note that in our previous work \cite{abdi2016}, we only considered GPUs. On both systems, we measure performance of the IMEX kernels in terms of two metrics: the rate of floating point operations per second (GFLOPS/s) and the rate of data transfer (GB/s). For all tests we use  the second-order ARK method  of \citet{giraldo2013} unless specified otherwise. The double Butcher tableaux for this L-stable ARK2 scheme is given in Table\ \ref{ark2table}. This method has a second-order accuracy for the explicit part and about third-order accuracy for the implicit part. For the purpose of comparison against an explciit method, we use the strong stability preserving (SSP) third-order 5-stage Runge-Kutta (RK35) method of \citet{spiteri2002}. The RK35 method has a larger stability region than the classic fourth-order RK method (see, e.g., \citet{giraldoRestelli2010} ), while being closer to ARK2 with regard to accuracy of the methods.

\subsection{Performance on the GPU}
\begin{figure*}[htb!]

	\centering
	\begin{subfigure}[b]{\textwidth}
	\includegraphics[width=0.33\linewidth]{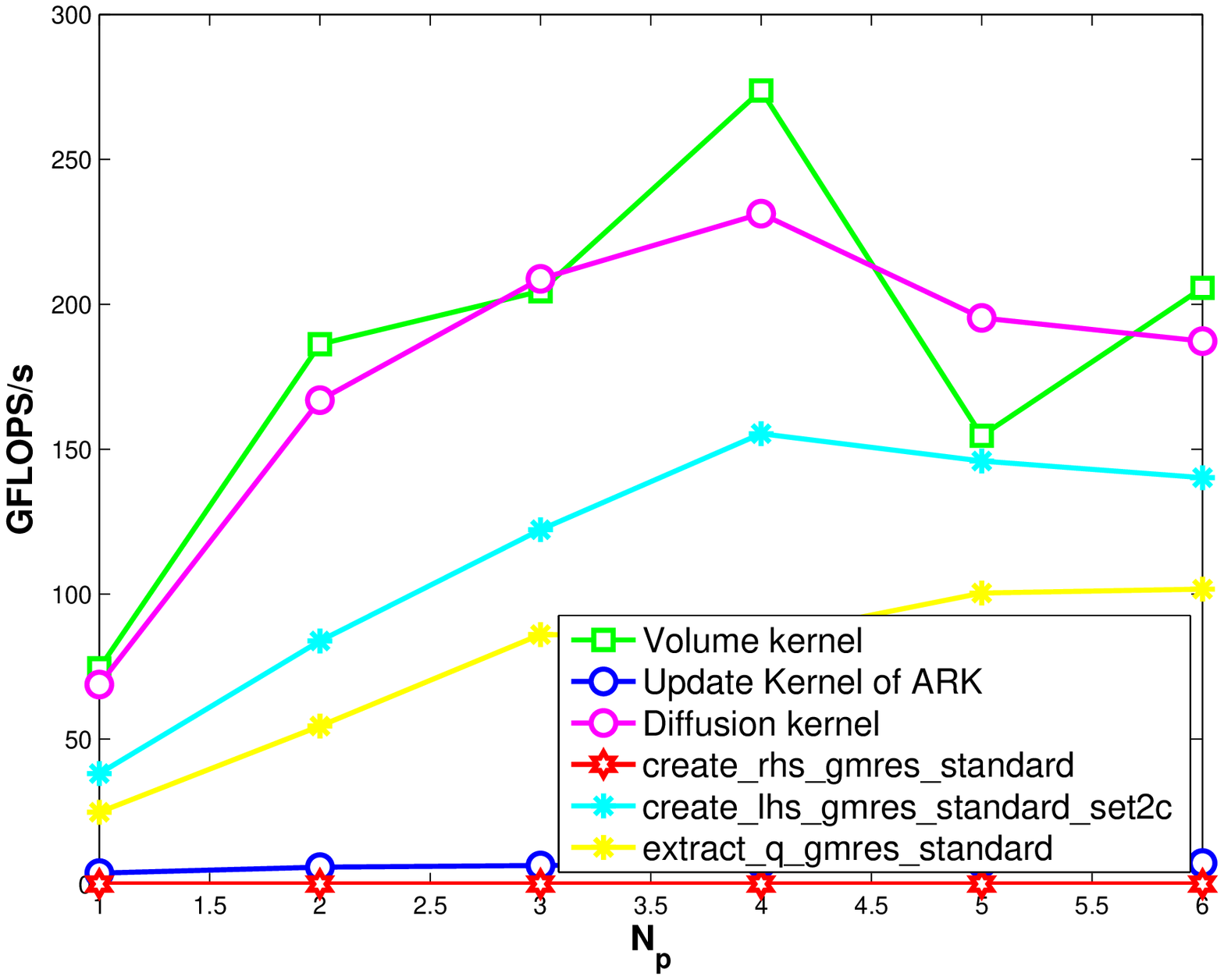}
	\includegraphics[width=0.33\linewidth]{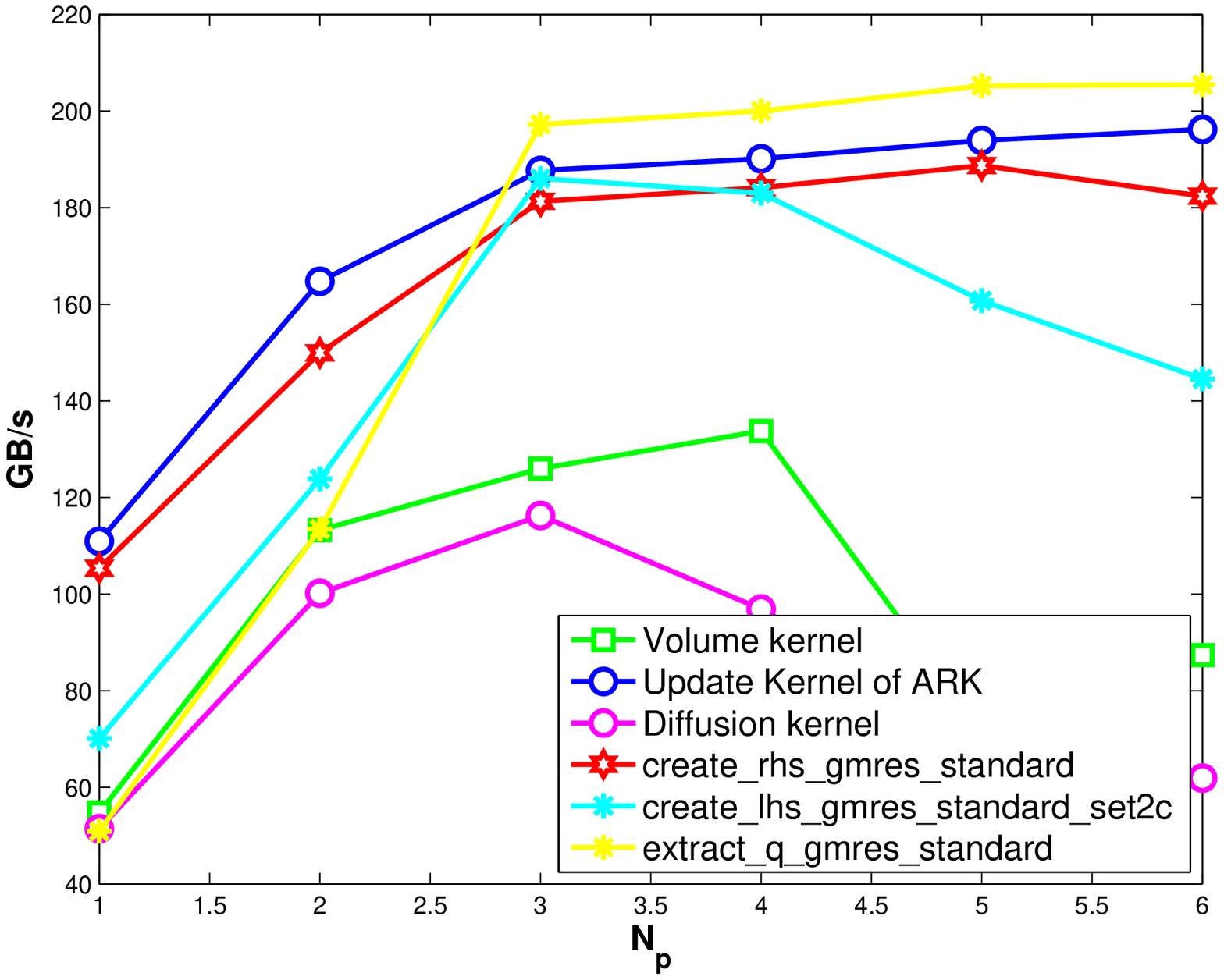}
	\includegraphics[width=0.33\linewidth]{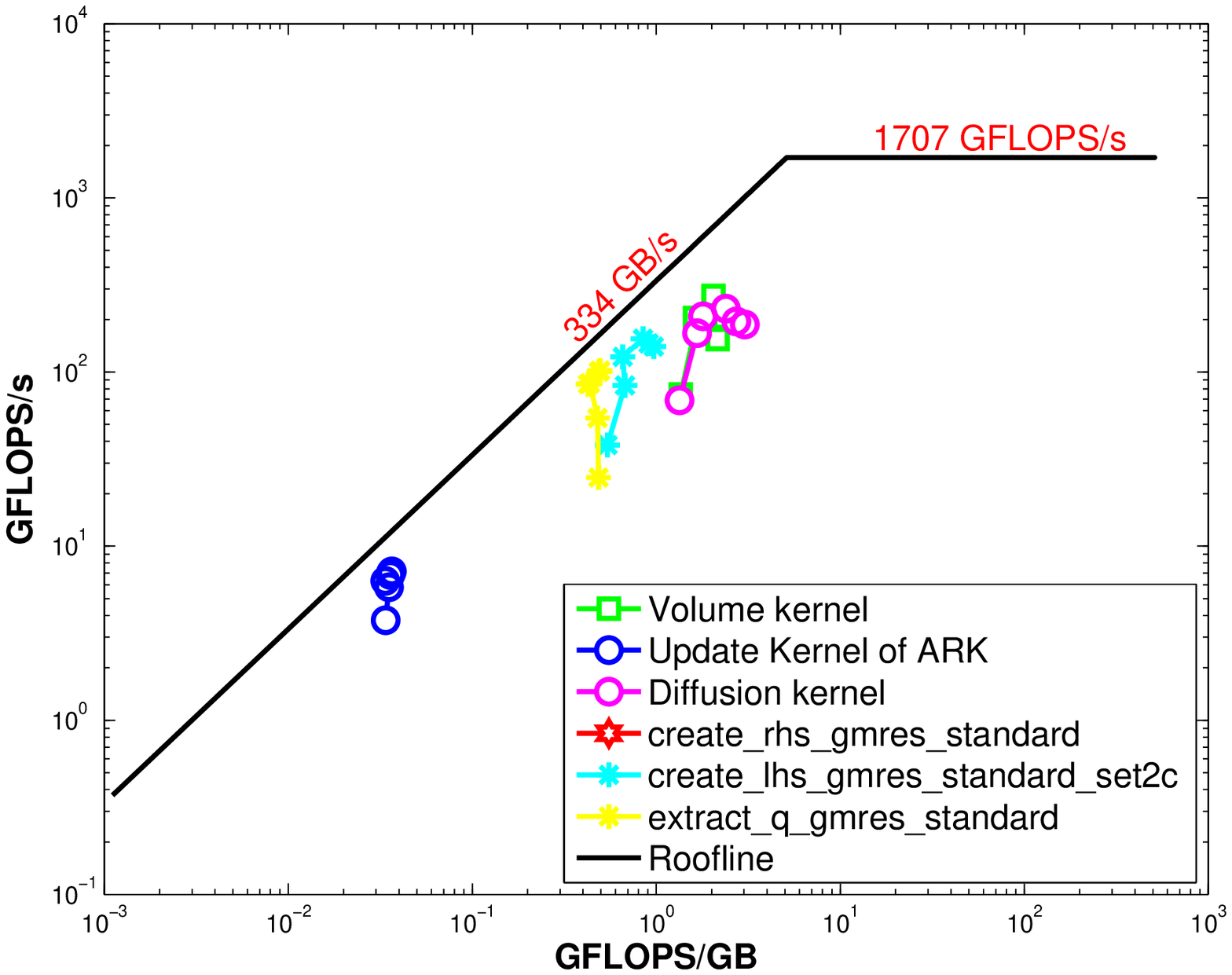}
	\caption{cG \nsr kernels performance}
	\label{gflopsCGn}
	\end{subfigure}
%
%
	\centering
	\begin{subfigure}[b]{\textwidth}
	\includegraphics[width=0.33\linewidth]{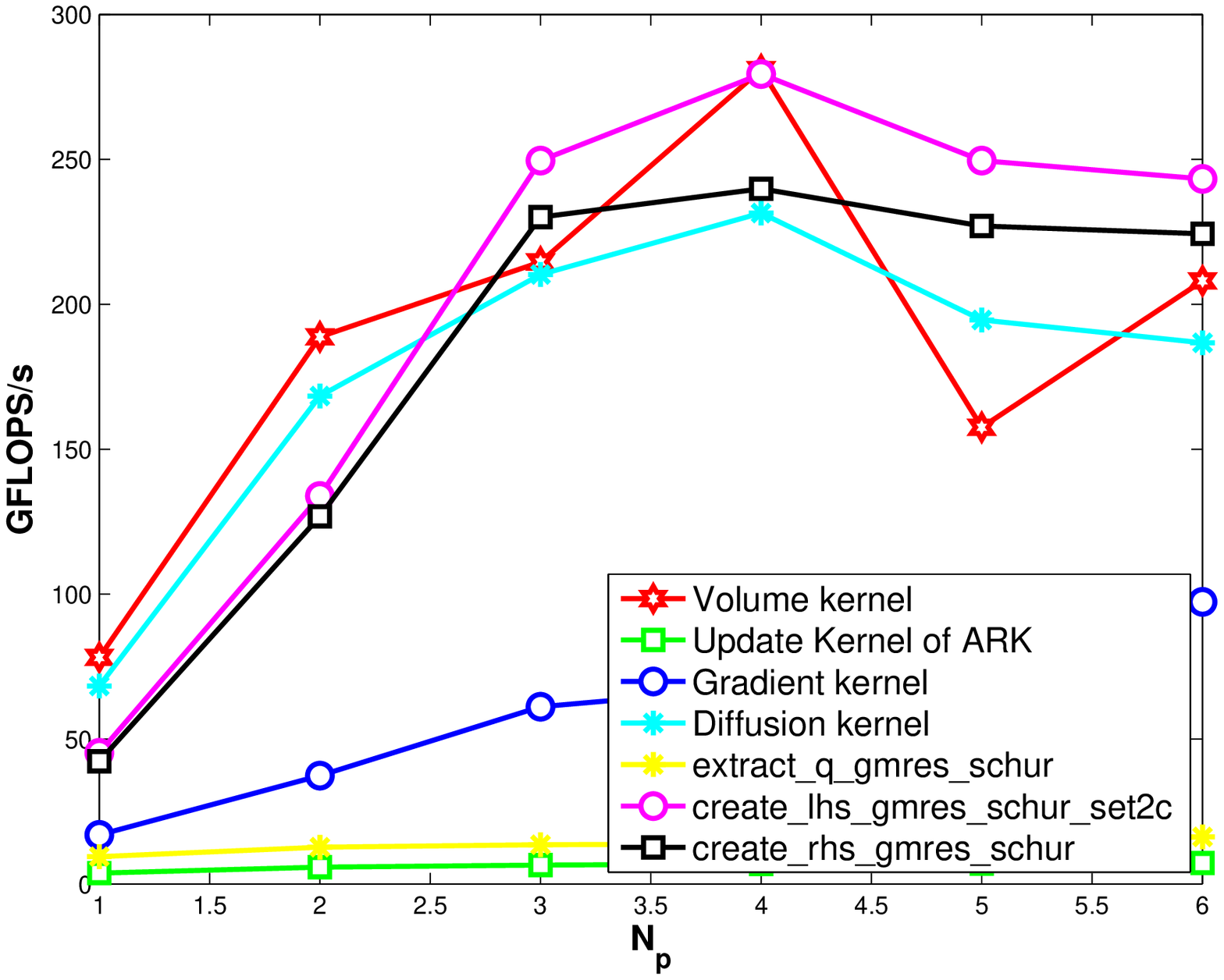} 
	\includegraphics[width=0.33\linewidth]{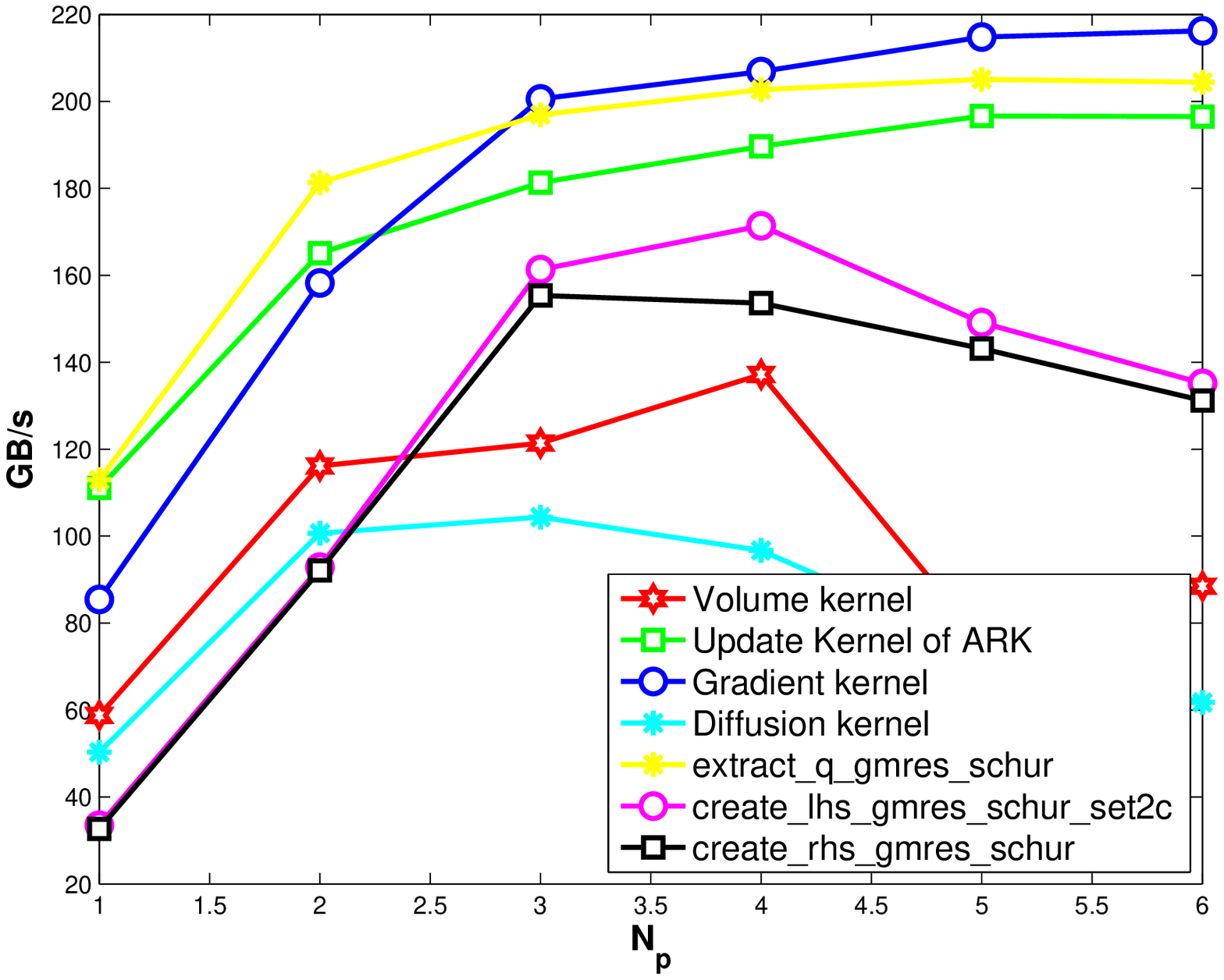}
	\includegraphics[width=0.33\linewidth]{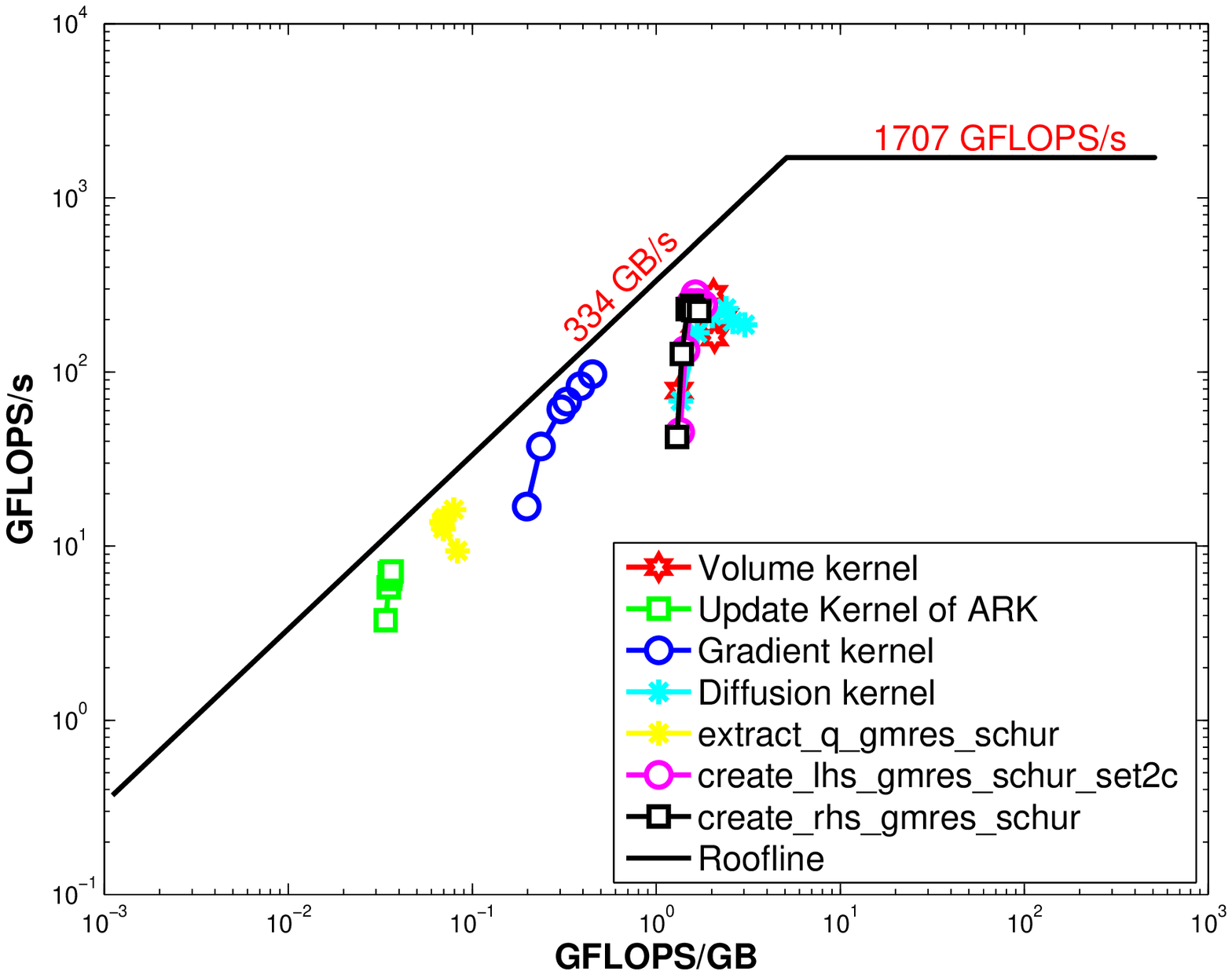}
	\caption{cG \sr form kernels performance}
	\label{gflopsCGs}
	\end{subfigure}
%
	\centering
	\begin{subfigure}[b]{\textwidth}
	\includegraphics[width=0.33\linewidth]{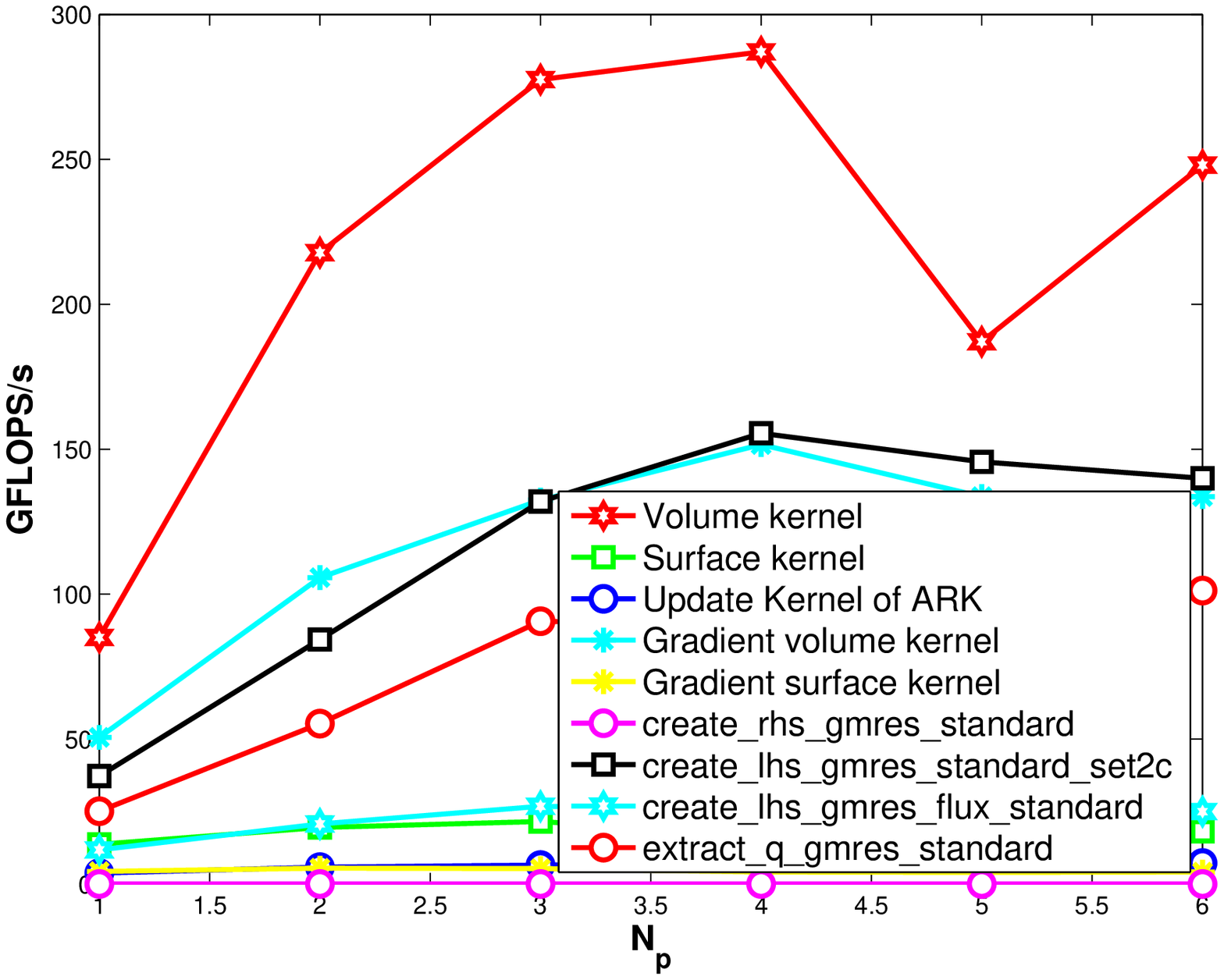} 
	\includegraphics[width=0.33\linewidth]{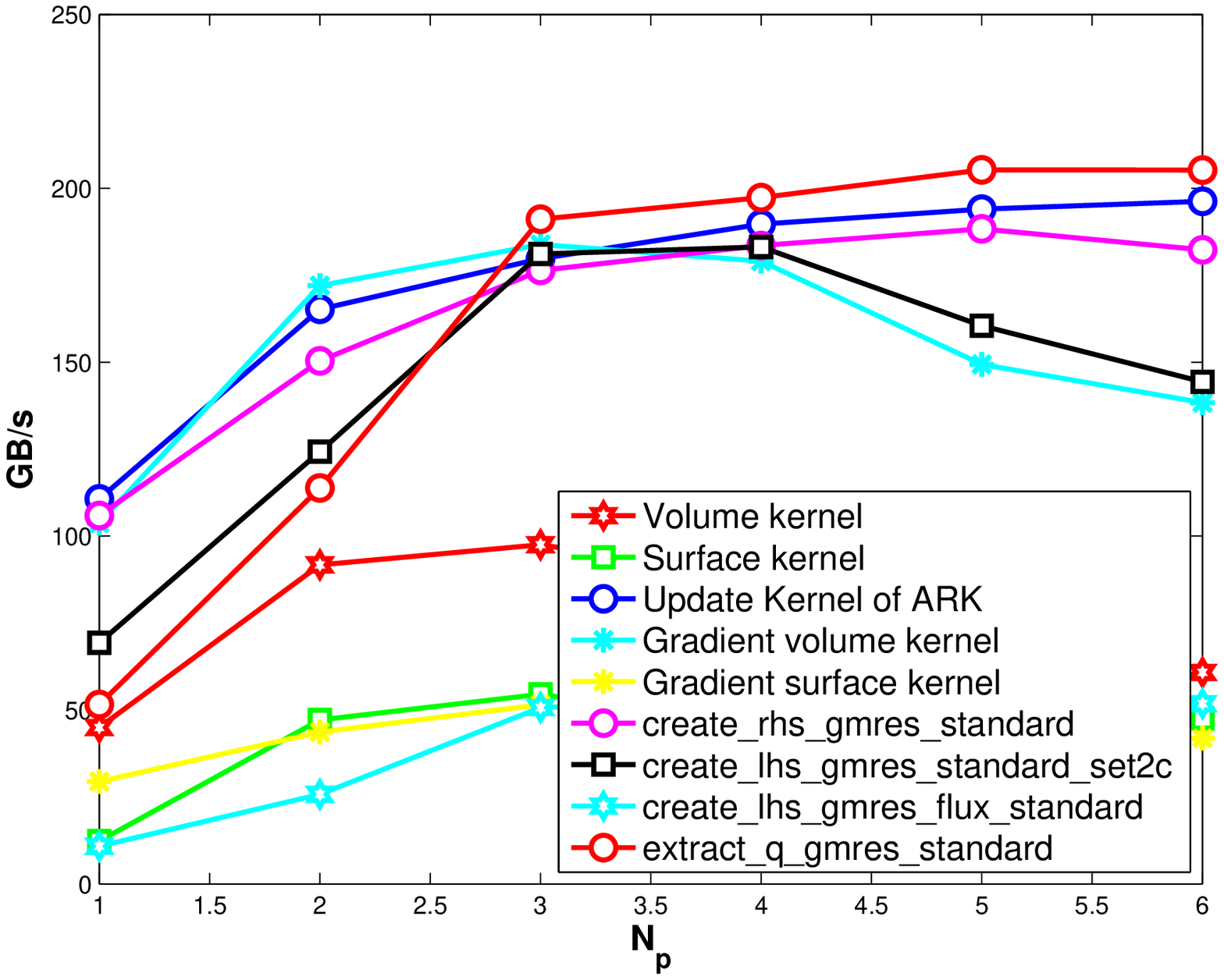}
	\includegraphics[width=0.33\linewidth]{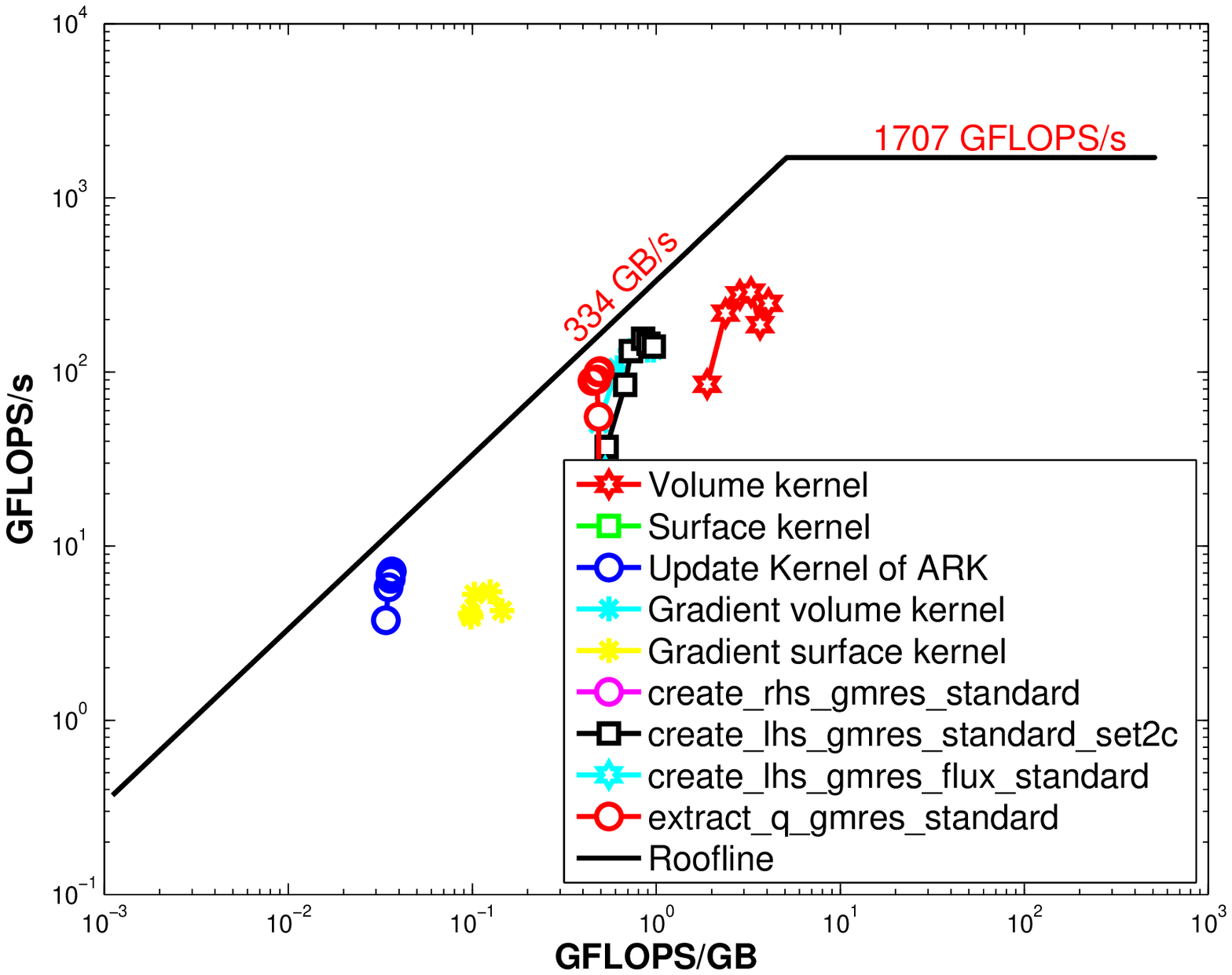}
	\caption{dG \nsr kernels performance}
	\label{gflopsDGn}
	\end{subfigure}

	\caption{Performance of IMEX kernels:  The measurements for this test are done using \textit{nvprof}: effective memory bandwidth = \textit{dram\_read\_throughput + draw\_write\_throughput}, and effective arithmetic throughput = \textit{flop\_dp/sp\_efficiency}. The Double Precision (DP) performance of the main kernels in cG and dG, \sr/\nsr forms of IMEX, are shown in terms of GFLOPS/s, GB/s and roofline plots to illustrate their efficiency. The GPU is a GTX Titan Black.
	}
\end{figure*}

The performance of the kernels involved in the IMEX time integration method is measured using the NVIDIA profiler counters \emph{dram\_throughput} and \emph{flop\_efficiency}. The results of these two metrics at different polynomial orders for the kernels involved in the IMEX time integration are shown in Figs.\ \ref{gflopsCGn}-\ref{gflopsDGn}. We also show roofline plots to easily identify whether a kernel is memory- or compute- bound. For the continuous Galerkin (cG) method, the IMEX volume kernels required for evaluating the left- and right- hand sides, for both \sr and \nsr forms, perform more or less the same as the standard right-hand side volume and diffusion kernels required for the explicit time integration step. The peak values are about 720 GFLOPS/s and 210 GB/s. The \sr forms are more compute intensive than the \nsr forms and this is revealed in the larger GFLOPS/s of the the right- and left- hand side IMEX volume kernels. The volume kernel for the discontinuous Galerkin (dG) discretization in Fig.\ \ref{gflopsDGn} shows a markedly superior GFLOPS/s performance that is better than the other kernels for two reasons: first, because the diffusion kernel is split into two in the LDG scheme, second, we have not formulated the dG IMEX \sr form that would have been more compute intensive than the \nsr form.  The right-hand side volume kernels for the \nsr IMEX forms for both cG and dG have very low GFLOPS/s as expected because the kernels do not perform many calculations. The explicit time update kernel of the ARK method has the highest bandwidth usage; most other kernels other than the volume and diffusion kernels also show very good bandwidth usage. Our kernels are mostly memory-bound but the compute-intensive volume kernels are very close to being compute-bound.

\begin{figure*}
        \centering
	\includegraphics[width=0.36\linewidth]{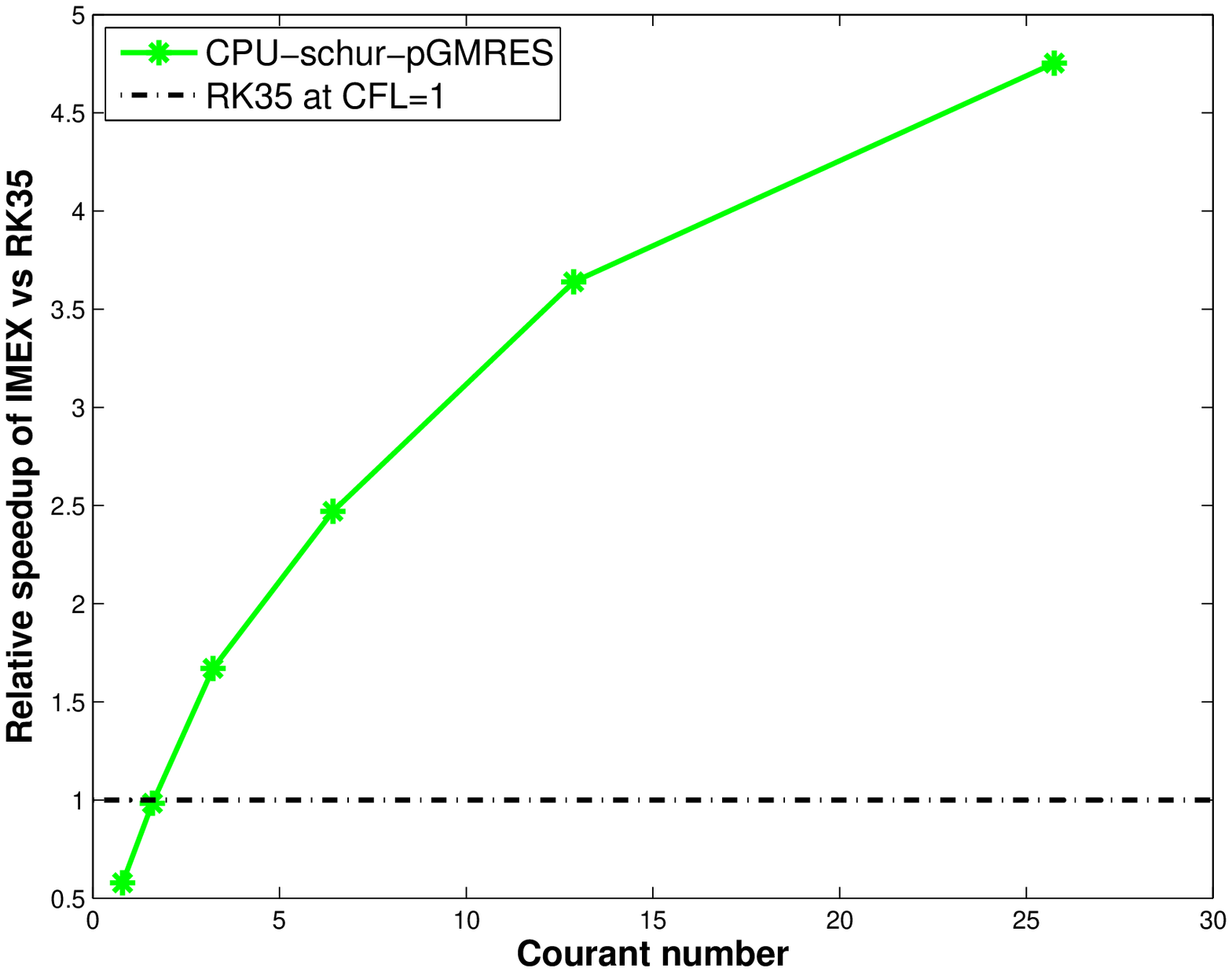} 
	\includegraphics[width=0.36\linewidth]{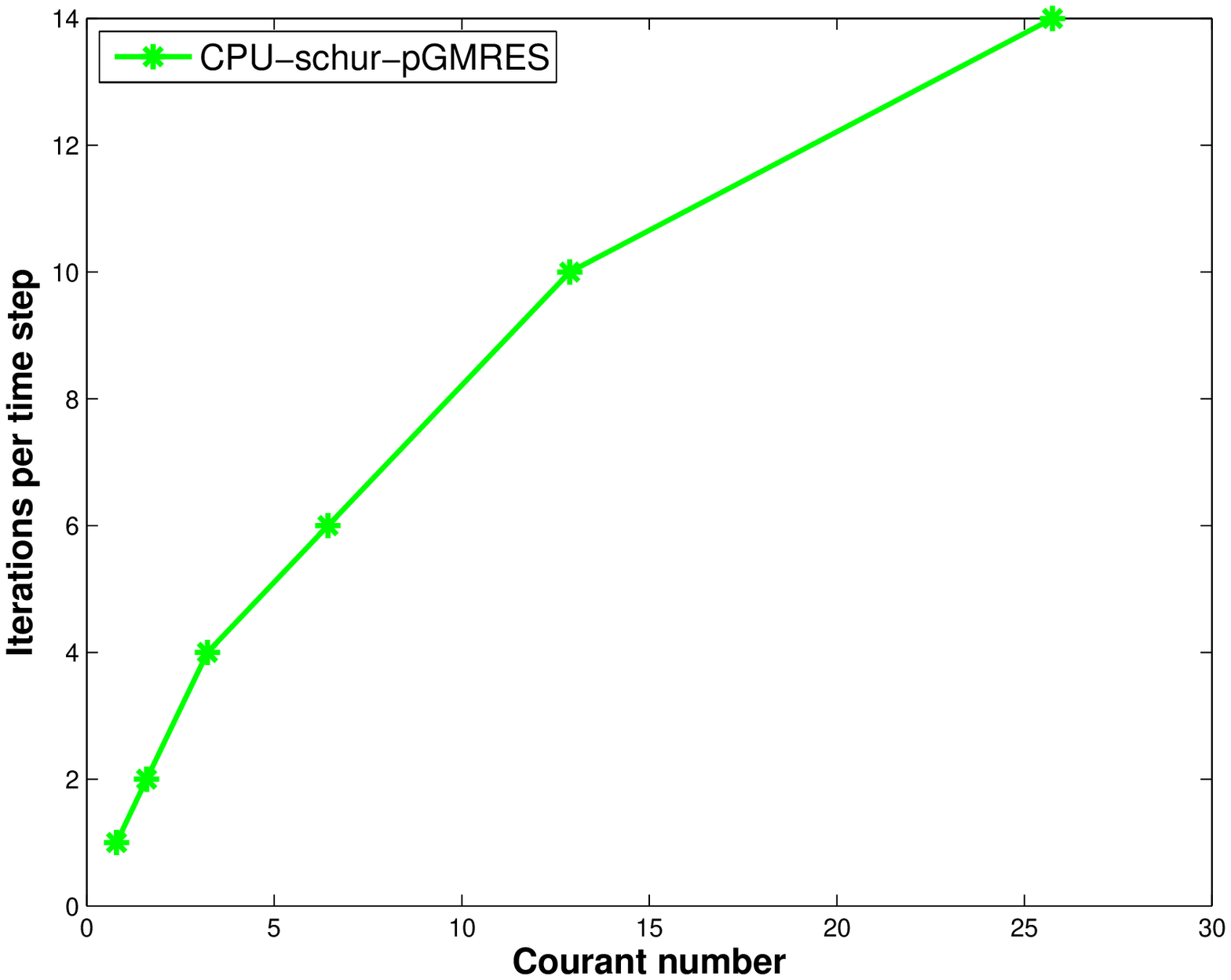} 
	\caption{CPU results: (left) Speedup of IMEX over RK35 at Courant number = 1 is shown for the \sr form using preconditioned GMRES iterative solver. The 2D rising thermal bubble problem is solved using a grid of 10 x 10 elements with polynomial order of 4. A maximum speedup of about 5X is observed for Courant numbers greater than 15. (right) The number of GMRES iterations required are shown.} 
	\label{speedupimex-cpu}
%
        \centering
	\includegraphics[width=0.32\linewidth]{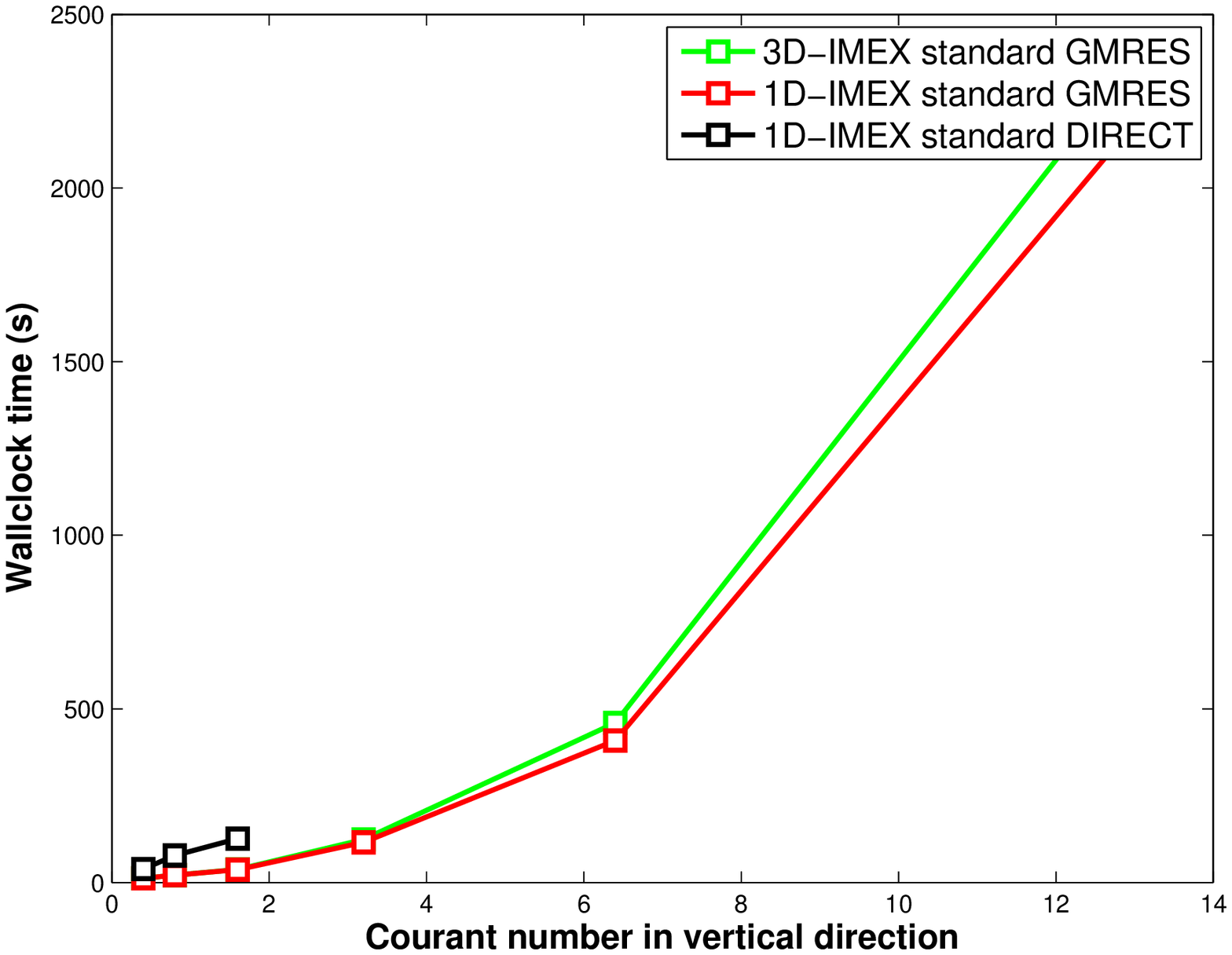} 
	\includegraphics[width=0.32\linewidth]{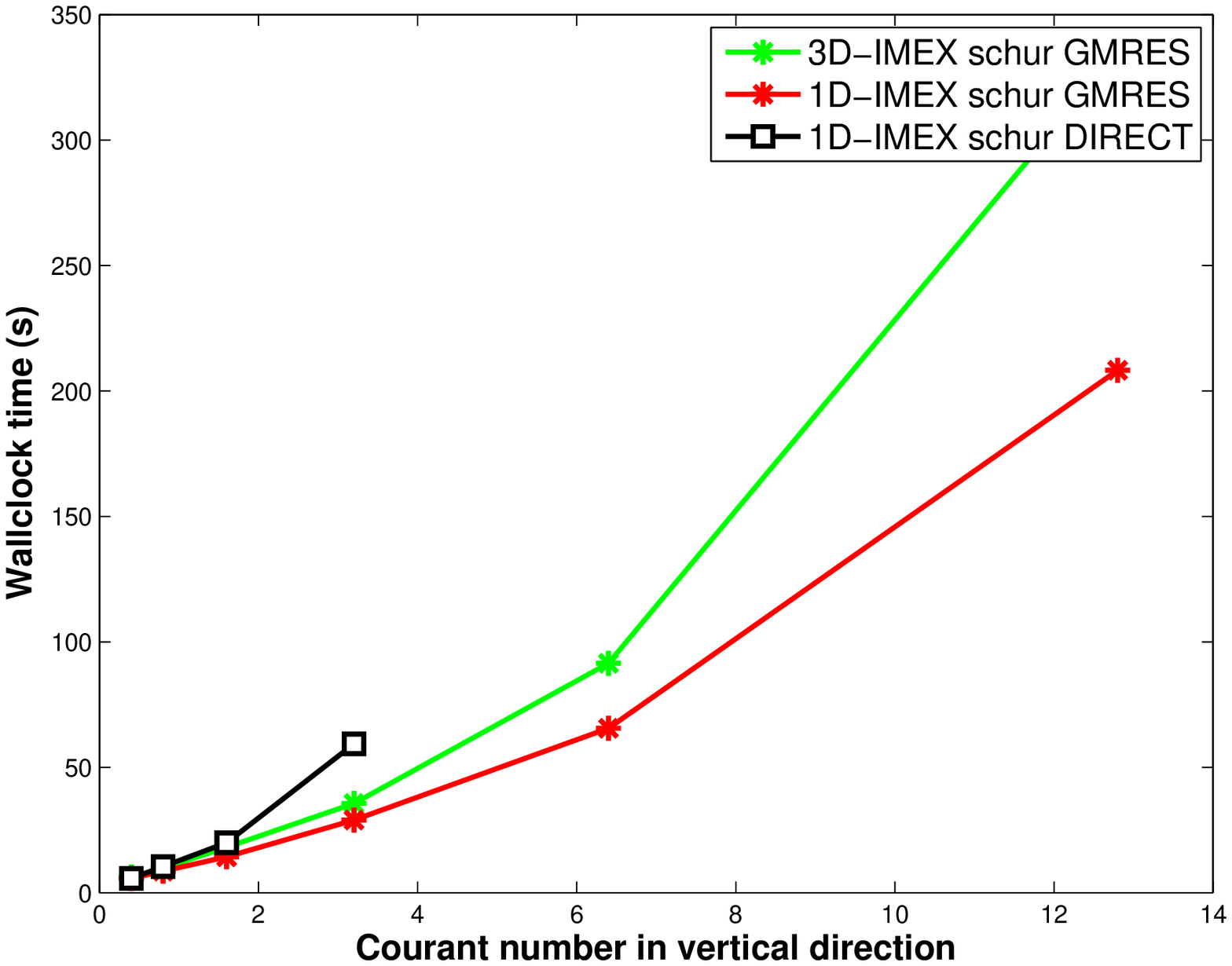} 
	\includegraphics[width=0.32\linewidth]{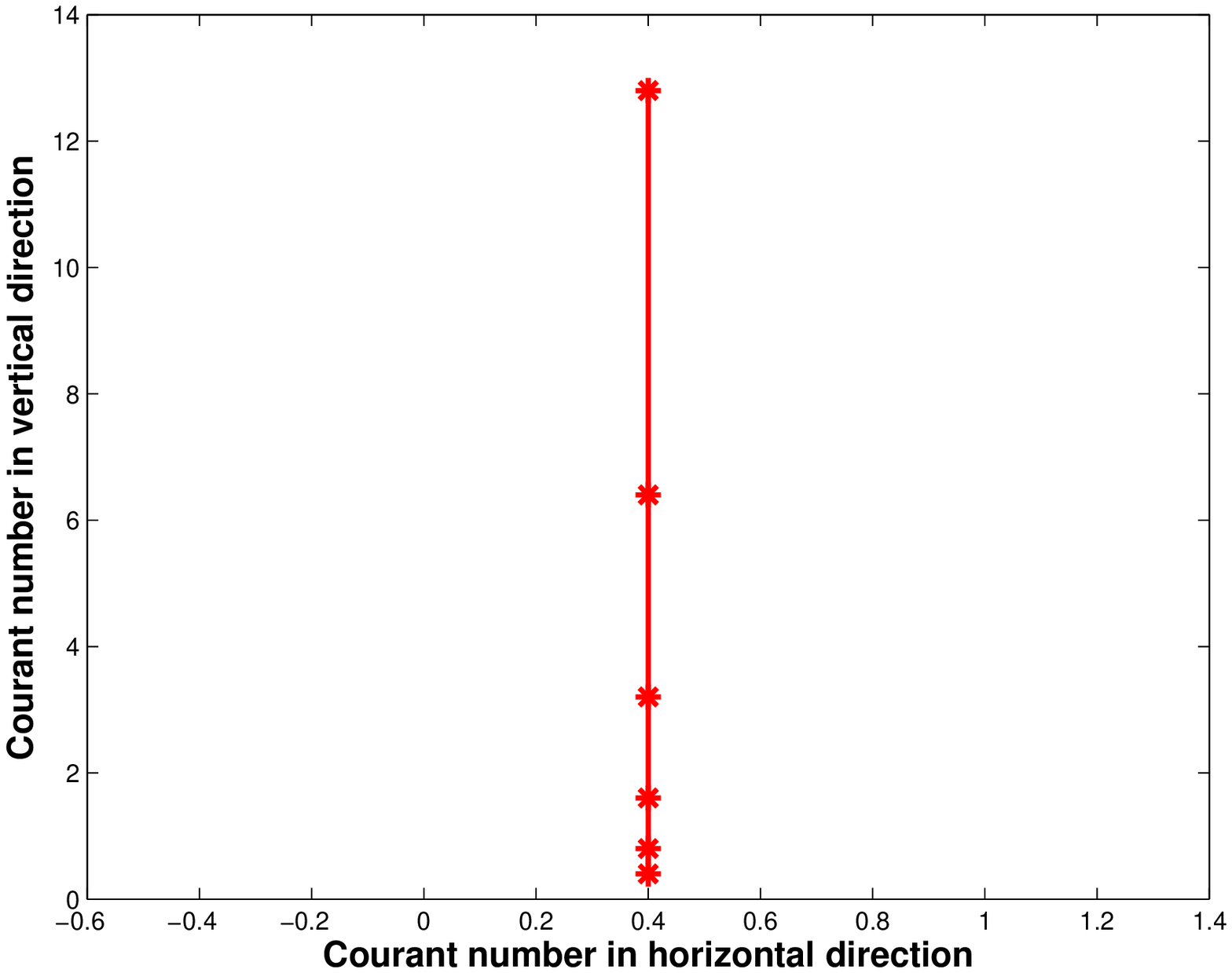} 
	\includegraphics[width=0.32\linewidth]{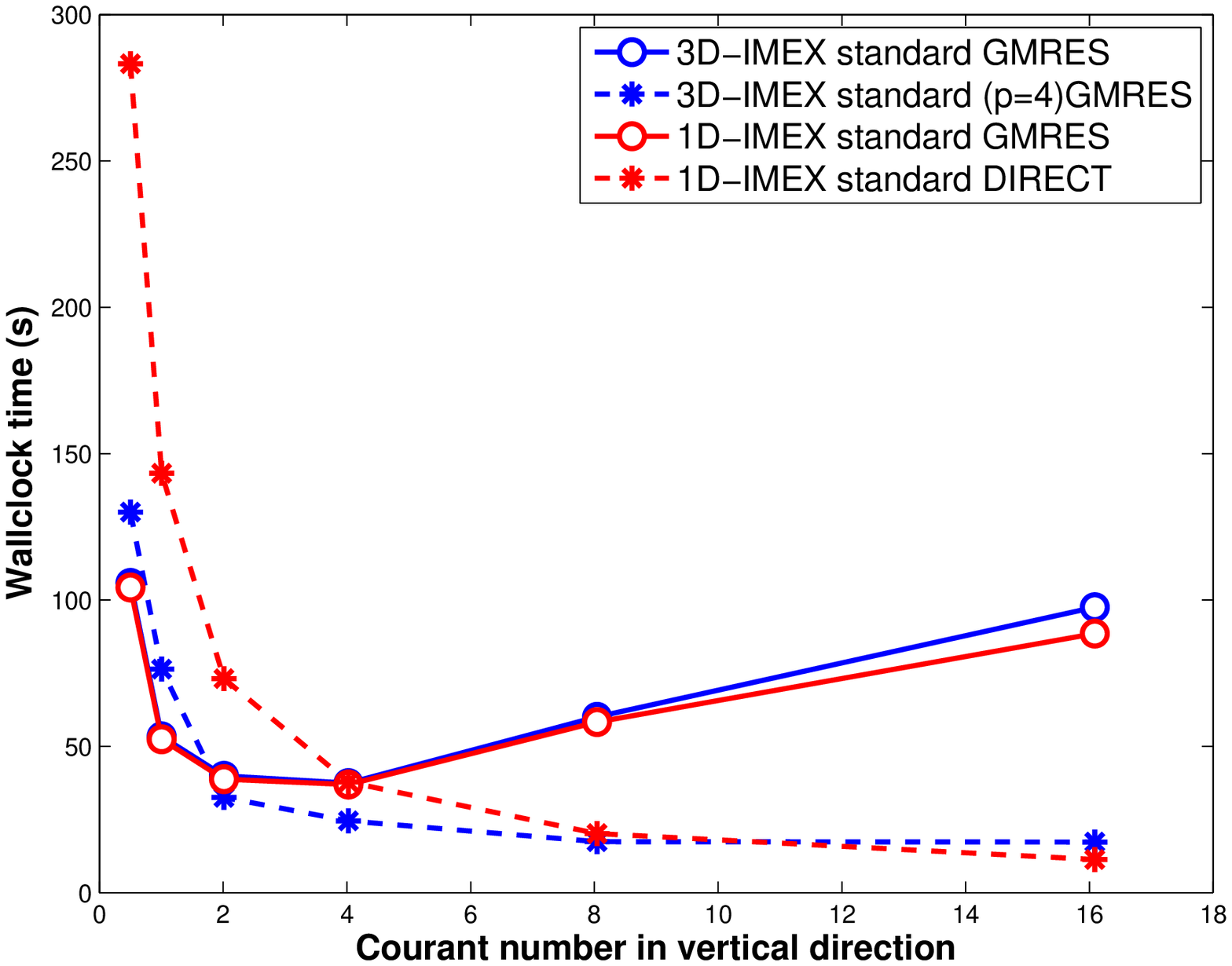} 
	\includegraphics[width=0.32\linewidth]{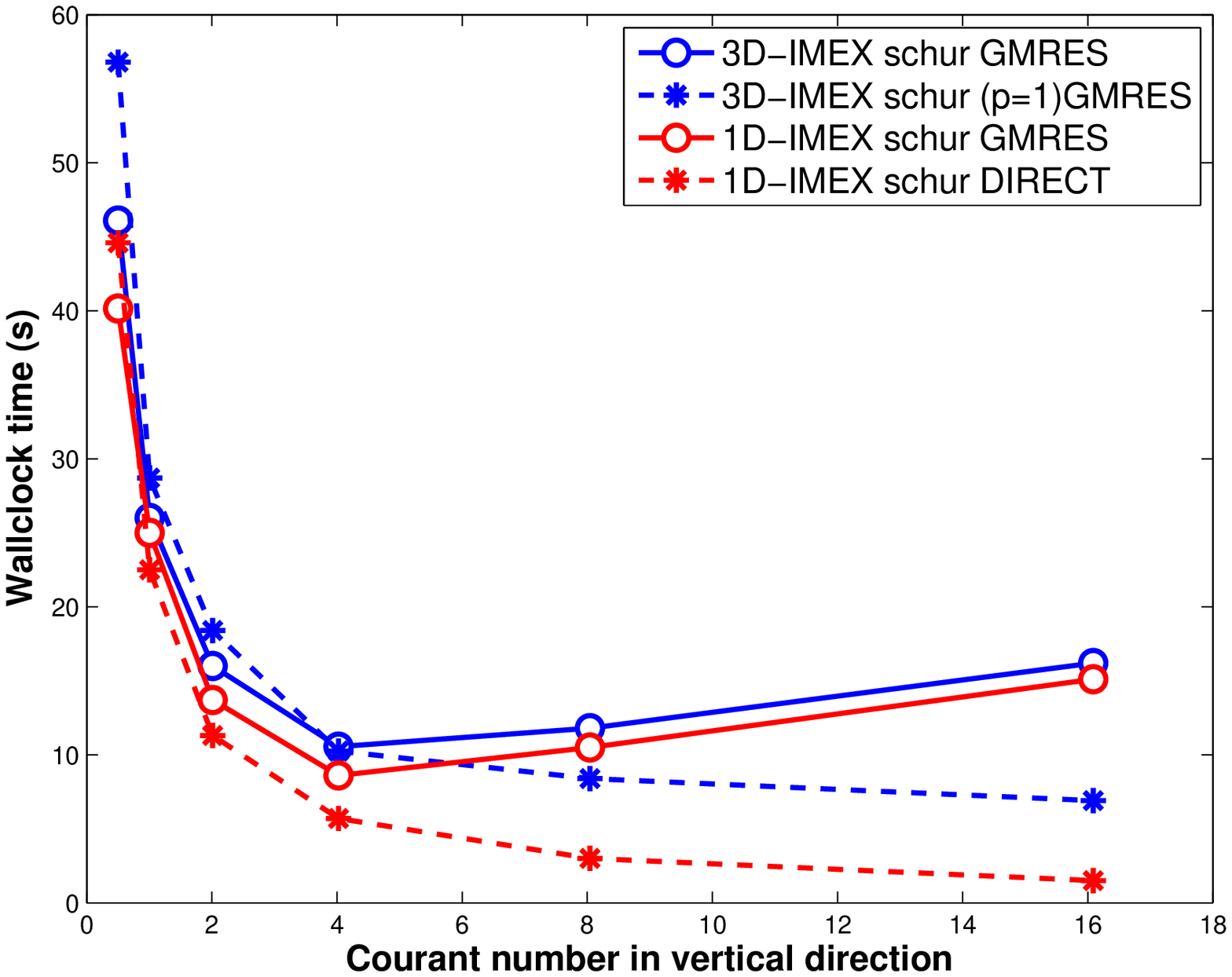} 
        \includegraphics[width=0.32\linewidth]{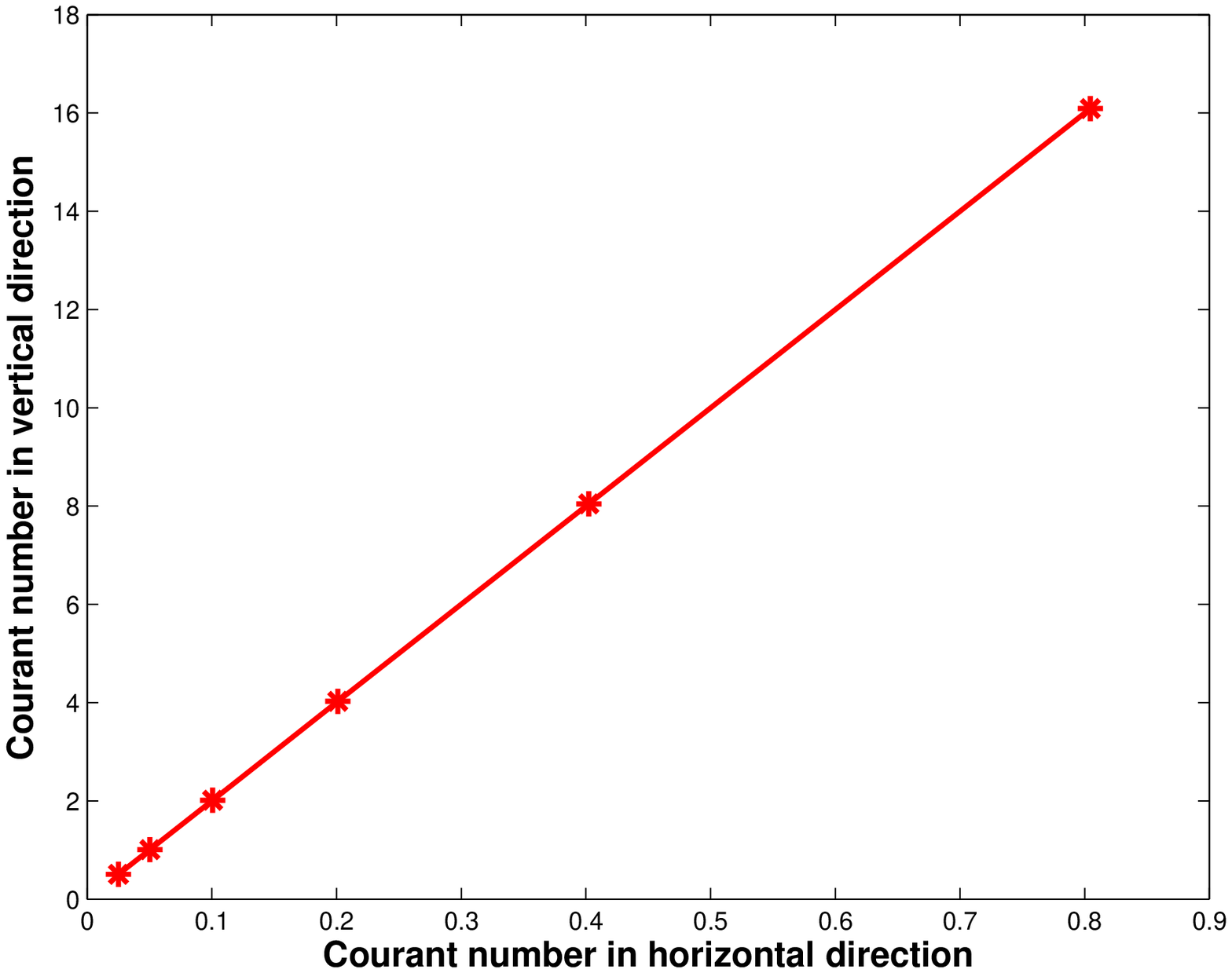} 
	\caption{Comparison of 1D- and 3D- IMEX using  GMRES iterations for the \nsr and \sr forms. The test case is a 2D rising thermal bubble problem. For the top plots, the grid in the vertical direction
	is refined to get different Courant numbers in the vertical direction while using the same time step; a constant $C_H$ = 0.4 is used in the horizontal direction. The bottom plots are produced by a more realistic approach in practical NWP -- changing the time step while
	using the same grid, hence, $C_H \le $ 1 varies for each simulation.}
	\label{comparison-1d3d}
\end{figure*}

The most important question we are interested in answering in this work is whether  IMEX time integration is worth the effort on the GPU, i.e., if it will be faster than using explicit time integrators. On the CPU,  NUMA can benefit from a relative speedup of 5X using IMEX time-integrators compared to the explicit RK35 time integrator when using the \sr forms of IMEX; the \nsr forms often do not give any speedup over an explicit time integrator.  We reproduce an example run on the CPU to demonstrate this advantage in Fig.\ \ref{speedupimex-cpu}, where we get a speedup of about 4.7X.  However, it is not clear whether NUMA could benefit equally well on the GPU because the most efficient iterative solvers (GMRES and BiCGstab) require dot-products that are very slow to compute on the GPU. This has also been the case on a multi-node CPU cluster for similar reasons.

\begin{figure*}
        \centering
	\includegraphics[width=0.4\linewidth]{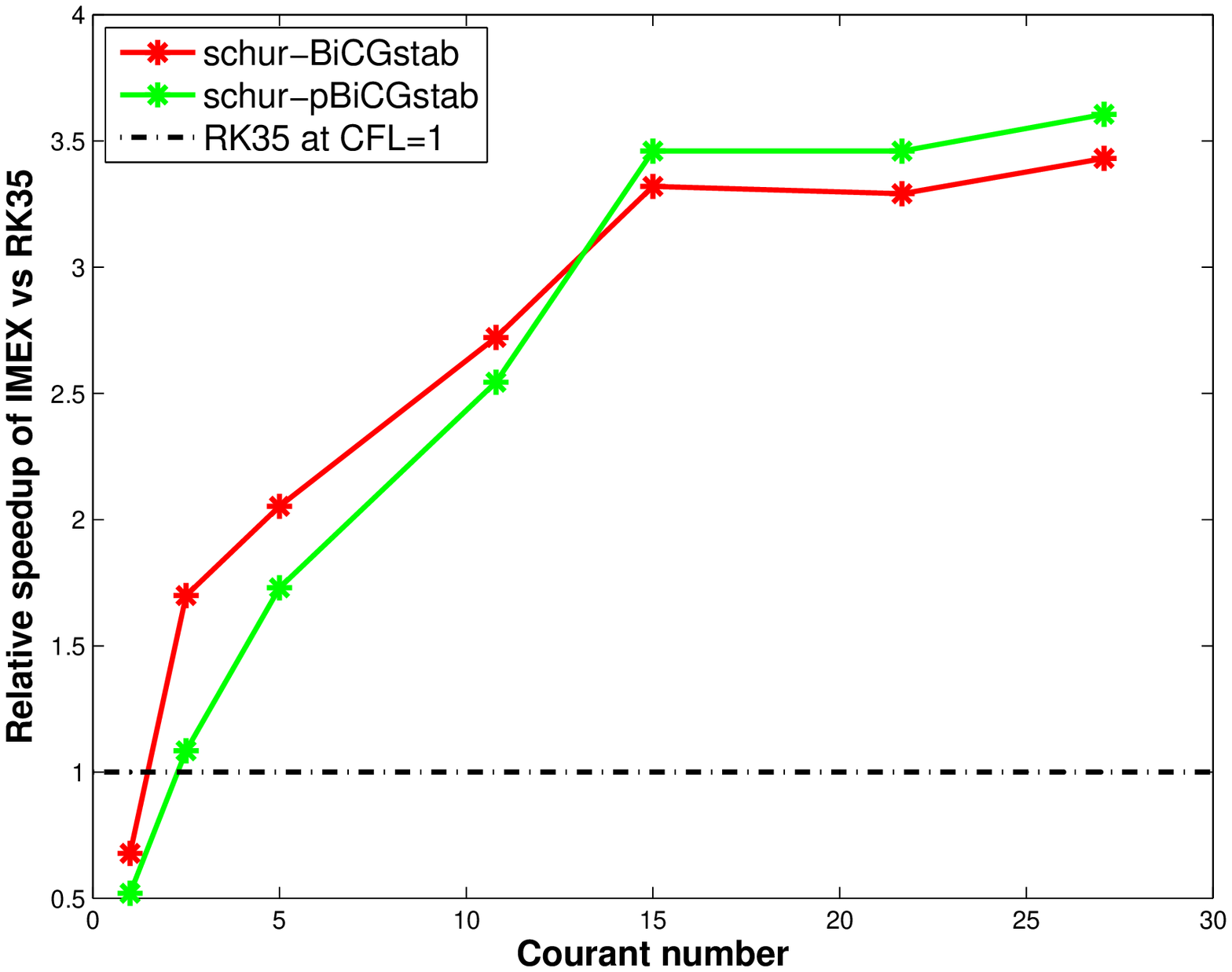} 
	\includegraphics[width=0.4\linewidth]{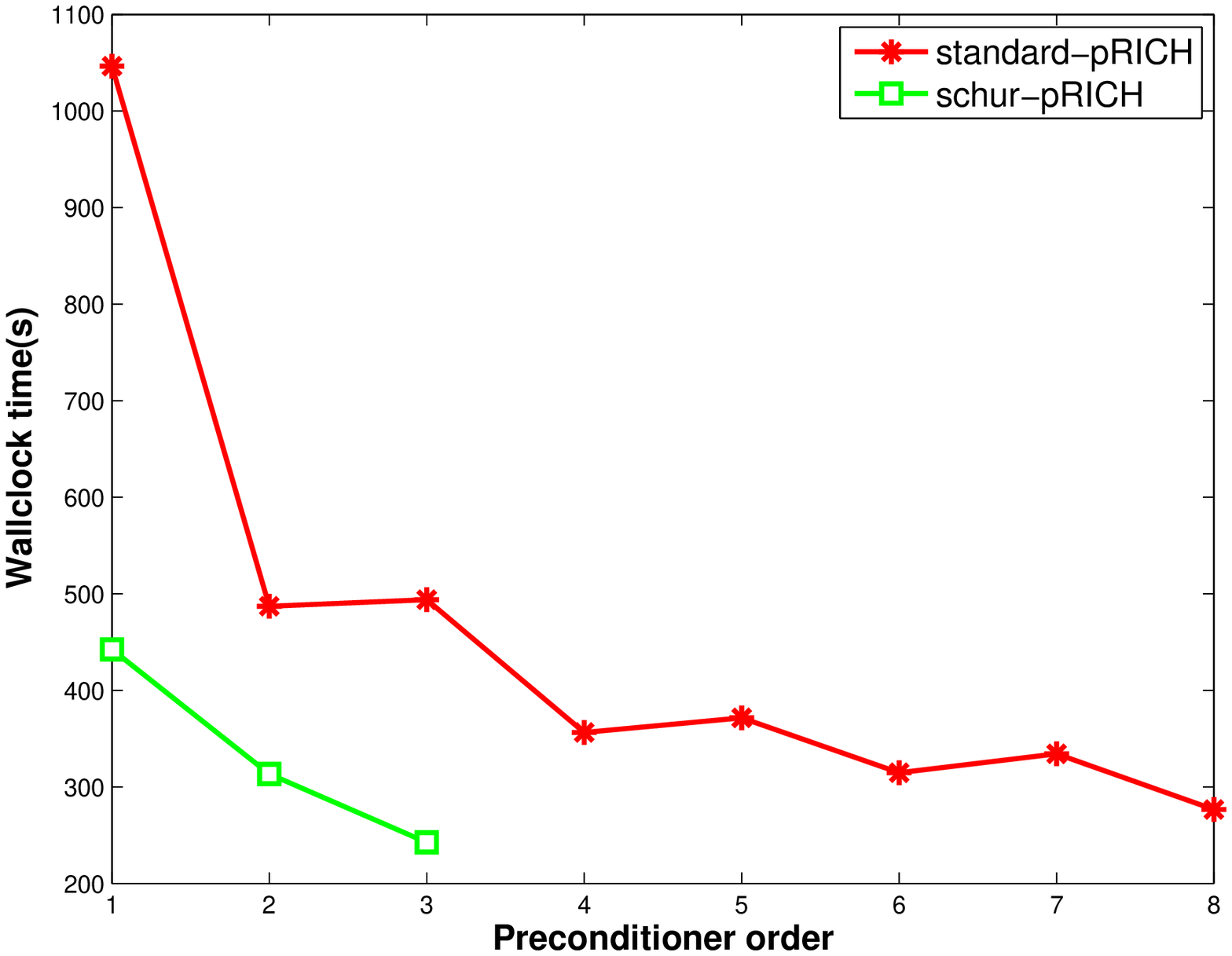} 
	\caption{(left) Speedup of IMEX over RK35 at Courant number = 1 is shown for the \sr form using un-preconditioned and preconditioned BiCGstab iterative solver. The 2D rising thermal bubble problem is solved using a grid of 20 x 20 elements with polynomial order of 7. A maximum speedup of about 3.5X is observed for Courant numbers greater than 15. (right) The effect of the PBNO preconditioner on Richardson iteration is shown at different polynomial order of the preconditioner.} 
	\label{speedupimex}
%
        \centering
	\includegraphics[width=0.4\linewidth]{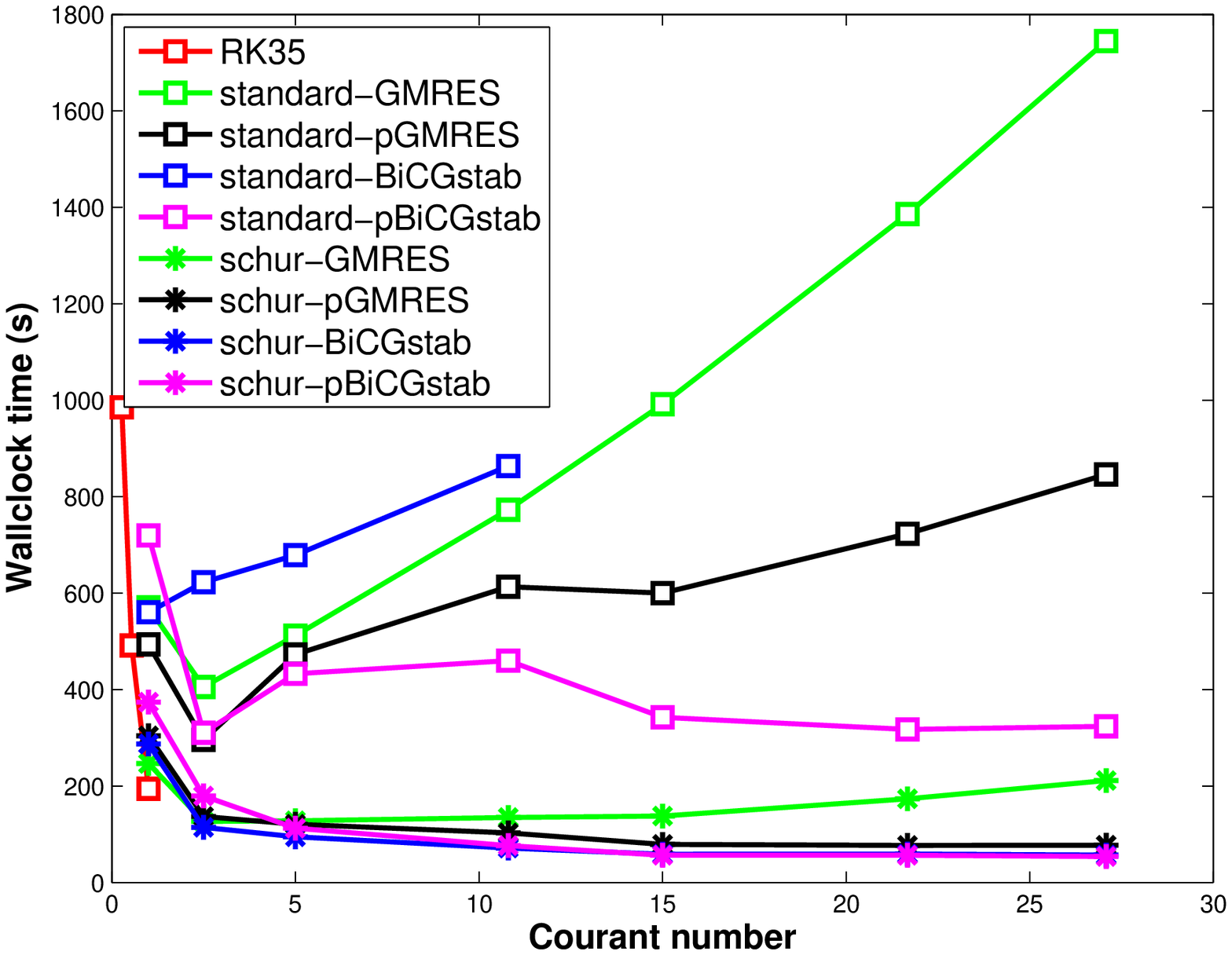} 
	\includegraphics[width=0.4\linewidth]{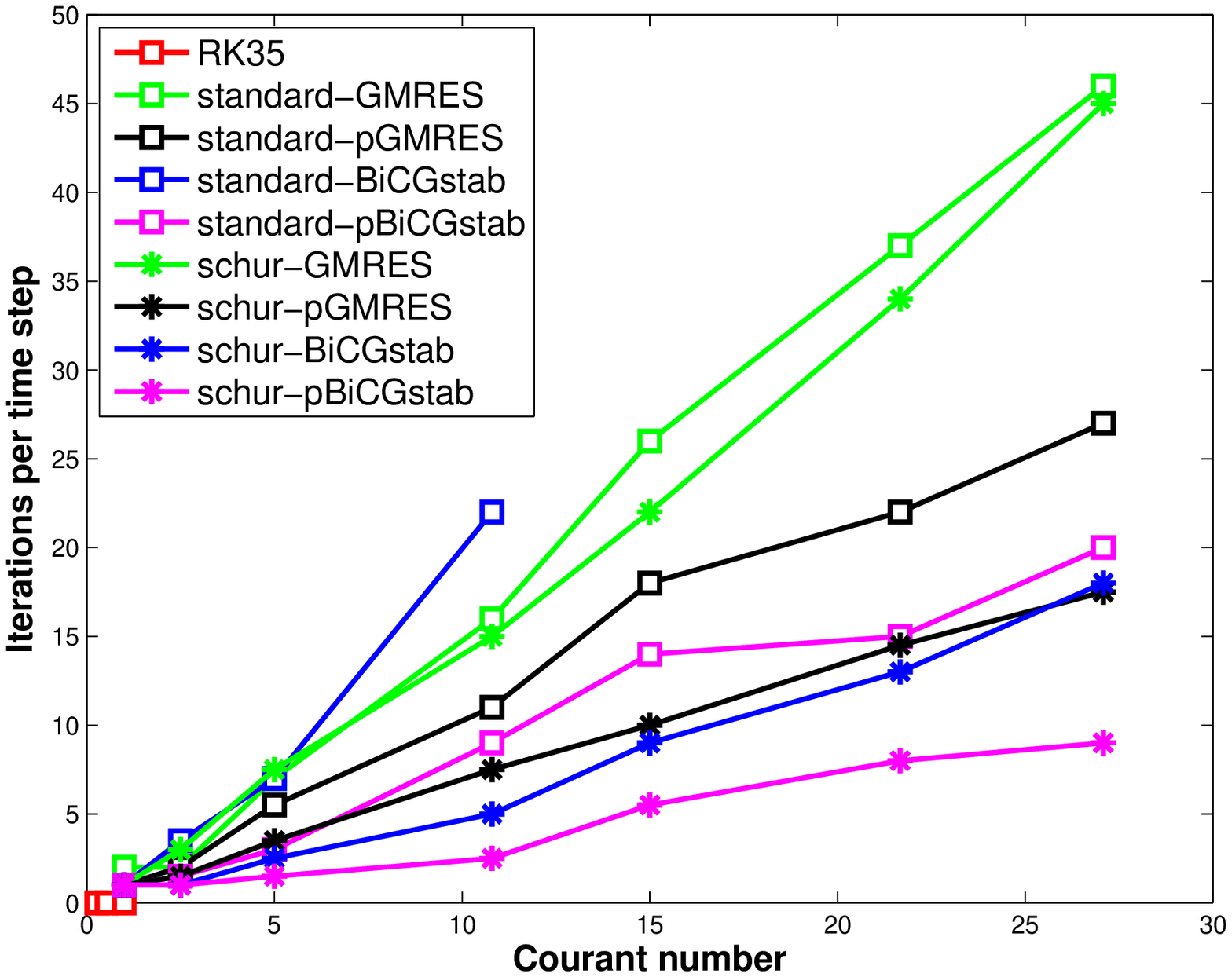} 
	\caption{The \nsr and \sr forms are compared using GMRES and BiCGstab iterative solvers. The rising thermal bubble problem is solved using a grid  of 20 x 20 elements at polynomial order 7. The order of the preconditioner is p=1. The \nsr forms are much slower than the \sr forms especially at higher Courant numbers; however, the \nsr forms can significantly benefit from increasing the preconditioner order.} 
	\label{comparisonimex}
%
        \centering
	\includegraphics[width=0.4\linewidth]{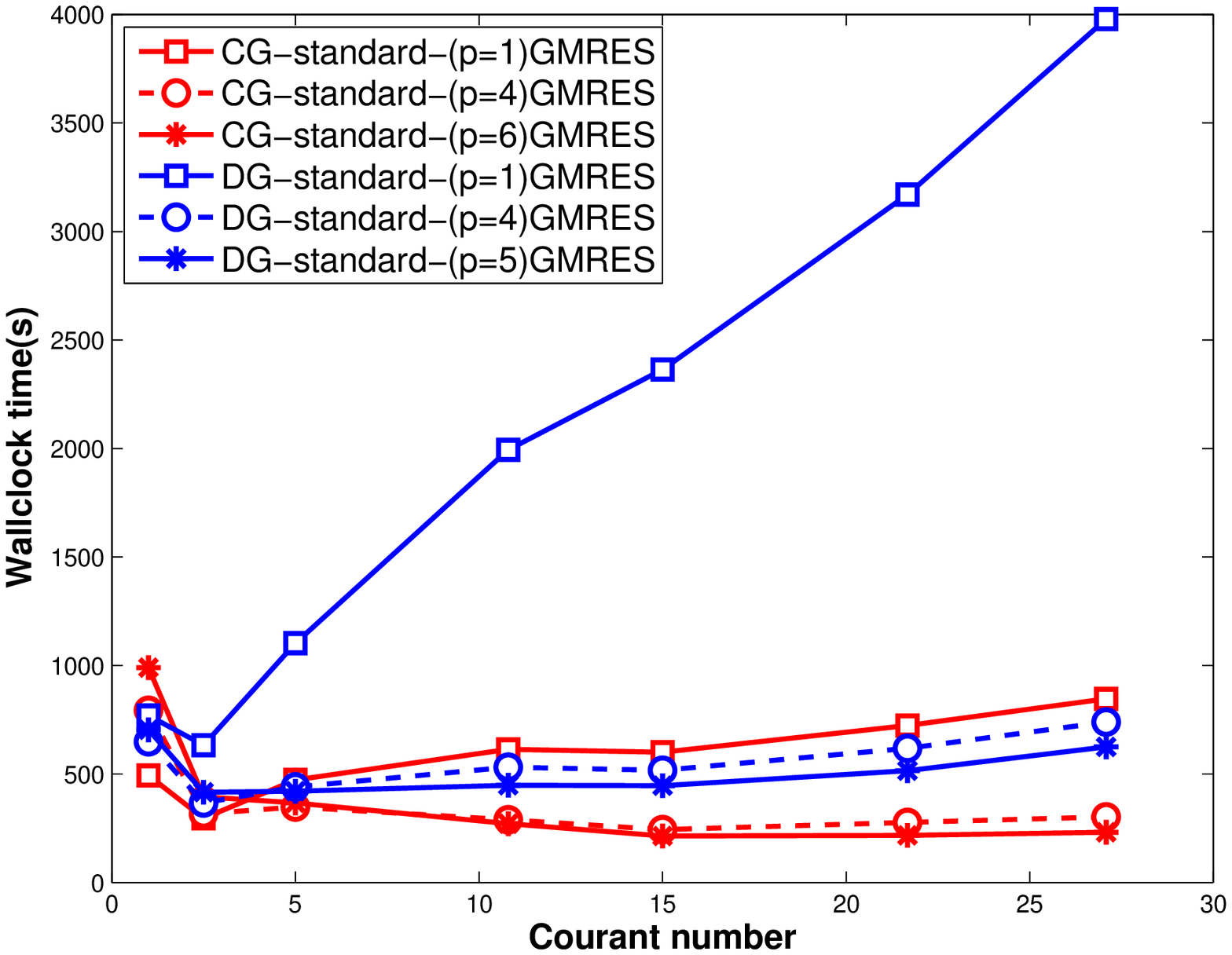} 
	\includegraphics[width=0.4\linewidth]{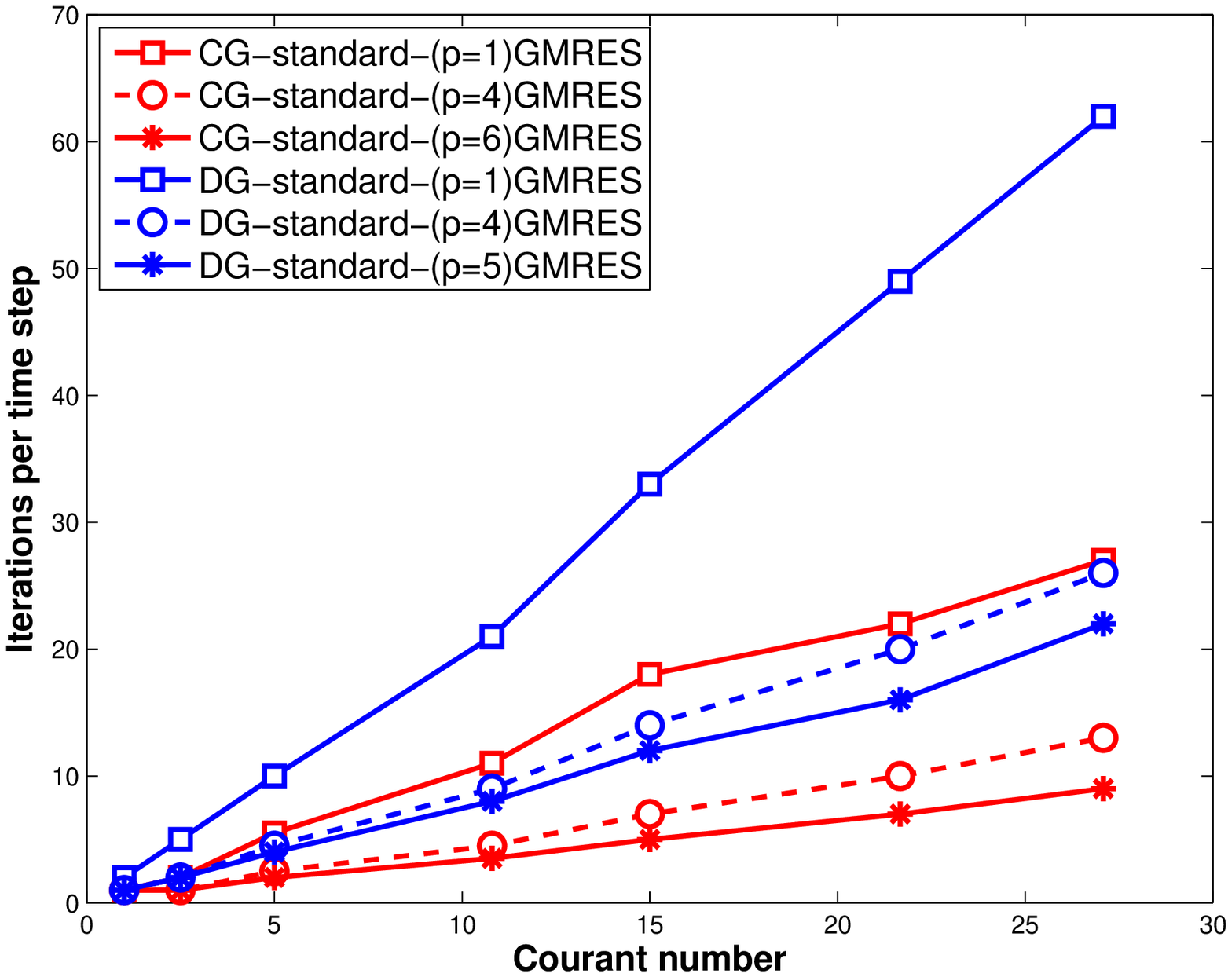} 
	\caption{The \nsr form cG and dG discretizations are compared. With the same preconditioner order, say p=1, dG requires more iterations than cG  which results in larger wall clock times. Increasing the preconditioner order to p=4 improves the performance of dG significantly and brings it closer to cG performance. The rising thermal bubble problem is solved using a grid  of 20 x 20 elements at polynomial order 7.} 
	\label{comparison-cgdg}
\end{figure*}

In Fig.\ \ref{speedupimex}, we show the relative speedup result of IMEX over RK35 on the GPU using the \sr form and preconditioned BiCGstab iterative solver. We get a speedup of more than 3.5X at Courant numbers greater than 15. For this particular test case, a 2D rising thermal bubble problem, IMEX also beats RK35 starting from Courant numbers as lows as 2. It seems that our concern regarding inefficiency of dot-products on GPUs did not  significantly impact our results. On the right, we show results of the dot-product free Richardson iteration using different polynomial order of the PBNO preconditioner. For the Richardson method to be competitive with the Krylov solvers, one should use a high-order polynomial preconditioner. The results show that Richardson iteration is unable to match the results of the Krylov solvers with the \sr form even though it is dot-product free and that we used orders of up to 7 for the preconditioner. However, it yielded more or less similar results as the Krylov solvers for the \nsr form in which none of the solvers were able to beat RK35.

In Fig.\ \ref{comparisonimex}, we compare the \nsr and \sr forms, GMRES and BiCGstab iterative solvers, and the effect of preconditioning. We clearly see that the \nsr forms are much slower than the \sr forms in all cases especially at higher Courant numbers. None of the \nsr form runs were able to beat RK35 which is run with the maximum allowable Courant number of 1; this  result highlights the value of formulating the \sr form. The \sr forms were able to give speedup over RK35 starting from Courant numbers as low as 2. In general, BiCGstab seems to perform better than GMRES at least for this particular test case; however, the un-preconditioned BiCGstab  was not able to converge for the \nsr forms unlike  GMRES. The wallclock times of the solvers in general follow the trend of the average number of iterations required per time step, shown on the right of Fig.\ \ref{speedupimex}. One downside of GMRES on the GPU is that it is much more memory intensive than BiCGstab. The memory required by GMRES for storing intermediate Krylov vectors grows with the maximum allowable number of iterations, while for BiCGstab it is constant. Moreover, the number of dot-products required by GMRES grows quadratically with the number of iterations $\mathcal{O}(N^2/2 + 3N/2)$ while that of BiCGstab grows linearly with number of iterations $\mathcal{O}(5N)$. For these  reasons, BiCGstab is preferable on the GPU as long as convergence can be ensured. 

In Fig.\ \ref{comparison-cgdg}, we compare the \nsr forms for cG and dG discretizations. We can see that the number of iterations per time step increases at a much faster rate for dG than cG when the same level of preconditioning is applied (p=1); therefore, the wallclock time of the dG solution is much larger than that of cG at p=1. This is likely due to the weaker coupling between elements in dG than cG; however, using a higher order preconditioner p=4, dG was able to get close to the performance of cG. The cG \nsr form also benefits from higher order preconditioning but not as much as the dG \nsr form. Both \nsr forms are slower than RK35, hence, this result again highlights the value of formulating the \sr form for dG. Unfortunately, our attempt at formulating the \sr form dG with centered fluxes does not converge but other efforts in this direction are underway.

\begin{table*}
\centering
\caption{Memory usage of the IMEX methods measured on:  1) a 2D rising thermal bubble of 10x10 elements at N=4;  2) acoustic wave problem on a cubed sphere grid of 6x10x10x3 elements at N=3. In both cases, the \nsr forms show a much larger memory usage than the \sr forms as expected. For 3D-IMEX, GMRES iterations need space for storing intermediate Krylov vectors  which makes it more expensive than BiCGstab. The direct solution of 1D-IMEX problem requires space for storing the LU decomposed Jacobian matrix and can be expensive especially for the \nsr form.}
\begin{tabular}{ |l|cr|cr| }
\hline
Model & Memory (MiB) of 2D RTB  &  Ratio & Memory (MiB) of acoustic& Ratio\\
\hline
RK & 14 & 1.00 & 129 & 1.00\\
3D-IMEX \nsr GMRES & 406 & 29.00 & 3746 & 29.00  \\
3D-IMEX \sr GMRES & 69 & 4.90 & 639 & 4.90 \\
3D-IMEX \nsr BiCGstab & 30 & 2.14 & 280 & 2.17 \\
3D-IMEX \sr BiCGstab & 22 & 1.57 & 206 & 1.59 \\
1D-IMEX \nsr & 330 & 23.60 & 909 & 7.04\\
1D-IMEX \sr & 26 & 1.85 & 214 & 1.65\\
\hline
\end{tabular}
\label{memory}
\end{table*}

The 1D-IMEX time integration scheme has a potential to reduce computation time as compared to 3D-IMEX for three reasons: a) one needs to consider only one direction (radial direction on the sphere) when computing gradient and divergence, $\mathcal{O}(N)$ for 1D-IMEX instead of $\mathcal{O}(3N)$; b) 1D-IMEX does not require communication between  columns during parallel computation of the implicit solve stage; and  c) direct solver for the implicit solve stage has a constant cost per time step for any Courant number. This can lead to significant savings as compared to iterative solvers which typically require larger number of iterations with large Courant numbers, especially when used without preconditioning. The additional cost of 1D-IMEX comes from rotating vectors, such as velocity, to- and from- the spherical coordinate system.  We compare 1D- and 3D- IMEX running on a single GPU in Fig.\ \ref{comparison-1d3d} using a 2D-rising thermal bubble problem test case. To modify the Courant number in the vertical direction, we use two approaches. For the top plots in Fig.\ \ref{comparison-1d3d}, we refine the grid in the vertical without changing the time step. This increases the workload in each column while keeping the number of columns constant, hence, the Courant number in the horizontal direction $C_H$ is constant. A more realistic setup, which is in line with operational NWP on the sphere, is to keep the workload in a column the same while changing the number of columns by refining in the horizontal direction. The plots in the bottom half are produced using the second method. From the top plots, we see that 1D-IMEX gives some speedup over 3D-IMEX for the \sr forms but  the benefit for the \nsr forms is not significant despite the fact that the \nsr forms are slower. The plots in the bottom half show that the gain from 1D-IMEX is more or less the same for both the \sr and \nsr forms. In Fig.\ \ref{comparison-1d3d}, we also compare the use of a direct solver, i.e. LU decomposition of the Jacobian matrix followed by backward and forward substitution, against using GMRES iterative solver with no preconditioning. The 1D-IMEX run with the GMRES iterative solver is much faster than using a direct solver for the first setup in which we increase the workload in a column while keeping a constant $C_H$. The direct solver, conducted by constructing the Jacobian matrix and storing its LU decomposition, requires much more memory than using an iterative solver; as a result, we were not able to go above $C_V$=4 using this approach. As discussed before, this setup is rather unrealistic, so we look to the bottom plots for a more realistic evaluation of the direct solver. Here, the direct solver is much slower than the iterative solver for the \nsr forms at low $C_V \le 4$; for the \sr forms, the direct solver is faster than the iterative solvers even at low Courant numbers. This is because the \sr form has only one degree of freedom (pressure) which makes the matrix size  $\mathcal{O}(N^2)$, while it is $\mathcal{O}(25N^2)$ for the \nsr form which has 5 degrees of freedom. However, as the Courant number in the vertical direction is increased while keeping constant workload in each column, the direct solver becomes much faster than the iterative solvers. The reason is that with higher $C_V$ the iterative solvers require many iterations especially when used without preconditioning, while the work done by the direct solver is the same regardless of the Courant number. When the 3D-IMEX GMRES iterative solver is used with a preconditioner, we can see that the time to solution decreases in a similar way as the direct solver. The \nsr form requires a larger polynomial preconditioner (p=4) compared to the (p=1) used for the \sr form.
 
Finally, we look at the memory usage of the IMEX methods that often require significantly more memory than explicit methods, and thus lead to a reduction in the maximum problem size that can be run on the device. This can influence, for instance, strong scalability on cluster of  GPUs that have high communication latency between the CPU and GPU. In Table\ \ref{memory}, we show the memory usage of the different IMEX methods on two test cases: a 2D rising thermal bubble problem solved using 10x10 elements at polynomial order of 4, and a global simulation problem solved on a cubed sphere grid of 6x10x10x3 elements at polynomial order of 3. We can see that the \nsr forms require much more memory than the \sr forms in all cases. The 3D-IMEX \nsr forms using GMRES iterations require about 29X more memory than the explicit time stepping method. However, we can reduce the memory usage by using BiCGstab iterations instead of GMRES, because the former does not require additional space for storing intermediate Krylov vectors. With BiCGstab, the difference between the \nsr and \sr forms is not as significant. The 1D-IMEX \nsr form also requires much more memory than the \sr form but for a different reason. Because we use direct solver for 1D-IMEX, we need to store the LU decomposed Jacobian matrix. For eight prognostic variables (5 fundamental state variables + 3 dummy variables), the \nsr form requires about 64X larger memory than the \sr form to store the Jacobian matrix. Comparing the results of the two test cases, we can see that the normalized ratios of memory requirement are more or less the same, except for 1D-IMEX. This is due to the difference in the number of elements in a column -- 10 elements for the 2D rising thermal bubble problem and 3 elements for the global simulation problem. The more elements and higher polynomial order in the vertical direction, the larger the size of the Jacobian matrix.

\begin{figure*}

	\centering
	\begin{subfigure}[b]{\textwidth}
	\includegraphics[width=0.33\linewidth]{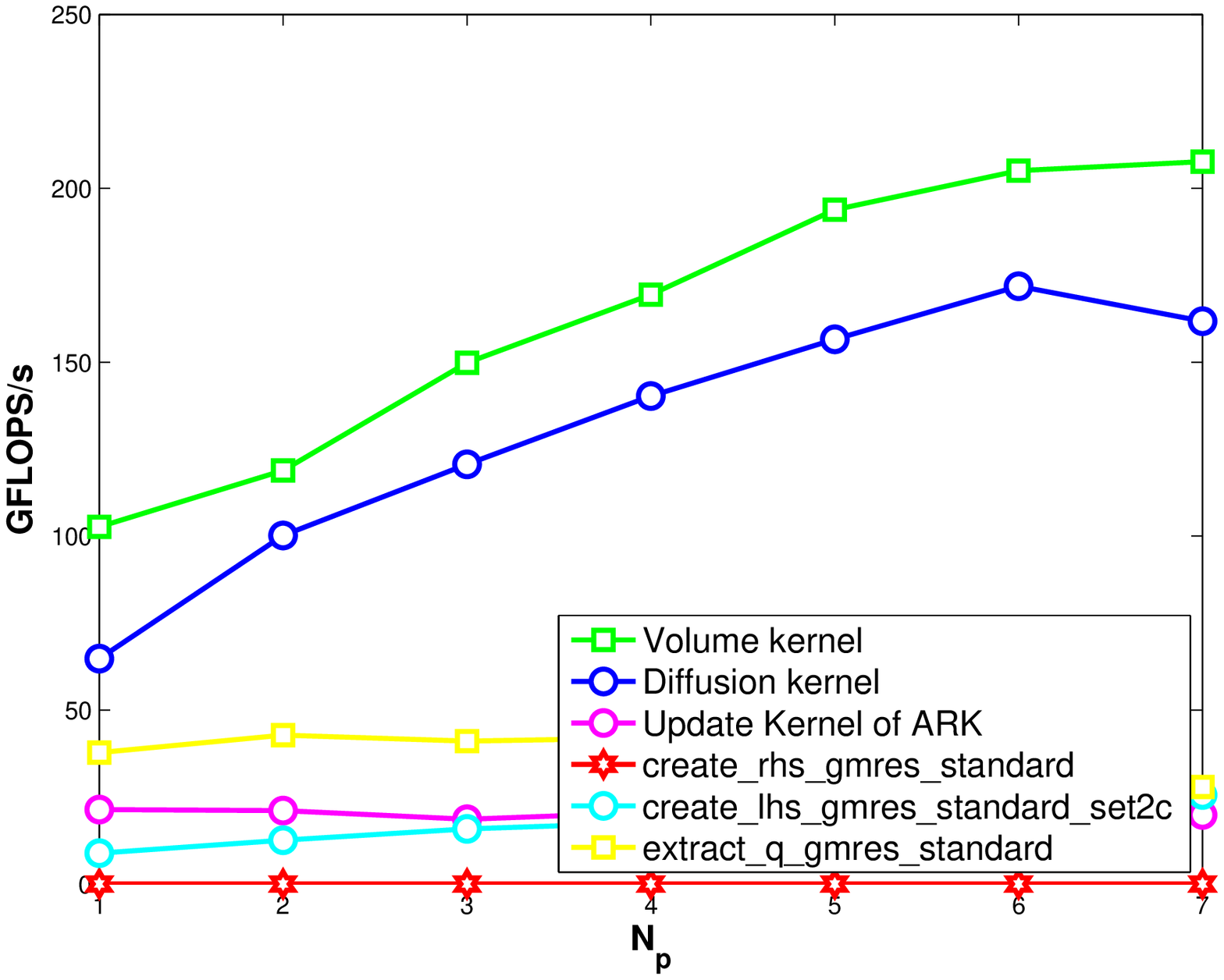}
	\includegraphics[width=0.33\linewidth]{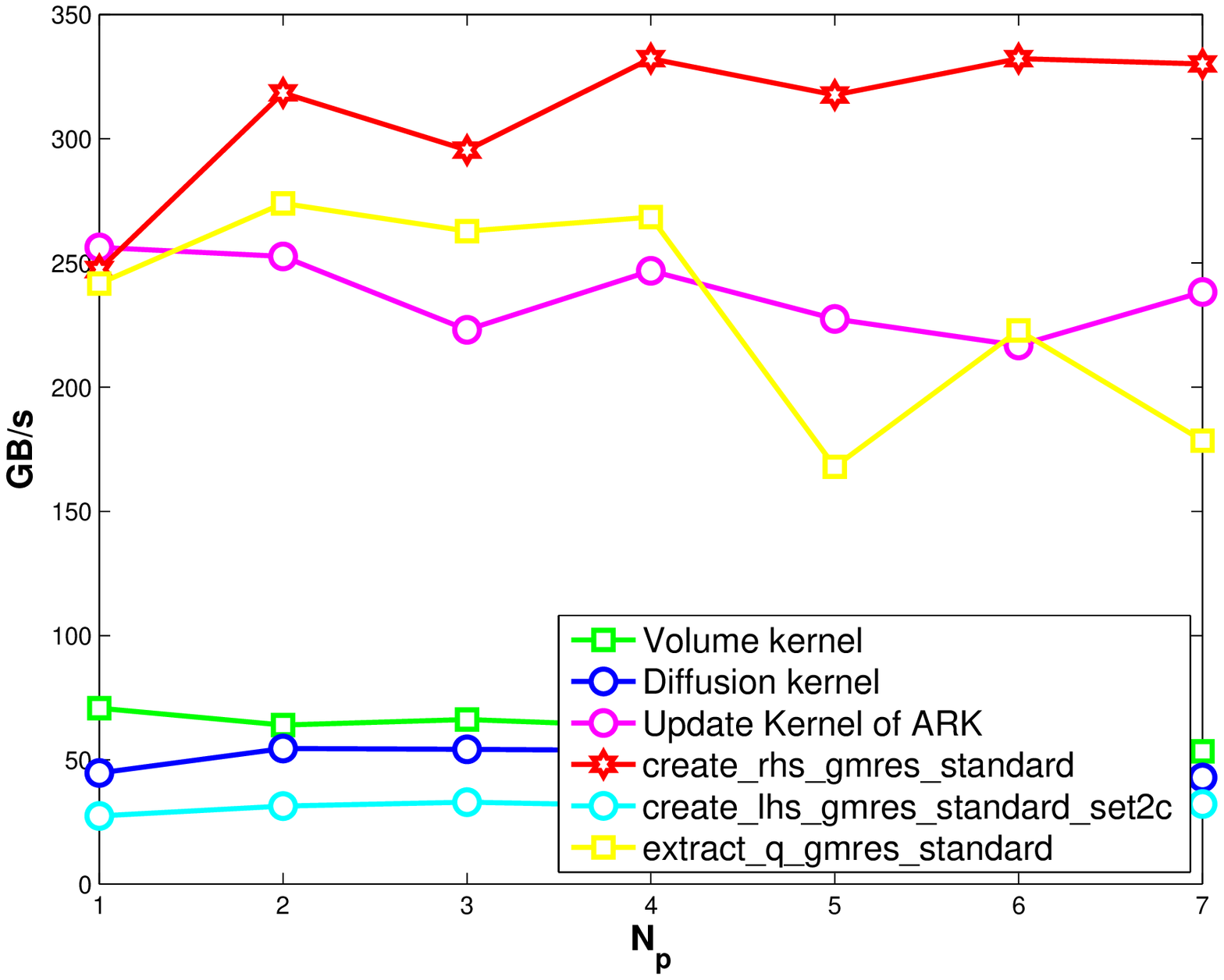}
	\includegraphics[width=0.33\linewidth]{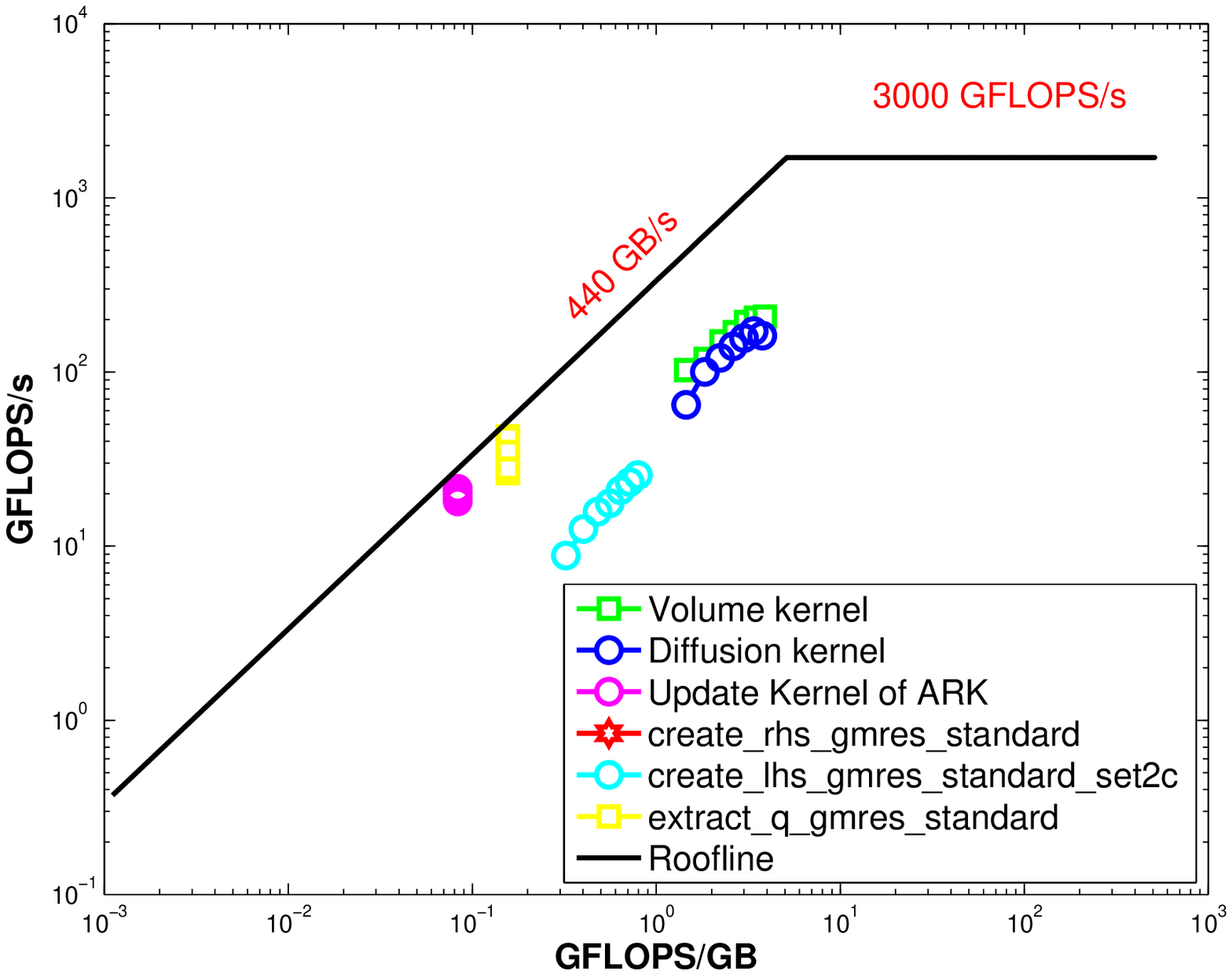}
	\caption{3D-IMEX in \nsr form kernels performance}
	\label{gflopsCGn-knl}
	\end{subfigure}
%
	\centering
	\begin{subfigure}[b]{\textwidth}
	\includegraphics[width=0.33\linewidth]{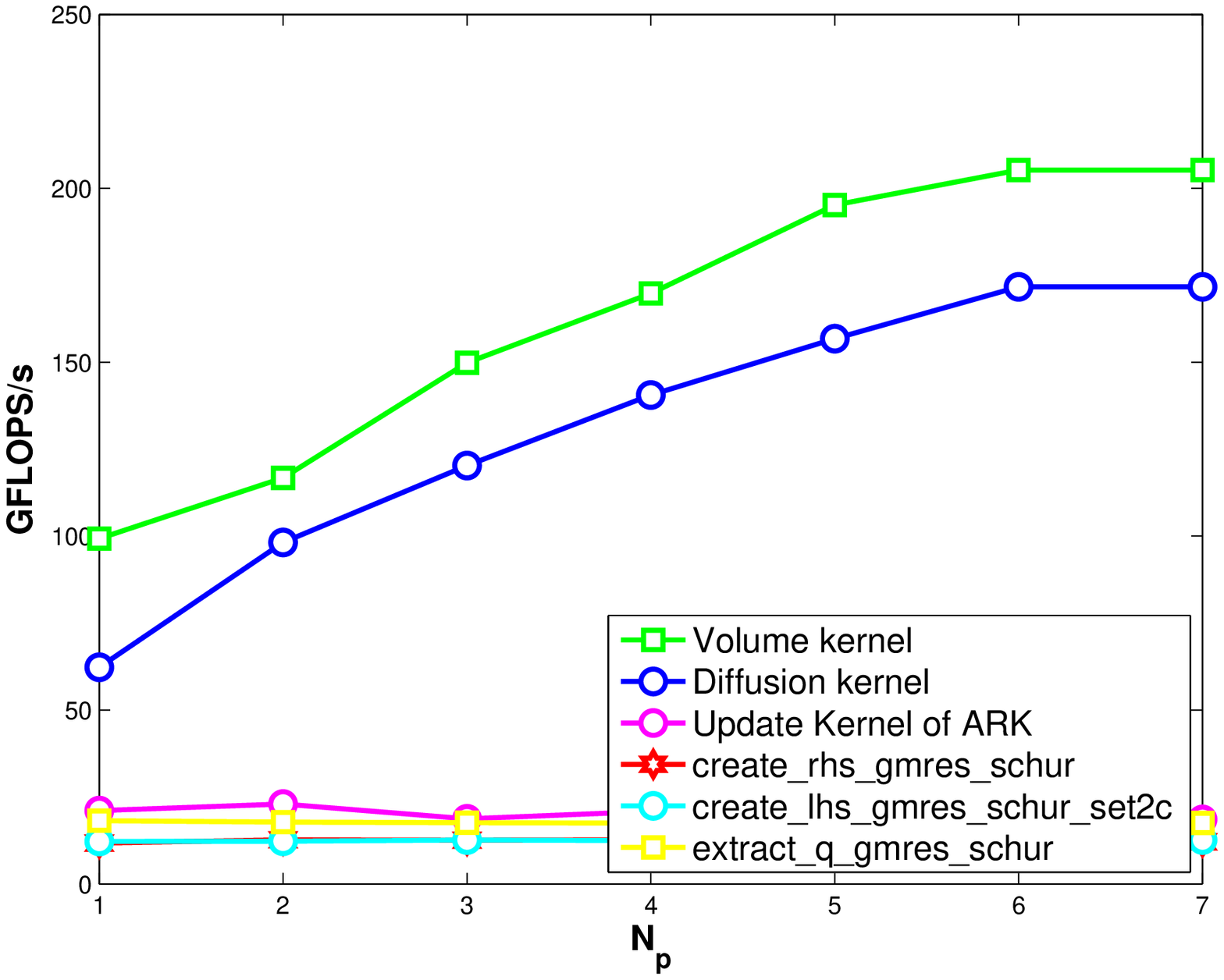}
	\includegraphics[width=0.33\linewidth]{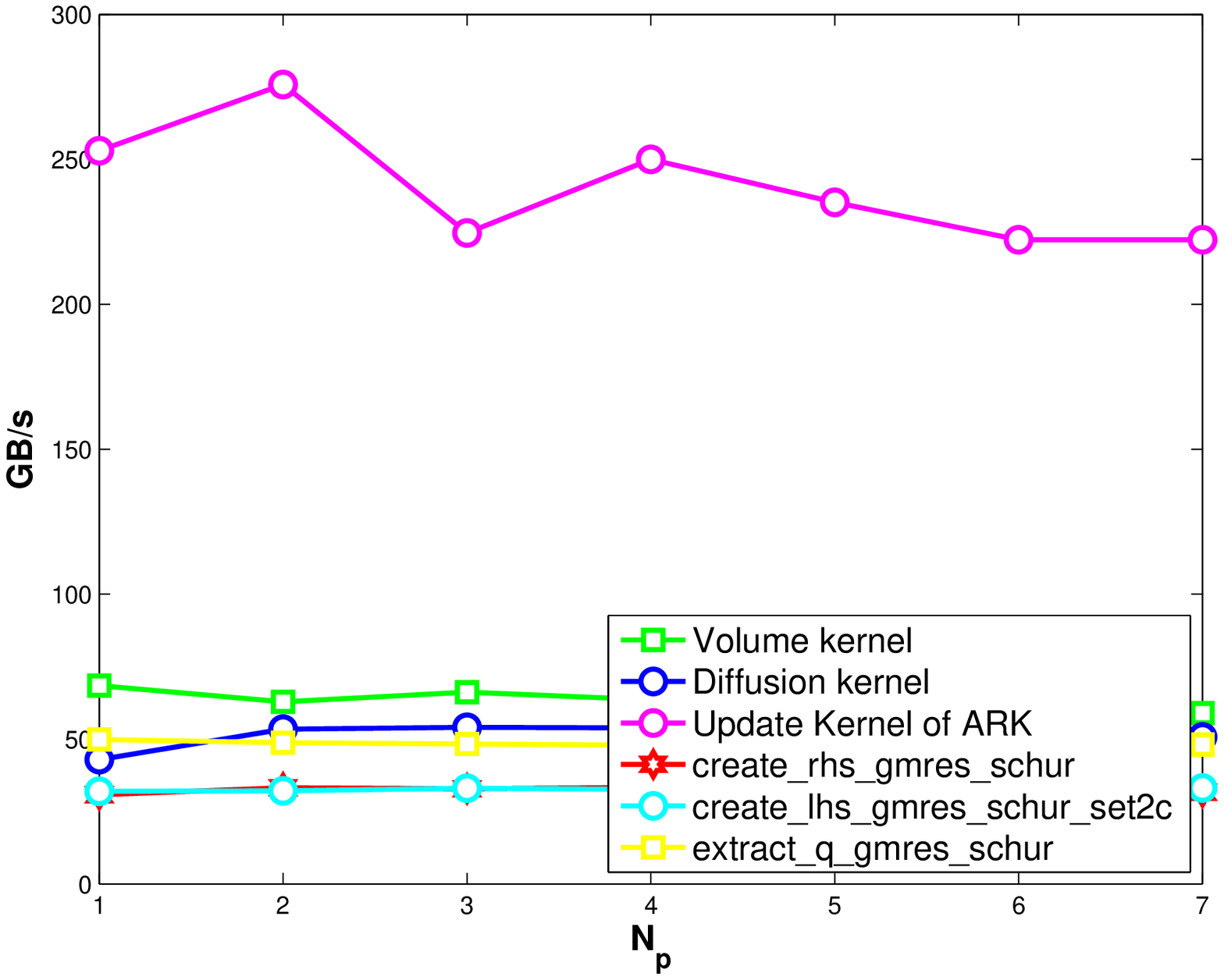}
	\includegraphics[width=0.33\linewidth]{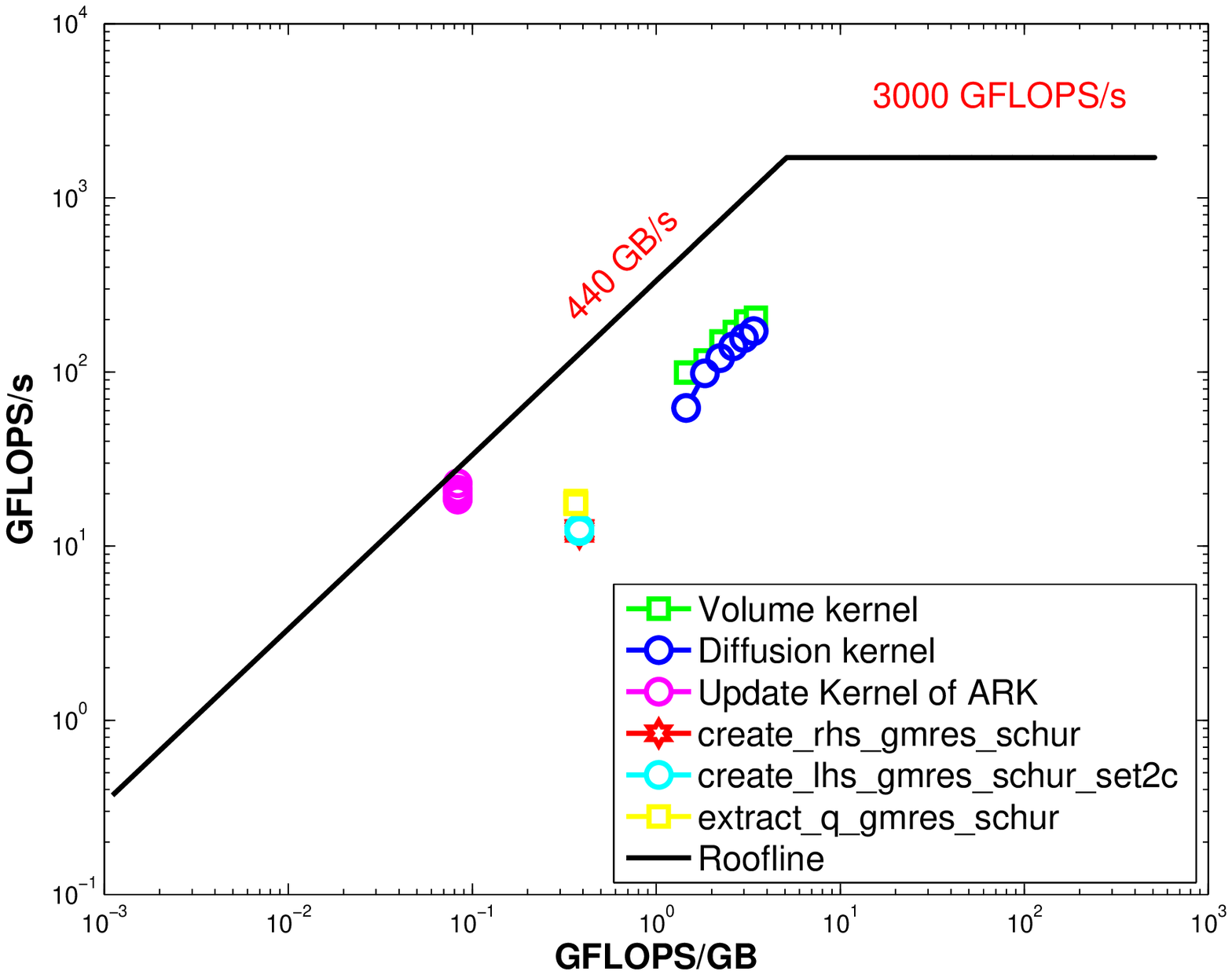}
	\caption{3D-IMEX in \sr form kernels performance}
	\label{gflopsCGns-knl}
	\end{subfigure}
	\caption{Performance of  kernels on the KNL are shown in terms of GFLOPS/s, GB/s and roofline plots to illustrate their efficiency. The test is done with MCDRAM usage mode set as a cache. }
%
	\centering
	\includegraphics[width=0.32\linewidth]{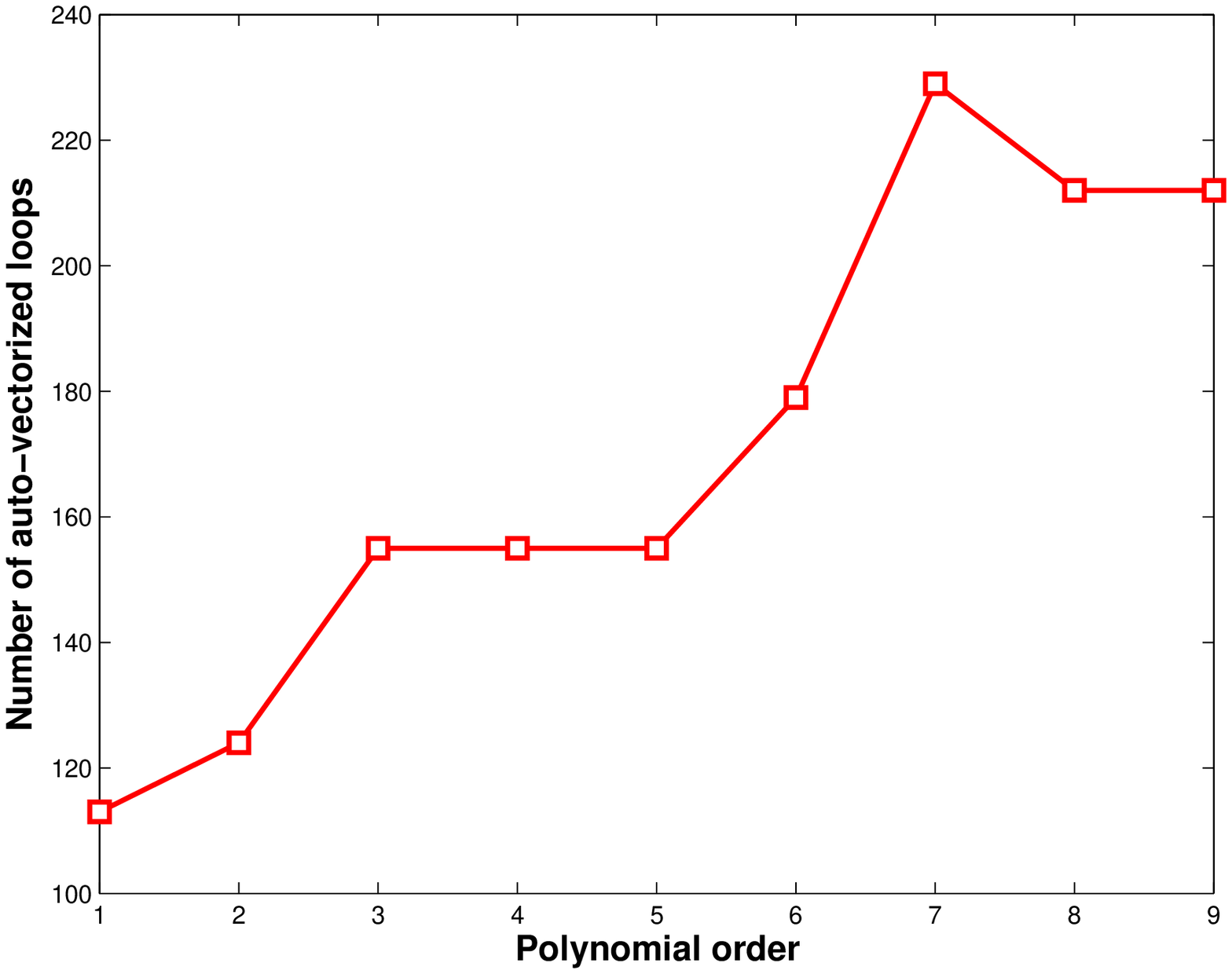} 
	\includegraphics[width=0.32\linewidth]{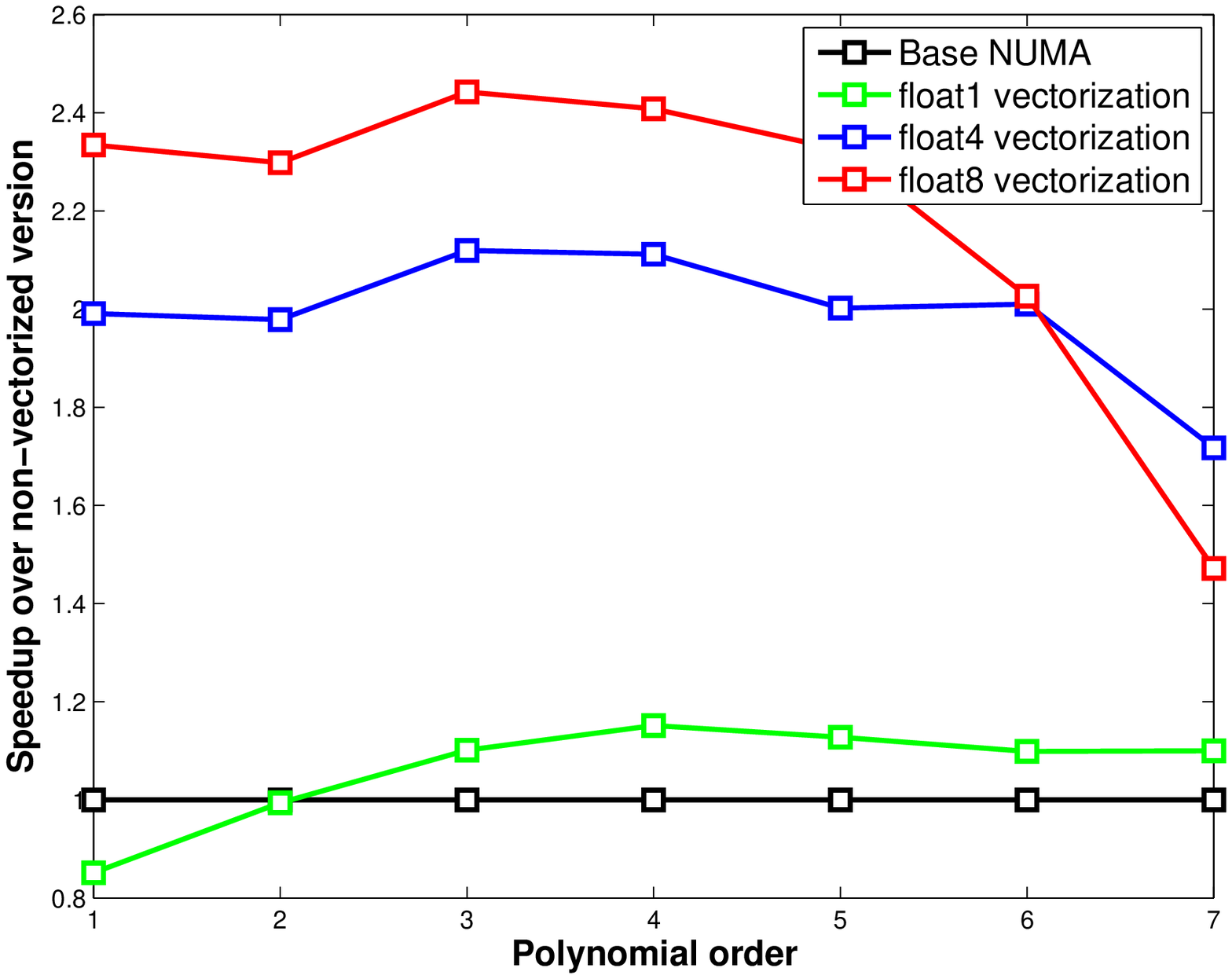} 
	\includegraphics[width=0.32\linewidth]{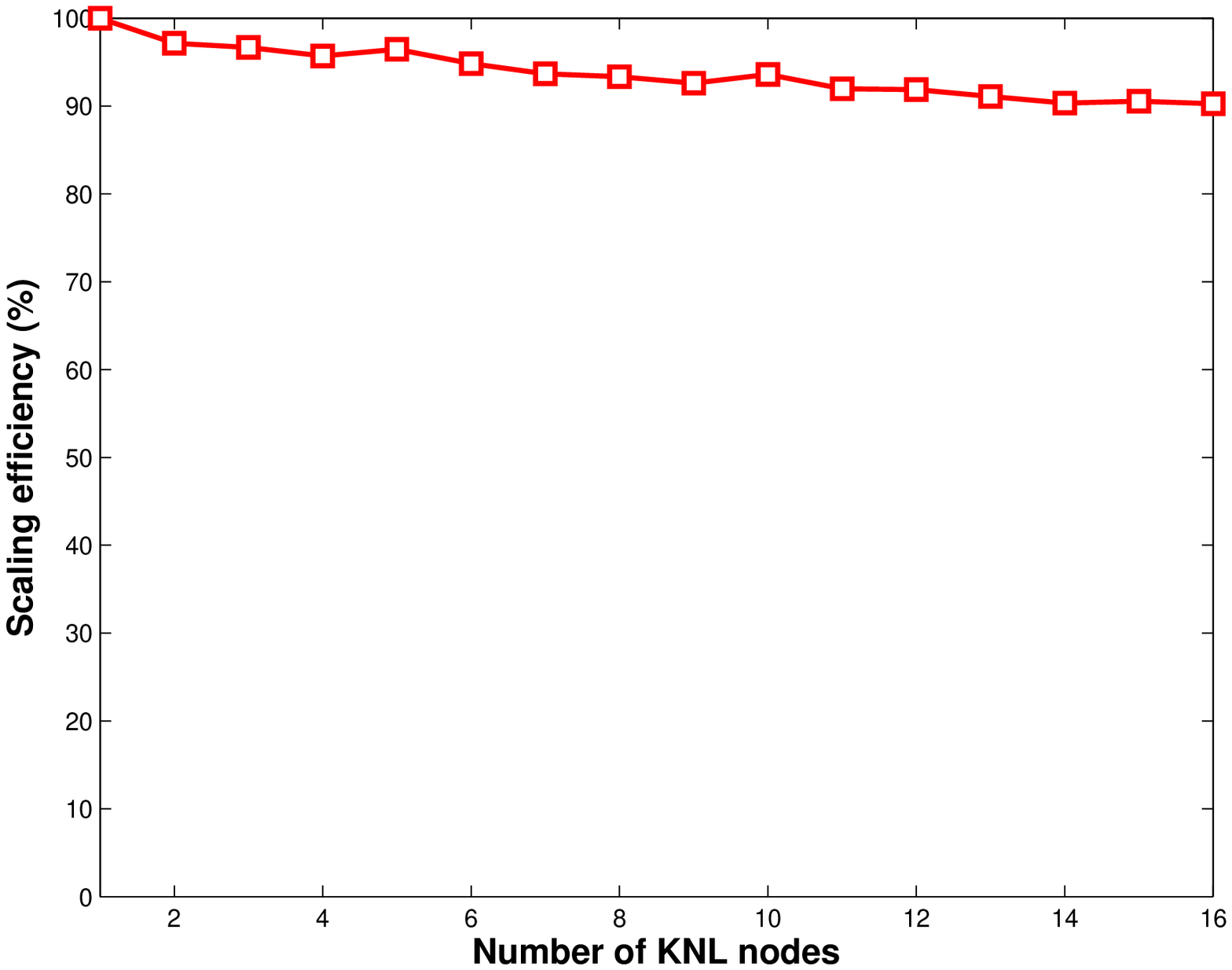} 
	\caption{(left) Auto-vectorization : the number of vectorized loops  is counted to produce this plot. (middle) Speedup of vectorized
	version is shown for float1/float4/float8 vectorization. (right) Strong scaling scalability using 16 KNL nodes in which each KNL node runs 64x2=128 tasks. The problem size for this test is about 32 million nodes.} 
	\label{vectorization}
\end{figure*}

 \begin{table*}
\centering
\small
\caption{Knights Landing (KNL) test: 64 physical cores with upto 4 hyper-threads per core. Time to solution of a 3D rising thermal bubble problem is given in seconds.
The optimal number of threads for NUMA is running 2 tasks (threads or processes) per core -- 128 tasks overall on a single KNL processor. No difference in performance is observed between threads and cores.}
\label{speedup}
\begin{tabular}{r|*{9}{r}{r}}
\hline
Processes & \multicolumn{8}{c}{Threads per process}  \\ 
& 1 & 2 & 4 & 8 & 16 & 32 & 64 & 128 & 256\\ \hline
128 & \trd{83.2} & 96.82 & 100 & & & & & &\\
64 & 101 & \trd{82.4} & 94.94  & & & & & &\\
32 & 184 & 101 & \trd{81.72} & 93.7& & & & &\\
16 &  &  & 101 & \trd{82.49} & 93.88 & 97 & & &\\
8 &  &  &  & 101 & \trd{82.3} & 94 & & &\\
4 &  &  & & &102 &\trd{84} & 98.17 & &\\
2 &  &  &  & & & 104 & \trd{86.44} & 102 &\\
1 &  &  & & & & & 107 & \trd{90.7} & 110\\ \hline
\end{tabular}
\end{table*}

\subsection{Performance on the KNL}
We conduct tests on the  Xeon Phi `co-processor' Knights Landing (KNL); the KNL is more like a standard CPU with many cores that allows standard threaded CPU programs to run without code modifications. 
The KNL has 64 or more cores that are able to execute 4 threads per core simultaneously, and also a new integrated on-package memory that delivers upto 5X the bandwidth of standard DRAM memory.
Because we used OCCA for programming both the CPU and the GPU, we did not need to rewrite the GPU kernels for KNL. However, we added vectorization support for the KNL that was not fundamentally necessary for the GPU. Moreover, the way OCCA kernels are run on the CPU is through what is known as an element-per-node approach where one or more elements are processed with a single OpenMP thread, while a node-per-thread approach is used on the GPU. This design gives best performance on both the CPU and GPU.

In Figs.\ \ref{gflopsCGn-knl} - \ref{gflopsCGns-knl}, we show the GFLOPS/s and GB/s of our kernels on the KNL using the OpenMP translation of the OCCA kernels.  Due to the lack of tools for evaluating individual kernel performance on the KNL, we use the `hand-counting' approach for measuring the number of FLOPS performed and bytes transferred per second. As expected, the volume and diffusion kernels show the highest rate of floating point operations  at about 210 GFLOPS/s, which is about 7\% of peak performance in double precision (3 TFLOPS).  The ARK time update kernel shows maximum bandwidth usage along with other similar kernels such as the \nsr form right hand side evaluation kernels that are bottlenecked by the data transfer rate. The bandwidth usage differs significantly depending on the type of RAM used for storing the data, for which two types are available on the KNL. The high-bandwidth RAM (440 GB/s MCDRAM) of the KNL is about 5X faster than the standard DRAM4 (90 GB/s peak). The results shown here is with all the 16GB MCDRAM used in cache mode. This means we did not need to modify our code, however, better results could potentially be obtained by managing memory ourselves. We can see that the time update kernel, the \nsr form right-hand side and result `extraction' kernels  get close to peak bandwidth usage using the cache mode while most other kernels are far from reaching the peak.  Looking at the roofline plots, most kernels are memory-bound with the volume and diffusion kernels close to being compute bound.

One can run on the KNL with a standard MPI program without OpenMP; however, a hybrid MPI-OpenMP program could perform better as threads could be more `lightweight' than processes and also because
communication between threads is faster through shared memory. The results in Table \ref{speedup} show that there is not much difference between using threads or processes on the Knights Landing (KNL) -- which was somewhat a surprise for us as we hoped to gain some benefit from adding OpenMP support.  The least amount of time to solution is obtained on the diagonal of the table where 128 tasks are launched per KNL device. Deviating from this optimum by using either 1 thread per core or 3 threads per core results in a significant increase in time-to-solution. Focusing on the diagonal of the table, we can see that the largest time-to-solution obtained using 128 MPI processes (90.7 sec) is  only about 7.7\% slower than that obtained using pure OpenMP with 128 threads (83.2 sec).

Performance improvement from vectorization is shown in Fig.~\ref{vectorization}. This is achieved through:  a) automatic vectorization of inner loops of ($i,j,k$) over the LGL points and b) using
float4/float8 vectorization over the field variables ($U,V,W,p,\rho,\Theta,-,-$). The first approach yields the best vectorization result at polynomial order 7, i.e., where number of nodes = $8^3$ in an element. 
The volume kernels are re-written to compute gradients and divergence of the 8 field variables all at once for the float8 approach using all the 512-bit wide SIMD units of the KNL. Therefore, we can potentially obtain  a speedup of 8X from float8 vectorization, however, this can only be achieved if the program is compute-bound. None of the kernels we showed so far are compute-bound so we do not expect this kind of result. Indeed, we can see in Fig.~\ref{vectorization} that the float4/float8 vectorization both show more than 2X speedup improvement; however, the float8 vectorization, that uses all the SIMD units, does not show a significant improvement over the float4 vectorization confirming the memory-bound nature of our program.

In Fig.~\ref{vectorization}, we show strong scalability results -- which is on the upper side of 90\% -- using 16 KNL nodes. The problem size used for this scalability test does not fill up the memory available on one KNL node let alone all 16 nodes, hence, we expect better scalability results with bigger problem sizes for both strong and weak scalability tests. 


\section{Acoustic wave problem}
For the purpose of validating the IMEX time-integrators, we consider the benchmark problem given in \cite{tomita2004}, namely,  an acoustic wave traveling around the globe. Though the wave travels at the speed of sound, the horizontal grid resolution is such that $C_H \le 1$. On the other hand,  the vertical $C_V$ can exceed 1 thereby leading to a potential performance gain using an IMEX  time integration method. The initial state  is hydrostatically balanced with an isothermal background potential temperature of $\theta_0$=300 K. A perturbation pressure $P'$ is superimposed on the reference pressure
\[
P'=f(\lambda,\phi) g(r)
\]
where
\[
f(\lambda,\phi)  = 
\begin{cases}
        0  &  \text{for } r > r_c \\
        \frac{\Delta P}{2} (1 + \cos(\frac{\pi r}{r_c})) &  \text{for } r \le r_c
\end{cases}
\]
and
\[
g(r)=\sin\left(\frac{n_v \pi r}{r_T}\right)
\]
where $\Delta P$ = 100 Pa, $n_v=1$, $r_c=r_e/3$ is one third of the radius of the earth $r_e=6371$ km, and the top of the model is $r_T$=10 km.  The geodesic distance $r$ is calculated as
\[
r = r_e \cos^{-1}[\sin\phi_0 \sin\phi + \cos\phi_0\cos\phi\cos(\lambda - \lambda_0)]
\]
where $(\lambda_0,\phi_0)$ is the origin of the acoustic wave. 

\begin{figure*}
        \centering
	\includegraphics[width=0.45\linewidth]{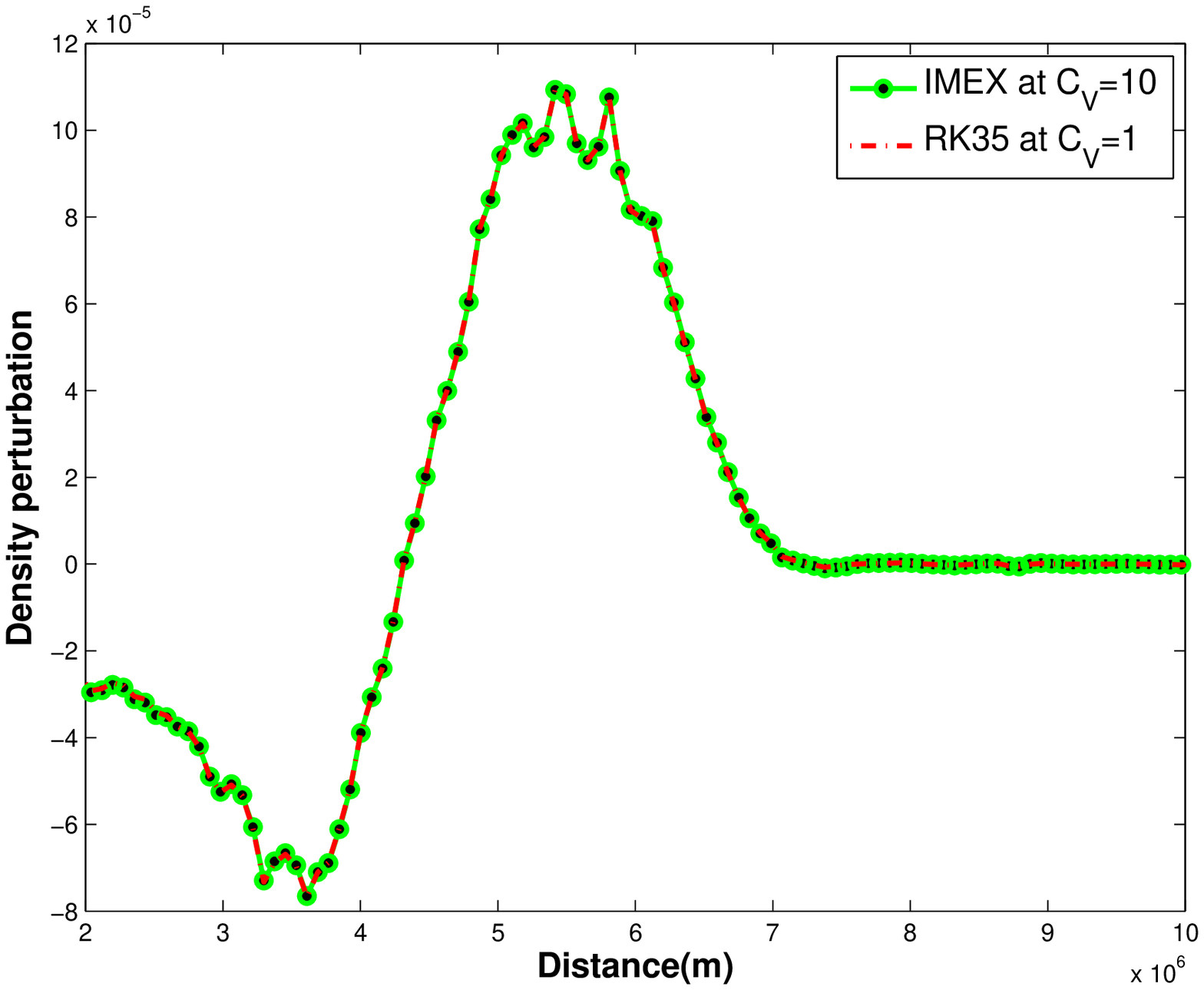} 
	\includegraphics[width=0.45\linewidth]{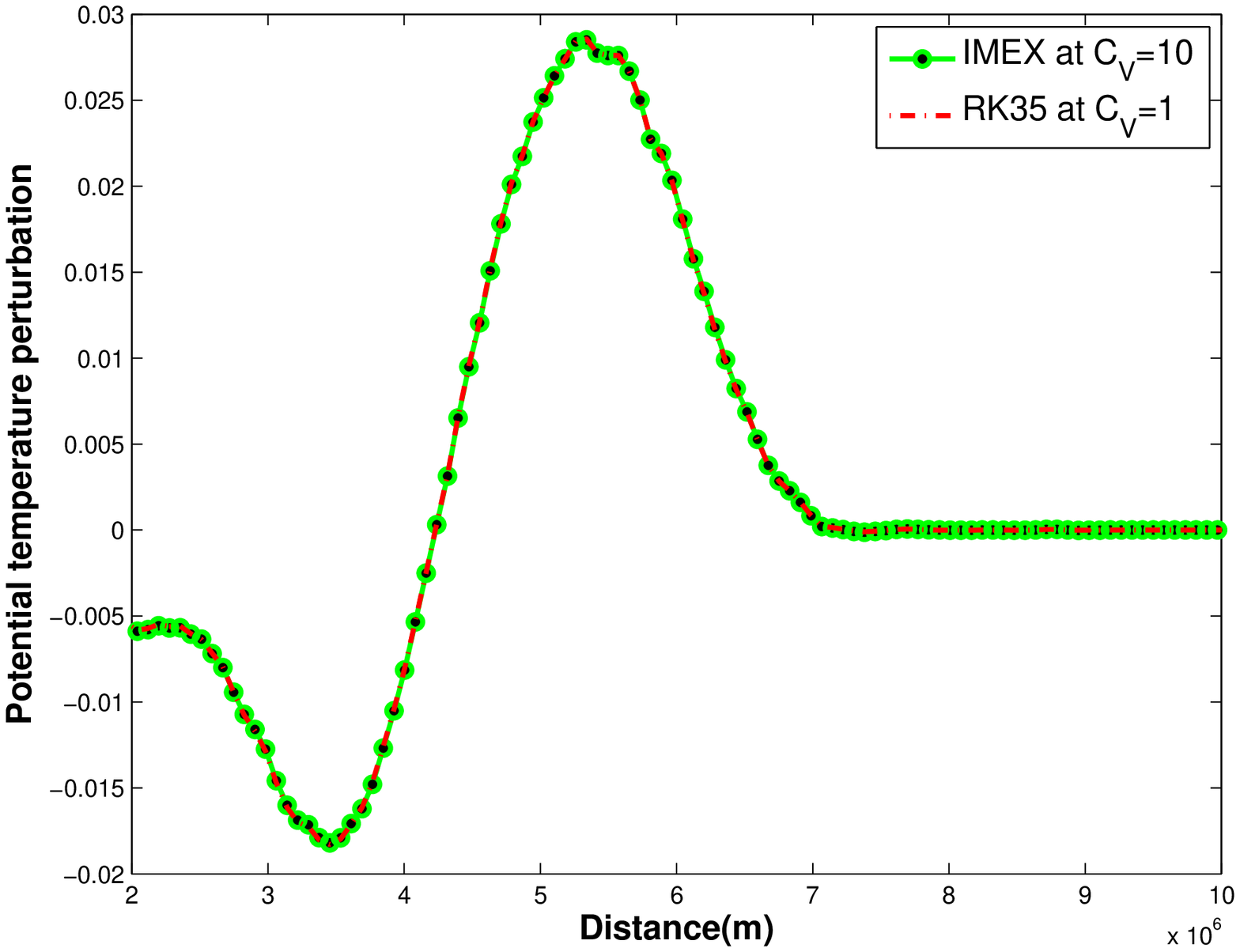} 
	\includegraphics[width=0.45\linewidth]{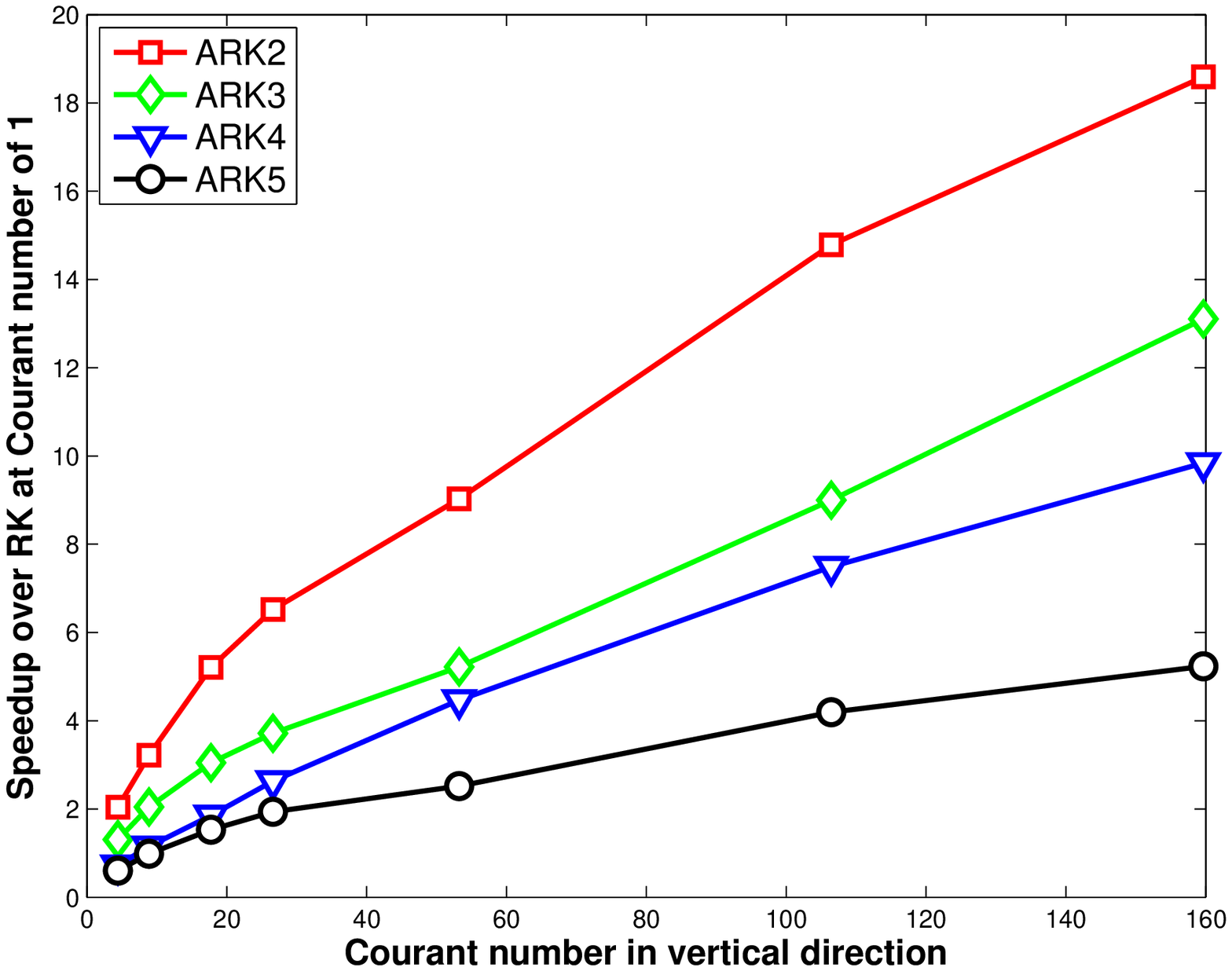} 
	\includegraphics[width=0.45\linewidth]{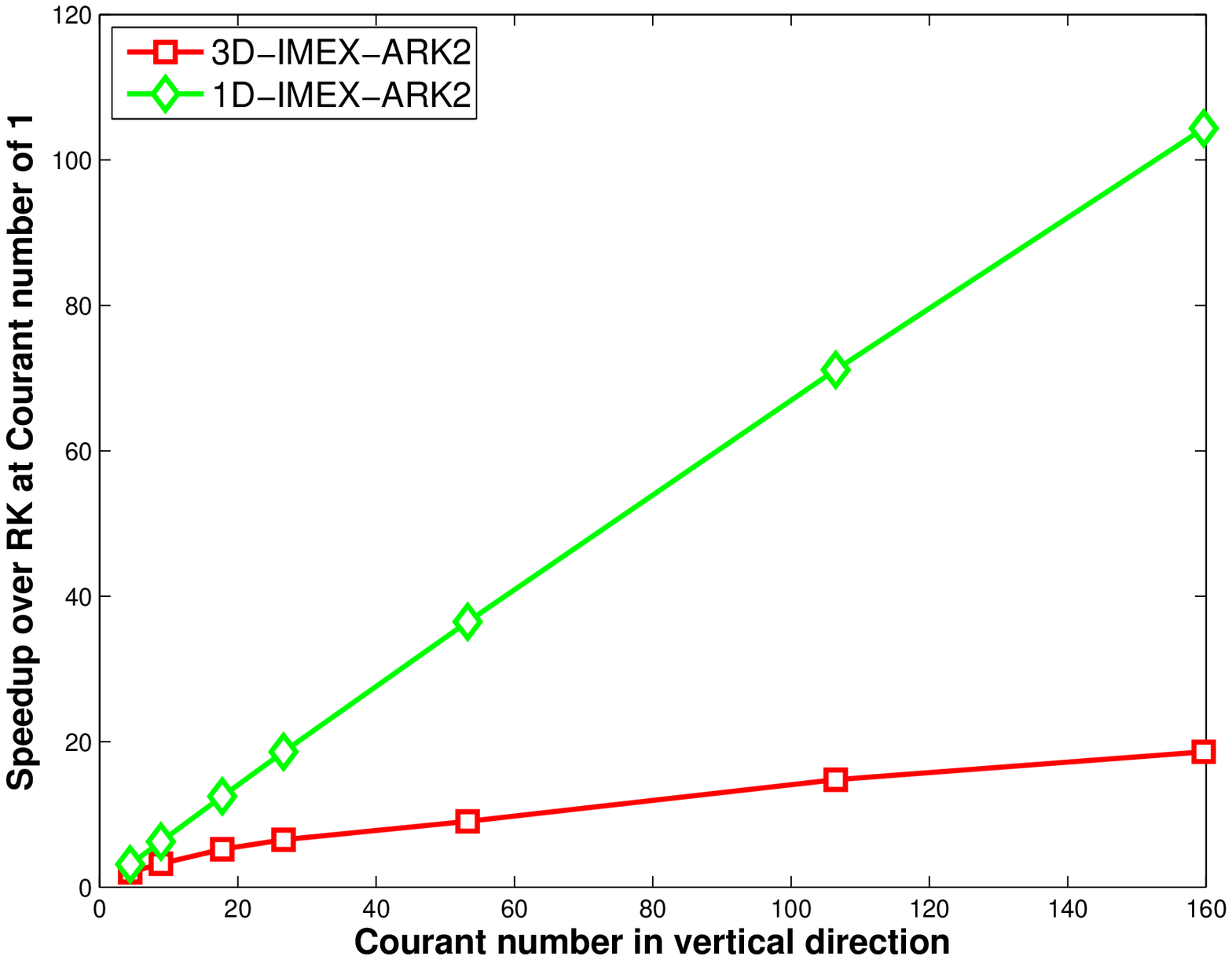} 
	\caption{The acoustic wave problem of \cite{tomita2004}:  Density perturbation (top left) and potential temperature perturbation (top right) results after 4 hours match exactly with each other, even though IMEX took 3X less time to complete. The bottom left figure shows the relative speedup of IMEX over RK at different Courant numbers. At $C_V=150$, a relative speedup of about 20X is observed using the second order ARK2 scheme. If we, instead, use the 1D-IMEX \sr form with a direct solver, we get a speedup of upto 100X.} 
	\label{acoustic}
\end{figure*}

The grid is a cubed sphere as such $6 \times 10 \times 10 \times3$ (6 panels, each panel has 10x10 elements on the surface with 3 elements in the vertical) for a total of 1800 elements with $5^{th}$ order polynomials.  No-flux boundary conditions are applied at the bottom and top surfaces. The explicit time stepping is conducted with a time step $\Delta t = 1$ sec so that $C_V ~= 1$ while $C_H$ is well below 1. The IMEX time integration uses a time step of $\Delta t=10$ sec so that $C_V~=10$. Even though IMEX starts out with a 10 time larger time step, we should not expect to recover all of the speedup because the implicit solve in each time step makes the IMEX time update much slower than explicit time updates of the same order. The 3D-IMEX \sr form with GMRES iterations using $1^{st}$ order PBNO preconditioner is used for this test. The explicit method took about 2356 sec to complete a 36 hours simulation while the IMEX method took about 796 sec, which is a relative speedup of about 3X in favor of IMEX (yielding an approximate gain of 33\% of the time step size increase).  However, a much larger speedup of upto 20X can be obtained for this particular problem by increasing the  Courant number in the vertical direction to 150. The increase in speedup with $C_V$ is shown in Fig.\ \ref{acoustic}.  If we, instead, use the 1D-IMEX \sr form with a direct solver, we get a speedup of upto 100X. The time-to-solution of the direct solver increases linearly with the time step, because the workload does not increase with the Courant number. This is unlike iterative solvers which require more iterations to converge with increasing Courant number. The polynomial order for the PBNO was increased from 1 to 3 for the range of Courant numbers required for the plot, hence, the preconditioning was not enough to bring the performance of the iterative solver close to that of the direct solver.

Higher order ARK methods can be used to improve the accuracy of the results especially when a large Courant number (larger time step) is used. We can see in the bottom left plot of Fig.\ \ref{acoustic} that the higher order ARK3, ARK4 and ARK5 still give significant speedups over using an explicit RK method. To compare the relative accuracy, we plot the density and temperature distributions obtained using $C_V=10$ for 3D-IMEX against those obtained using explicit time integration. We can see that both the density and potential temperature perturbations match exactly with each other after 4 hours into the simulation. The speed of sound after 16 hours, i.e., at the time the wave reaches the antipode, is calculated to be about 348 m/s. This result is in close agreement with the sound speed calculated from the initial conditions using the relation $a = \sqrt{\gamma p / \rho}=$ 347.32 m/s (less than 0.2\% relative error).

\section{Scalability of IMEX}

\begin{figure*}
        \centering
	\includegraphics[width=0.6\linewidth]{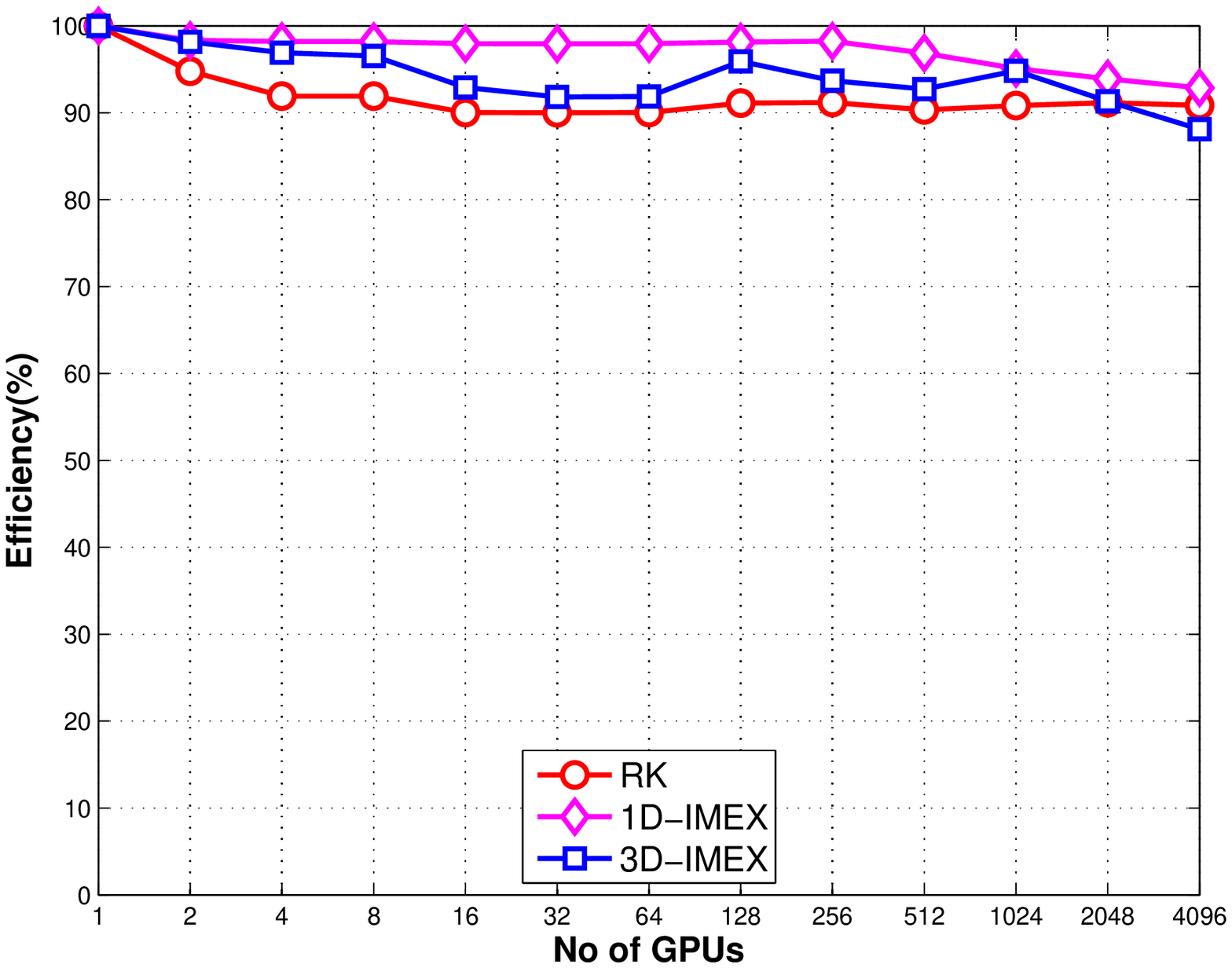} 
	\caption{Weak-scalability of RK, 1D-IMEX and 3D-IMEX on multiple GPUs. The test case is the 2D rising thermal bubble problem in which each GPU gets a subdomain of $30 \times 10$ elements with polynomial order of 6. The number of elements in the vertical direction is kept constant at 10. The same time step and grid is used for all simulations, therefore the IMEX simulations are slower than that of RK, which contributes towards a better scalability for IMEX. We can also conclude that 1D-IMEX has a better scalability than 3D-IMEX due to the absence of communication during the implicit solve stages.} 
	\label{scalability}
	
        \centering
	\includegraphics[width=0.8\linewidth]{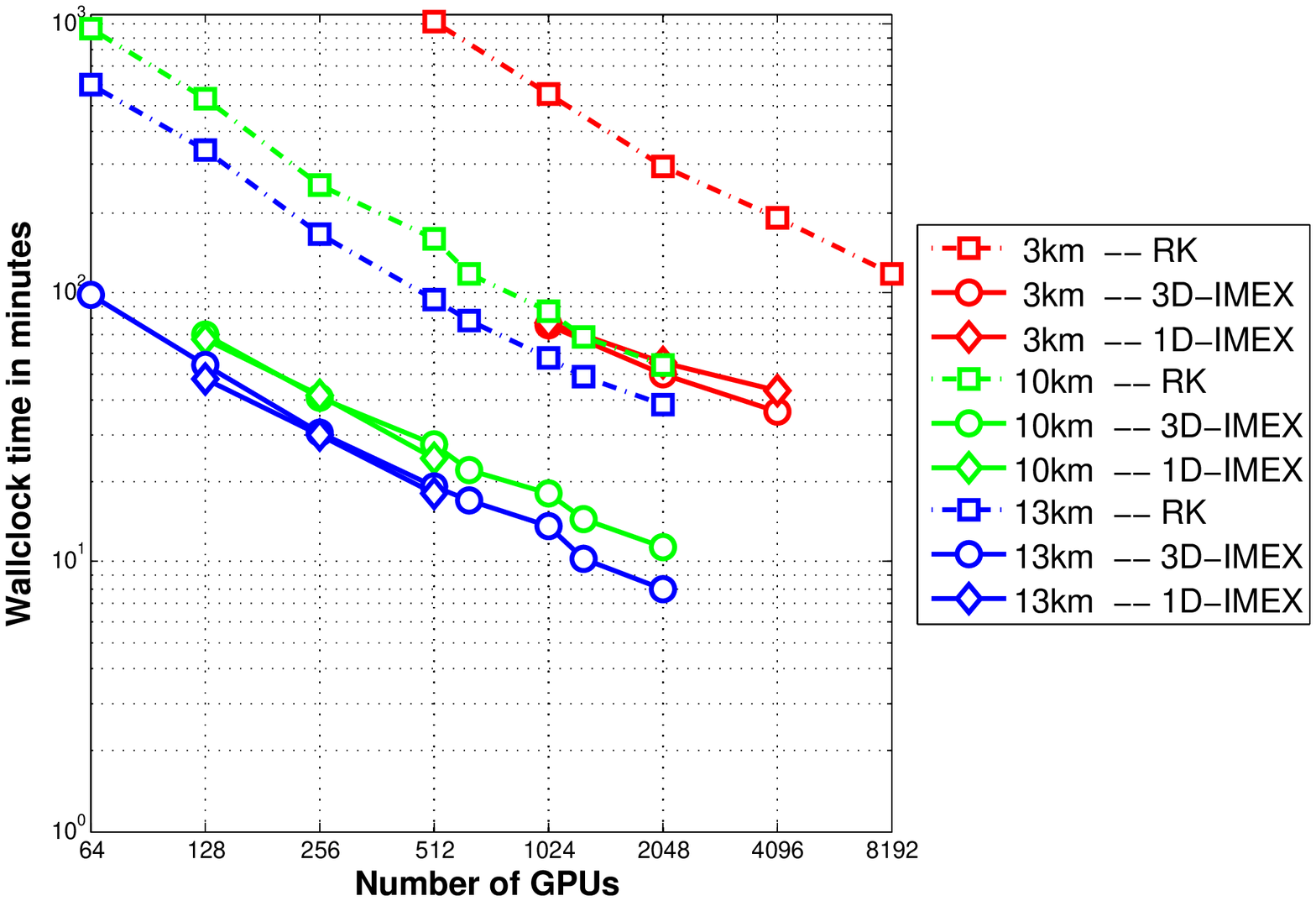} 
	\caption{Strong-scalability of 3D-IMEX on multiple GPUs is tested on the global acoustic wave simulation problem using grid resolutions of 13 km (6x112x112x4 elements at N=7), 10 km (6x144x144x4 elements at N=7) and 3 km (6x448x448x4 elements at N=7). 3D- and 1D- IMEX show  constant speedups of about 5X over the explicit method at all grid resolutions.} 
	\label{scalability2}
\end{figure*}

We show weak scalability results of the explicit and IMEX methods in Fig.\ \ref{scalability}. The 2D rising thermal bubble problem, instead of a global simulation problem, is used for this  test because of the convenience it offers in increasing problem size as required by weak scalability tests. We fix the problem size on each GPU to be a subdomain of  $30 \times 10$ elements with polynomial order of 6. The number of elements in the vertical direction is kept constant at 10, whereas we double the number of columns in the horizontal for each doubling in number of GPUs. The same grid and time step is used for the explicit and IMEX methods. This makes the IMEX methods to be  slower than the explicit method because the Courant number is kept below one for all methods to work properly. Hence, the fact that the IMEX method is slower in this test is not a reflection of its performance at a larger Courant number; however, IMEX could see better scalability results since it is slower. For 3D-IMEX, we fix the number of BiCGstab iterations per time step so that  the scalability results are not spoiled by non-constant number of iterations. This setup allows us to measure the relative cost of 1D- and 3D- IMEX for one pass of the algorithms; however, it does not provide a practically useful scalability result -- which we consider later in the strong scalability test. The average wallclock times are 50 sec, 110 sec and 27 sec for the explicit RK, 1D-IMEX and 3D-IMEX, respectively. The main thing to take from this plot is that 1D-IMEX scales better than 3D-IMEX because the latter requires additional communication during the implicit solve stages. Moreover, we can conclude that IMEX scales as well as the explicit methods despite the difference in wallclock times.

The strong scalability of IMEX on a multi-GPU cluster is evaluated using larger grid resolutions of 13 km, 10 km and 3km achieved using 6x112x112x4, 6x144x144x4 and 6x448x448x4 elements, respectively, at  polynomial order of 7. We can see in the scalability plot in Fig.\ \ref{scalability} that both IMEX methods show scalability that is as good as the explicit method for all  grid resolutions, while achieving about a 5X speedup. IMEX achieves  target wall clock minutes at fewer number of GPUs than the explicit time integration method. For example, on the coarse grids, IMEX requires about 64 GPUs to achieve a target 100 min wall clock time while RK requires 512 GPUs for the same target; IMEX brings down the wall clock minutes below 10 min if we use the same number of GPUs for both, i.e., 512. Therefore, we can conclude that IMEX plays an important role in decreasing both the wallclock minutes on a single-node machine and also the number of nodes required to achieve a specified target wall clock minutes on multi-node clusters. 

\section{Conclusions}
This work presented the porting of  IMEX time-integrators to manycore processors for the solution of the governing equations in NWP, namely, the Euler equations. Several improvements were required to better the performance of explicit time stepping methods. The first improvement comes from rewriting the equations in \sr form to reduce the number of degrees of freedom and also the condition number of the resulting linear system of equations (as well as positioning the eigenspectra in a friendlier region for the Krylov methods). As a result, the iterative solution of the \sr form is about two times faster than the \nsr form. The second improvement comes from formulating a 1D-IMEX method that is in better agreement with the horizontal/vertical length and time scales of real global atmospheric simulations than 3D-IMEX. In 1D-IMEX, the vertical motion is treated implicitly via columns that are solved independently of one another, while the horizontal motion is treated explicitly. The ARK schemes are used for the IMEX methods in which we conduct an explicit update followed by a corrective stage that involves implicit solves. The rest of the improvements come from accelerating the solution of the implicit solve stage. The resulting system of equations from 3D-IMEX are often solved iteratively which often requires preconditioning. A simple but effective PBNO preconditioner is used to accelerate the solution of both the \nsr and \sr forms. Furthermore, 1D-IMEX  can  benefit from direct solution of the small matrices in a column to accelerate the solution significantly, especially at high Courant numbers where iterative solvers struggle without preconditioning.

We first presented the IMEX formulation of the Euler equations conducted by removing the fast waves from the explicit update stage for two equation sets, Set2NC that is in non-conservation form, and Set2C that is in conservation form. An exact linear operator for the implicit terms is constructed using background reference states for the prognostic variables and then solved for the perturbation components. The implementation of the spatial discretization of the linear operator for cG is described for both the \nsr (first order) form and \sr (second order) forms; the \sr form for dG is currently not available. Several iterative solvers (GMRES, BiCGstab, etc.) for 3D-IMEX and a direct solver for batched solutions of the small matrices  in 1D-IMEX is implemented on manycore processors using the unified OCCA language. Performance of the IMEX kernels is evaluated using  rate of floating point operations (GFLOPS/s) and rate of data transfer  (GB/s) as metrics. Roofline plots are also provided for a K20X GPU and Intel's KNL hardware. We conducted several parametric studies to evaluate the performance improvements required for IMEX to beat an explicit RK scheme. On the GPU the 3D-IMEX method yielded an average performance improvement of about 4X over the explicit scheme at Courant numbers of about 15; while 1D-IMEX provides even more speedup when using a direct solver for the implicit stage. Using OpenMP (CPU) mode on the KNL, the performance of IMEX reaches about 7\% of peak for the volume kernel. To get this performance on the KNL, the kernel codes needed additional optimizations a) vectorization to process upto eight double precision calculations, b) taking advantage of fast MCDRAM.

For the purpose of validation on a realistic NWP problem, we considered the case of an acoustic wave problem traveling around the globe. The results from simulations conducted using IMEX match those obtained with explicit time stepping methods. The most important result we are interested in is that IMEX achieves a relative speedup of about 3X over the explicit method without degrading the quality of the results at $C_V=10$. Larger Courant numbers are encountered in operational NWP and this can yield significantly larger speedups. To demonstrate this point, we showed a speedup of upto 100X using the 1D-IMEX scheme  at $C_V=150$. The coarser the horizontal grid resolution is, the more speedup can be obtained using 1D-IMEX. One should not expect any speedup from 1D-IMEX if the horizontal and vertical resolutions are equal. 

Finally, we performed both weak and strong scalability tests on the Titan supercomputer which is a cluster of GPUs. The weak scalability test showed that the IMEX methods scale as well as the explicit methods -- with the 1D-IMEX showing better scalability than 3D-IMEX due to the absence of inter-processor communication during the implicit solve stage. The strong scalability of IMEX was also as good as the explicit methods using upto 4096 GPUs, while yielding a relative speedup of 5X over the explcit methods.
The strong scalability of the IMEX methods is linear on a log-log plot for upto 4096 GPUs. We expect this scalability to continue as long as there is enough work per GPU card. As regards the scalability on Intel's KNL, we measured strong scalability of IMEX to be about 90\% using upto 16 KNL nodes. 

\section{Acknowledgement}
This research used resources of the Oak Ridge Leadership Computing Facility at the Oak Ridge National Laboratory, which is supported by the Office of Science of the U.S. Department of Energy under Contract No. DE-AC05-00OR22725. The authors gratefully acknowledge support from the Office of Naval Research through PE-0602435N. We would also like to thank Intel for letting us conduct tests on the Knights Landing cluster Endeavor, and AFOSR Comp Math.

\bibliographystyle{sageh}
\bibliography{Giraldo_refs,references}

\end{document}